\documentclass[11pt]{article}
\usepackage{graphicx}
\usepackage{amssymb,amsfonts,amsthm}

\usepackage{amsmath,amsthm,amssymb}
\usepackage{amscd}
\usepackage{amsfonts}
\usepackage{color}
\newcommand{\Sk}{{\mathrm{Sk}}}

\newcommand{\la}{\langle}
\newcommand{\ra}{\rangle}

\newcommand{\vk}{{van Kampen}\ }
\newcommand{\dv}{{\mathrm{Div}}}

\newtheorem{theorem}{Theorem}[section]
\newtheorem{lemma}[theorem]{Lemma}
\newtheorem{notation}[theorem]{Notation}
\newtheorem{problem}[theorem]{Problem}

\newtheorem{cor}[theorem]{Corollary}
\newtheorem{proposition}[theorem]{Proposition}
\theoremstyle{definition}
\newtheorem{definition}[theorem]{Definition}

\newtheorem{example}[theorem]{Example}
\newtheorem{remark}[theorem]{Remark}
\newtheorem{prob}[theorem]{Problem}

\newcommand{\card}{{\mathrm{card}}\, }

\newcommand{\ds}{{\mathrm{dist}}_S}

\newcommand{\cg}{{\mathcal{G}}}

\newcommand{\qq}{{\mathcal Q}}

\newcommand{\Lab}{\mathrm{Lab}}
\newcommand{\Con}{{\mathrm{Con}}}

\newcommand{\dist}{{\mathrm{dist}}}

\newcommand{\pp}{{\mathcal P}}
\newcommand{\calr}{{\mathcal R}}

\newcommand{\rank}{{\mathrm{rank}}}

\newcommand{\Ball}{{\mathrm{Ball}}}
\def\calp{\mathcal{P}}   
\def\calr{\mathcal{R}}   
\def\calt{\mathcal{T}}   
\def\calc{\mathcal{C}}   

\newcommand {\N}{\mathbb{N}} 

\newcommand {\Z}{\mathbb{Z}}            
\newcommand {\R}{\mathbb{R}} 
\newcommand {\free}{\mathbb{F}} 
\newcommand {\q}{q} 
\newcommand {\g}{ g} 
\newcommand {\pgot}{ p}
\newcommand {\cf}{c}

\newcommand {\iv}{^{-1}}
\newcommand{\lio}[1]{\lm^\omega\left(#1\right)}

\newcommand{\cc}{{\cal C}}

\newcommand{\al}{\alpha}
\newcommand{\e }{\varepsilon }
\newcommand{\C}{\mathcal C}
\newcommand{\CX}{{\rm Con}^{\omega } (X; e, d)}
\newcommand{\CG }{{\rm Con} ^{\omega } (G, d)}

\newcommand{\lo }{\lim ^\omega }
\newcommand{\oom }{o_\omega }
\newcommand{\Oo }{O_\omega }
\newcommand{\To }{\Theta _\omega }
\newcommand{\oas}{$\omega $--almost surely}
\renewcommand{\dh }{{\rm dist}_{Hau} }
\renewcommand{\d }{{\rm dist}}

\newcommand{\CN}{\Con ^\omega (N, d)}
\newcommand{\CGN}{\Con ^\omega (G/N, d)}
\newcommand{\God }{G_e^\omega(d)}
\newcommand{\Nod }{N_e^\omega(d)}
\newcommand{\h }{\varkappa}

\renewcommand{\lio}{\lo }

\newcommand{\ch}{lacunary hyperbolic\,}
\newcommand{\gsc}{graded small cancellation\,}

\begin{document}

\title{Lacunary hyperbolic groups}
 \author{A.Yu. Ol'shanskii, D. V. Osin, M.V. Sapir\thanks{The first and the third
authors were supported in part by the NSF grants DMS 0245600 and DMS
0455881. The second author was supported in part by the NSF grant
DMS 0605093. In addition, the research of the first and the second
author was supported in part by the Russian Fund for Basic Research
grant 05-01-00895, and the research of the third author was
supported by a BSF (USA-Israeli) grant.}}
\date{}
\maketitle
\begin{center}
\large{\textit{with an Appendix by Michael Kapovich and Bruce
Kleiner}}\footnote{The first author of the Appendix was supported in
part by the NSF Grant  DMS-04-05180, the second author was supported
in part by the NSF Grant  DMS-07-01515.}
\end{center}

\begin{abstract}
We call a finitely generated group \ch if one of its asymptotic
cones is an $\R$-tree. We characterize \ch groups as direct limits
of Gromov hyperbolic groups satisfying certain restrictions on the
hyperbolicity constants and injectivity radii. Using central
extensions of \ch groups, we solve a problem of Gromov by
constructing a group whose asymptotic cone $\calc$ has countable
but non-trivial fundamental group (in fact $\calc$ is homeomorphic
to the direct product of a tree and a circle, so
$\pi_1(\calc)=\Z$). We show that the class of \ch groups contains
non-virtually cyclic elementary amenable groups, groups with all
proper subgroups cyclic (Tarski monsters), and torsion groups. We
show that Tarski monsters and torsion groups can have so-called
graded small cancellation presentations, in which case we prove
that all their asymptotic cones are hyperbolic and locally
isometric to trees. This allows us to solve two problems of Dru\c
tu and Sapir, and a problem of Kleiner about groups with cut
points in their asymptotic cones. We also construct a finitely
generated group whose divergence function is not linear but is
arbitrarily close to being linear. This answers a question of
Behrstock.
\end{abstract}

\tableofcontents

\section{Introduction}

\subsection{Asymptotic cones and \ch groups}
\label{acachg}

Asymptotic cones were introduced by Gromov in \cite{Gro81}, a
definition via ultrafilters was given by van den Dries and Wilkie
\cite{DW}. An asymptotic cone of a metric space is, roughly
speaking, what one sees when one looks at the space from infinitely
far away. More precisely, any asymptotic cone of a metric space $(X,
\dist)$ corresponds to an ultrafilter $\omega$, a sequence of
observation points $e=(e_n)_{n\in \N}$ from $X$ and a sequence of
scaling constants $d=(d_n)_{n\in \N}$ diverging to $\infty$. The
cone $\Con^\omega(X;e,d)$ corresponding to $\omega$, $e$ and $d$ is
the ultralimit of the sequence of spaces with basepoints
$(X,\dist/d_n , e_n)$ (see Section \ref{asco} for a precise
definition).

In particular, if $X$ is the Cayley graph of a group $G$ with a word
metric then the asymptotic cones of $X$ are called asymptotic cones
of $G$. For every finitely generated group $G$, its asymptotic cones
are complete geodesic homogeneous metric spaces. Since asymptotic
cones of a group do not depend on the choice of observation points,
we shall omit it from the notation.


The power of asymptotic cones stems from the fact that they capture
both geometric and logical properties of the group, since a large
subgroup $G^\omega_e$ of the ultrapower $G^\omega$ of the group $G$
acts transitively by isometries on the asymptotic cone
$\Con^\omega(G;d)$. Since a large part of the first order theory of
$G$ is inherited by $G^\omega_e$, the isometry group $G^\omega_e$ of
the asymptotic cone ``looks" like $G$. One of the simple but
fundamental applications of asymptotic cones is the following
statement by Gromov \cite{Gr3}: if all asymptotic cones of a group
are simply connected then the group has polynomial isoperimetric and
linear isodiametric functions.

On the other hand, the asymptotic cone captures the coarse
properties of the word metric in $G$. In particular, the asymptotic
cones of two quasi-isometric groups are bi-Lipschitz equivalent.
This makes asymptotic cones very useful tools in proving
quasi-isometric rigidity of some classes of groups \cite{KlL, KaL,
KK, Drutu, DS, D}.

Using asymptotic cones, one can characterize several important
classes of groups. For example, groups of polynomial growth are
precisely groups with all asymptotic cones locally compact
\cite{Gro81, Drutu, Point}.

Another well-known result of Gromov is the following: a finitely
generated group is hyperbolic if and only if all its asymptotic
cones are $\R$-trees \cite{Gr}.

In fact, results of Gromov from \cite{Gr} imply that  a finitely
presented group is hyperbolic if just one of the asymptotic cones is
an $\R$-tree. It was discovered by Kapovich and Kleiner who give a
detailed proof in the Appendix to this paper (see Theorem
\ref{main}). On the other hand, there are non-hyperbolic finitely
generated (but not finitely presented) groups with one asymptotic
cone an $\R$-tree and another one not an $\R$-tree \cite{TV}. We
call a group {{\em \ch} if one of its asymptotic cones is an
$\R$-tree. The term is originated in \cite{Grrand} where sparse
sequences of relations satisfying certain small cancellation
condition as in \cite{TV} are called lacunary.
 Thus a finitely presented \ch group
is hyperbolic.

The following theorem characterizes \ch groups as certain direct
limits of hyperbolic groups. The proof is not too difficult modulo
Theorem \ref{main}, but the result has never been formulated before.

Let $\alpha\colon G\to G'$ be a homomorphism, $G=\la S\ra$. The {\em
injectivity radius} of $\alpha$ is the maximal radius of a ball in
the Cayley graph $\Gamma(G,S)$ where $\alpha$ is injective.

\begin{theorem}\label{dl2}
Let $G$ be a finitely generated group. Then the following conditions
are equivalent.

\begin{enumerate}

\item[1)] $G$ is {\ch}.

\item[2)] There exists a scaling sequence $d=(d_n)$ such that $\CG$ is an $\mathbb R$--tree
for any non--principal ultrafilter $\omega $.

\item[3)] $G$ is the direct limit of a sequence of
hyperbolic groups $G_i=\la S_i\ra$ ($S_i$ is finite) and
epimorphisms
$$
G_1\stackrel{\al_1}\longrightarrow
G_2\stackrel{\al_2}\longrightarrow \ldots, $$ where
$\alpha_i(S_i)=S_{i+1}$, and the hyperbolicity constant of $G_i$
(relative to $S_i$) is ``little $o$" of the injectivity radius of
$\alpha_i$.
\end{enumerate}
\end{theorem}

Note that not every direct limit of hyperbolic groups is \ch. For
example, the free non-cyclic Burnside group of any sufficiently
large odd exponent and the wreath product $(\Z/n\Z)\, {\rm wr}\, \Z$
are direct limits of hyperbolic groups (see \cite{Ivanov} and
\cite{Os02}, respectively) but are not \ch \cite{DS}.

Groups constructed by Thomas and Velickovic \cite{TV} and more
general small cancellation groups given by relations whose lengths
form lacunary sequences of numbers from Gromov \cite[Section
1.7]{Grrand} are \ch (see also Section \ref{lgatcscc} below).

In this paper, we prove that the class of \ch groups is very large:
non-virtually cyclic groups in that class can be elementary
amenable, can have infinite centers, can have all proper subgroups
cyclic (Tarski monsters), and can be torsion groups.

We also show (Theorem \ref{classG7}) that the class of \ch groups
contains all groups given by {\em graded small cancellation}
presentations, a notion originated in \cite{book, Ols}. Moreover,
all asymptotic cones of groups given by graded small cancellation
presentations are hyperbolic and locally isometric to trees (Theorem
\ref{classG7}). Thus methods from \cite{book, Ols} can be used to
construct \ch groups with unusual properties (see Sections
\ref{cscactac} and \ref{chad}).

Theorem \ref{dl2} implies that the torsion-free group $G$ with all
proper subgroups cyclic from the paper \cite{Ol79} of the first
author is \ch. Indeed, $G$ is a direct limit of hyperbolic groups
$G_i$ by \cite[Lemma 9.13]{Ol79}. The same lemma shows that the
injectivity radius of the homomorphism $G_{i-1}\to G_i$ can be
chosen arbitrary large relative to the hyperbolicity constant of
$G_{i-1}$. Similarly, the finitely generated infinite torsion group
with all proper subgroups of prime orders constructed in \cite{Ol80}
is \ch as well.

Although the class of \ch groups is very large, these groups share
some common algebraic properties (see Sections \ref{rach},
\ref{twor}). In particular (Theorem \ref{sub}, Corollary
\ref{corsub}, Theorem \ref{thsub2}), we show that

\begin{itemize}
\item an undistorted subgroup of a \ch group is \ch itself,
\item a \ch group cannot contain a copy of $\Z^2$, an infinite
finitely generated subgroup of bounded torsion and exponential
growth, a copy of the lamplighter group, etc.,
\item  every \ch group is embedded into a relatively hyperbolic 2-generated
\ch group as a peripheral subgroup,
\item any
group that is hyperbolic relative to a \ch subgroup is \ch itself.
\end{itemize}

Theorem \ref{dl2} implies that \ch groups satisfy the Strong Novikov
Conjecture (that is the Baum-Connes assembly map with trivial
coefficients is injective) since all direct limits of hyperbolic
groups satisfy that conjecture (hyperbolic groups satisfy it by
\cite{STY}, and direct limits respect the conjecture by
\cite[Proposition 2.4]{Ro}).

It is also easy to see that the class of \ch groups is closed under
quasi-isometry: indeed, asymptotic cones of quasi-isometric groups
are bi-Lipschitz equivalent if they correspond to the same
ultrafilter and the same sequence of scaling constants. Hence if an
asymptotic cone of one of these groups is an $\R$-tree, then the
other group also has an asymptotic cone that is a $\R$-tree.

\subsection{Central extensions of \ch groups and fundamental groups of asymptotic cones}

One of the interesting properties of hyperbolic groups was
established by Gersten \cite{Ger}: every finitely generated central
extension $H$ of a hyperbolic group $G$ is quasi-isometric to the
direct product of $G$ and the center $Z(H)$. Using some general
properties of asymptotic cones of group extensions (Theorem
\ref{factors}) we establish an asymptotic analog of this result for
central extensions of \ch groups.

\begin{theorem}[Theorem \ref{product}]
Let $N$ be a central subgroup of a finitely generated group $G$
endowed with the induced metric. Suppose that for some
non--principal ultrafilter $\omega$ and some scaling sequence
$d=(d_n)$, $\CGN $ is an $\mathbb R$--tree. Then $\CG $ is
bi--Lipschitz equivalent to $\CN\times \CGN $ endowed with the
product metric.
\end{theorem}

This theorem opens many opportunities to construct asymptotic cones
of groups with unusual properties. Recall that one of the main
problems about asymptotic cones of groups is the following question
by Gromov.

\begin{problem}[Gromov \cite{Gr3}]\label{gr} Is it true that the fundamental
group of an asymptotic cone of a finitely generated group is either
trivial or of cardinality continuum?
\end{problem}

Here is what was known about Problem \ref{gr} before.

\begin{itemize}
\item As we have mentioned above, the triviality of the fundamental
group of all asymptotic cones of a group $G$ implies that the group
is finitely presented, its Dehn function is polynomial and its
isodiametric function is linear. Thus fundamental groups of
asymptotic cones carry important algorithmic information about the
group.

\item By \cite{Pap} if the Dehn function of a finitely presented
group is at most quadratic, then all asymptotic cones are simply
connected. By \cite{OSgafa}, one cannot replace in the previous
statement ``quadratic" by, say, $n^2\log n$.

\item By \cite{Bur}, in many cases asymptotic cones of groups
contain $\pi_1$-embedded Hawaiian earring, and their fundamental
groups are of order continuum (that is true, for example, for
solvable Baumslag-Solitar groups and the Sol group).

\item Non-simply connected asymptotic cones are non-locally compact
\cite{Drutu, Point} but homogeneous, and the isometry groups act on
them with uncountable point stabilizers. Hence every non-trivial
loop in the asymptotic cone typically has uncountably many copies
sharing a common point. This makes a positive answer to Problem
\ref{gr} plausible.

\item In all cases when the non-trivial fundamental groups of
asymptotic cones of groups could be computed, these groups had
cardinality continuum. In \cite{DS}, it is proved that for every
countable group $C$ there exists a finitely generated group $G$
and an asymptotic cone of $G$ whose fundamental group is the free
product of continuously many copies of $C$.
\end{itemize}

Nevertheless, by carefully choosing a central extension of a \ch
group, we answer Problem \ref{gr} negatively.

\begin{theorem}[Theorem \ref{cext2}]\label{cext21}
There exists a finitely generated group $G$ and a scaling sequence
$d=(d_n)$ such that for any ultrafilter $\omega $, $\CG $ is
bi--Lipschitz equivalent to the product of an $\mathbb R$--tree and
$\mathbb S^1$. In particular, $\pi _1 (\CG )=\Z$.
\end{theorem}

\subsection{Cut points in asymptotic cones}

The examples of \ch groups constructed in this paper solve several
problems of Dru\c tu-Sapir \cite{DS} and of Kleiner (see below).

Recall that one of the main applications of asymptotic cones of
groups is the following: if a finitely generated group $H$ has
infinitely many homomorphisms into a locally compact (say, finitely
generated) group $G$ that are pairwise non-conjugate in $G$, then
$H$ acts on an asymptotic cone of $G$ without a global fixed point.
If the asymptotic cone is an $\R$-tree, this implies (using the
theory of groups acting on $\R$-trees due to Rips, Sela,
Bestvina-Feighn, Dunwoody, and others) that $H$ splits into a graph
of groups.

In \cite{DS1}, Dru\c tu and Sapir showed that similar conclusions
can be drawn if all asymptotic cones of $G$ are not trees but only
have global cut points (i.e. points whose removal makes the cones
disconnected). Such groups are called {\em constricted}. In that
case the asymptotic cone is tree-graded in the sense of \cite{DS},
and an action on a tree-graded space under some mild assumptions
leads to an action on an $\R$-tree. In \cite{DS1}, this program
has been carried out for relatively hyperbolic groups (all
asymptotic cones of relatively hyperbolic groups have cut points
by a result of Osin and Sapir \cite{DS}). It is quite plausible
that the program will work also for mapping class groups (where
existence of cut points in asymptotic cones has been proved by
Behrstock \cite{Behr}), fundamental groups of graph manifolds
(their asymptotic cones have cut points by a result of Kapovich,
Kleiner and Leeb \cite{KK}), groups acting $k$-acylindrically on
trees \cite{DMS},  and other groups.

On the other hand many groups do not have cut points in any of
their asymptotic cones. Such groups were called {\em wide} in
\cite{DS}. Among them are non-virtually cyclic groups satisfying
non-trivial laws \cite{DS}, lattices in classical semi-simple Lie
groups of rank $>1$ \cite{DMS}, groups having infinite cyclic
central subgroups \cite{DS}, direct products of infinite groups,
and so on.

Metric spaces whose asymptotic cones do not have cut points (i.e.,
wide spaces) are characterized internally in terms of divergence
in \cite{DMS}. A metric space is {\it wide} if and only if there
are constants $C$, $\e >0$ such that for every three points
$a,b,c$ there exist a path of length at most $C \dist(a,b)$
connecting $a$ with $b$ and avoiding a ball of radius $\e
\dist(c,\{a,b\}) $ about $c$.

One can formulate this condition more precisely in terms of the
divergence function of a metric space.

\begin{definition}\label{div} Let $(X, \dist)$ be a $1$-ended
geodesic metric space, $0<\delta<1$, $\lambda>0$. Let $a,b,c\in X$,
$\min(\dist(c,a), \dist(c,b))=r$. Define $\dv_\lambda(a,b,c;\delta)$ as
the infimum of lengths of paths connecting $a, b$ and avoiding the
ball $\Ball(c,\delta r-\lambda)$ (a ball of non-positive radius is defined to be empty).
Now define the {\em divergence function}
$\dv_\lambda(n;\delta)\colon \R\to \R$ of the space $X$ as the supremum of
all numbers $\dv_\lambda(a,b,c;\delta)$ where $\dist(a,b)\le n$.
\end{definition}

Clearly, the smaller $\delta$, and the bigger $\lambda$, the smaller the functions
$\dv_\lambda(a,b,c;\delta)$ and $\dv_\lambda(n;\delta)$. For 1-ended Cayley graphs, and any $\delta, \delta'<\frac12, \lambda,\lambda'>2$, the functions $\dv_\lambda(n;\delta)$ and $\dv_{\lambda'}(n,\delta')$ are equivalent \cite{DMS}. (Recall that two non-decreasing functions $f,g\colon\N\to \N$ are called equivalent if for some constant $C>1$, we have:
$$f(n/C-C)-Cn-C\le g(n)<f(Cn+C)+Cn+C$$ for every $n$.)

 Hence we can talk about the divergence function $\dv(n)$ of a 1-ended Cayley graph (setting $\delta=\frac13, \lambda=2$).

It is proved in \cite{DMS} that a 1-ended Cayley graph $X$ is wide if and only if
$\dv(n)$ is bounded by a linear function;
and asymptotic cones $\Con^\omega(X, (d_n))$ do not have cut
points for all $\omega$ if and only if for every
$C>1$, the divergence function $\dv(n)$ is uniformly (in $n$) bounded
by a linear function on the intervals $[\frac{d_k}{C}, Cd_k]$.

The divergence function $\dv(n)$ is an interesting
quasi-isometry invariant of a group. It is essentially proved in
\cite{Short} that for every hyperbolic group $\dv(n)$ is
at least exponential. On the other hand, for the mapping class
groups (that also have cut points in all asymptotic cones
\cite{Behr}), the divergence function is quadratic.
The following question was asked by J. Behrstock.
\begin{prob}\label{BehrQ}
Does there exist a group with strictly subquadratic but not linear
divergence function?
\end{prob}

The answer to this question is given below.

If some asymptotic cones of a group $G$ do not have cut points, then
the divergence function of $G$ can be estimated.

\begin{theorem}[Theorem \ref{divergence}]
\label{divergence1} Let $G$ be a 1-ended finitely generated group. Suppose
that for some sequence of scaling constants $d_n$ and every
ultrafilter $\omega$, the asymptotic cone $\Con^\omega(G,(d_n))$
does not have cut points. Let $f(n)\ge n$ be a non-decreasing
function such that $d_n\le f(d_{n-1})$ for all sufficiently large
$n$. Then the divergence function
$\dv(n)$ of $G$ does not exceed $Cf(n)$ for some constant
$C$ (and all $n$).
\end{theorem}

An a-priori stronger property than existence of cut points in
asymptotic cones is the existence of the so-called Morse
quasi-geodesics in the Cayley graph of the group \cite{DMS}. A
quasi-geodesic $\q$ is called {\em Morse} if every
$(L,C)$-quasi-geodesic $\pgot$ with endpoints on the image of $\q$
stays $M$-close to $\q$ where $M$ depends only on $L, C$. By the
Morse lemma, every bi-infinite quasi-geodesic of a hyperbolic
space is Morse. It is proved in \cite{DMS}, that a quasi-geodesic
$\q$ in a metric space $X$ is Morse if and only if in every
asymptotic cone $\cal C$, and every point $m$ in the ultralimit
$\bar\q$ of $\q$, the two halves of $\bar\q$ (before $m$ and after
$m$) are in two different connected components of ${\cal
C}\setminus \{m\}$ (the implication ``$\to$" of this statement was
proved in \cite{Behr}).

Note that similar divergence properties of geodesics have been
studied in the case of CAT(0)-spaces with a co-compact group
action by Ballmann \cite{Bal} and Kapovich-Leeb \cite{KaL}.
In particular, linear (and even subquadratic) divergence for a locally compact Hadamard metric space
(i.e. CAT(0), complete, geodesic, simply connected metric space) implies that
every periodic bi-infinite geodesic in it bounds a
flat half-plane (recall that a geodesic is called {\em periodic} if it is stable under an isometry of the space that acts non-trivially on the geodesic). Ballmann proved \cite[Theorem 3.5]{Bal} that if a
CAT(0)-space $X$ has at least 3 points on the boundary and
contains a bi-infinite periodic geodesic that does not bound a flat
half-plane, then
any sufficiently large (say, co-compact) group of isometries
of $X$ contains a free non-Abelian subgroup.

This leads to the following two problems from \cite{DS}. Recall that
a finitely generated group is {\it constricted} if all its
asymptotic cones have cut points.

\begin{problem}\cite[Problem 1.17]{DS}\label{117} Is every non-constricted
group wide (i.e. if one of the asymptotic cones of a group has no
cut points, does every asymptotic cone of the group have no cut
points)?
\end{problem}

\begin{problem}\cite[Problem 1.19]{DS}\label{119}
Does every non--virtually cyclic finitely generated constricted
group contain free non--abelian subgroups? Is there a constricted
group with all proper subgroups cyclic?
\end{problem}

Bruce Kleiner asked the following stronger questions:

\begin{problem}(Kleiner)\label{Kleiner}
Can a finitely generated group $G$ without free non-cyclic subgroups
contain a bi-infinite Morse quasi-geodesic that is periodic?
Is there a non-wide amenable non-virtually cyclic group?
\end{problem}

We show that the answers to these questions are affirmative.

The second part of Problem \ref{Kleiner} is answered by the
following result (because non-trivial trees have cut points).

\begin{theorem}[Theorems \ref{am}, Lemma \ref{newgr}]\label{am1}
There exists a finitely generated \ch non-virtually cyclic
elementary amenable group $G$. The group $G$ satisfies the following
additional properties:
\begin{itemize}
\item $G$ is 2-generated,
\item $G$ is (locally nilpotent $p$-group)-by-(infinite cyclic),
\item $G$ is residually (finite $p$-group), in particular it is residually nilpotent.
\end{itemize}
\end{theorem}

Note that since $G$ is not hyperbolic (being amenable and
non-virtually cyclic), not all of its asymptotic cones are
$\R$-trees. Hence we obtain the first example of an amenable group
with two non-homeomorphic asymptotic cones.

The following result gives a solution of the first half of that
problem and of Problem \ref{119}.

Recall that a geodesic metric space $X$ is called {\em tree-graded}
with respect to a collection of connected proper subsets $\calp$
\cite{DS} if any two distinct subsets from $\calp$ intersect by at
most one point, and every non-trivial simple geodesic triangle of
$X$ is contained in one of the sets from $\calp$. In particular, if
subsets from $\calp$ are circles (with the natural length metric) of
diameters bounded both from above and from below, we call $X$ a {\em
circle-tree}. It is easy to see that every circle-tree is a
hyperbolic space.

It is proved in \cite{DS} that every (non-singleton) space that is
tree-graded with respect to proper subspaces has cut points.
Conversely, every geodesic space with cut points is tree-graded
with respect to the collection of maximal connected subsets
without cut (their own) cut points \cite{DS}.

\begin{theorem}[Theorem \ref{ExoticQuotients}, Remark \ref{fsub}]\label{main5}
There exist two \ch non-virtually cyclic groups $Q_1$ and $Q_2$ such
that all asymptotic cones of $Q_i$ are circle-trees and
\begin{itemize}
\item[(1)] $Q_1$ is a torsion group.

\item[(2)] Every proper subgroup of $Q_2$ is infinite cyclic,
every infinite periodic path  in the Cayley graph of $Q_2$ is a
Morse quasi-geodesic;
\end{itemize}
\end{theorem}

Note that circle-trees are locally isometric to trees and are
hyperbolic, so all the asymptotic cones of all groups from Theorem
\ref{main5} are locally isometric hyperbolic spaces. But since some
of the asymptotic cones of $Q_i$ are trees and some are not, not all
of the cones are homeomorphic.

The group $Q_1$ is a torsion group but the exponents of elements
in $Q_1$ are not bounded. By the cited result from \cite{DS} about
groups satisfying a law, asymptotic cones of infinite torsion groups of
bounded exponent do not have cut points.

The following theorems give two solutions of Problem \ref{117}. The
first theorem uses central extensions of \ch groups again.

Although any finitely generated group with infinite central cyclic
subgroup is wide \cite{DS}, the next theorem shows that one can
construct \ch groups with infinite (torsion) centers. Such a group
has an asymptotic cone with cut points (a tree) and an asymptotic
cone without cut points.

In fact the information we get is much more precise.

\begin{theorem}[Theorem
\ref{cext1}]\label{cext11} For every $m\ge 2$, there exists a
finitely generated central extension $G$ of a \ch group such that
for any ultrafilter $\omega $ and any scaling sequence $d=(d_n)$,
exactly one of the following possibilities occurs and both of them
can be realized for suitable $\omega $ and $d$.
\begin{enumerate}
\item[(a)] $\CG $ is an $m$--fold cover  of a circle--tree, the
fibers of that cover are cut sets, and every finite cut set of
$\CG$ contains one of the fibers.

\item[(b)] $\CG $ is an $\mathbb R$--tree.
\end{enumerate}
In particular, in both cases $\CG $ is locally isometric to an
$\mathbb R$--tree.
\end{theorem}

\begin{theorem}[Theorem \ref{th2}]\label{th21}
There exists a finitely generated torsion \ch group $G$ such that
one of the asymptotic cones of $G$ does not have cut points.
\end{theorem}

The construction from the proof of Theorem \ref{th2} allows us to
answer Problem \ref{BehrQ}. Indeed, by carefully choosing
exponents of elements of the group $G$ we can control the scaling
constants in the asymptotic cones without cut points. Using
Theorem \ref{divergence1}, we prove (Corollary \ref{cormain}) that
for every function $f(n)$ with $f(n)/n$ non-decreasing, $\lim
f(n)/n=\infty$, there exists a finitely generated torsion group
$G$ whose divergence function $\dv(n)$ is
\begin{itemize}
\item not linear but bounded by a linear function on an infinite
subset of $\N$,
\item bounded from above by $Cf(n)$ for some constant $C$ and all $n$.
\end{itemize}
In addition, one can arrange that the orders of elements $x\in G$
grow with the length $|x|$ as $O(g(|x|)$ for any prescribed in
advance non-decreasing unbounded function $g(n)$. (Recall that
groups with bounded torsion are wide \cite{DS}.)

\subsection{Plan of the paper}

Section \ref{prel} (Preliminaries) contains the main properties of
tree-graded spaces (this makes this paper as independent of
\cite{DS} as possible), the definition and main properties of
asymptotic cones. It also contain some useful properties of
hyperbolic groups.

Section \ref{chgcaae} starts with the definition and basic
properties of \ch groups. In particular, we show (Lemma
\ref{tree-like}) that a group is \ch provided it has a hyperbolic
asymptotic cone, or an asymptotic cone that is locally isometric to
an $\R$-tree. Then we prove the characterization of \ch groups
(Theorem \ref{dirlim}).

The easiest examples of \ch non-hyperbolic groups are groups given
by certain infinite small cancellation presentations. Proposition
\ref{sparse} characterizes such presentations.

In Section \ref{rach} we present several observations connecting
relative hyperbolicity and lacunar hyperbolicity. In particular
Proposition \ref{thsub2} shows that a group that is hyperbolic
relative to a \ch subgroup is \ch itself. This implies that every
\ch group embeds into a 2-generated \ch group (quasi-isometrically,
malnormally, and even as a peripheral subgroup). This result cannot
be generalized to several subgroups: we show (Example \ref{315})
that even a free product of two \ch groups can be non-\ch.

In Section \ref{twor}, we provide several general properties of
subgroups of \ch groups already mentioned above in Section
\ref{acachg}.

In Section \ref{chag}, we construct elementary amenable \ch groups
and prove Theorem \ref{am1}.

Section \ref{cscactac} is devoted to several small cancellation
conditions and their applications. We start by introducing a small
cancellation condition $C(\e, \mu,\rho)$ for presentations over any
group $H$ (i.e. presentations of factor groups of $H$ instead of
just factor-groups of the free group as in the classical case). We
show (Lemma \ref{gamma-cell}) that if the group $H$ is hyperbolic
and the cancellation parameters are appropriately chosen, then the
factor-group satisfies an analog of the Greendlinger lemma, and is
hyperbolic again with a nice control on the hyperbolicity constant.

This allows us to use induction, and introduce direct limits of
groups $G_1\to G_2\to...$ where each $G_{i+1}$ is given by a
presentation over $G_i$ satisfying an appropriate
$C(\e,\mu,\rho)$-condition so that Lemma \ref{gamma-cell} holds.
The union of presentations of all $G_i$ gives us a presentation of
the limit group $G$. We say that such a presentation satisfies a
{\em graded small cancellation} condition.  We prove (Corollary
\ref{subset}) that many infinitely presented groups with classical
small cancellation conditions have \gsc presentations.

Theorem \ref{classG7} gives an important property of \gsc
presentations: every asymptotic cone of a group given by a \gsc
presentation is a circle-tree or an $\R$-tree. Moreover, given the
parameters of the cone, one can tell which of these options holds,
and what are the sizes of the circles in the circle-tree.

In Section \ref{gwfs}, we apply results from \cite{Ols} and show
that there are non--virtually cyclic groups with \gsc presentations that have all proper
subgroups infinite cyclic or all proper subgroups finite (Theorem
\ref{ExoticQuotients}, Remark \ref{fsub}).

In Section \ref{fb}, we notice that existence of cut points in all
asymptotic cones follows from the non-triviality of the Floyd
boundary of a group. The converse statement does not hold as
follows from Theorem \ref{ExoticQuotients}.

In Section \ref{ceochg}, we first establish some very general
results about asymptotic cones of group extensions. In particular
(Theorem \ref{factors}), if
$$1\to N\to G\to H\to 1$$ is an exact sequence then there exists a
continuous map from a cone of $G$ to a cone of $H$ (corresponding to
the same parameters) with fibers homeomorphic to the cone of $N$
(considered as a subspace of $G$). In the case of central
extensions, the situation is much nicer, and in particular when the
cone of $H$ is an $\R$-tree, the fibration becomes trivial and the
cone of $G$ becomes bi-Lipschitz equivalent to the direct product of
the $\R$-tree and the cone of $N$ (Theorem \ref{product}). As
applications of these general results, we give proofs of Theorems
\ref{cext21} and  \ref{cext11}.

Section \ref{chad} is devoted to torsion groups and the proofs of
Theorem \ref{th21} and Corollary \ref{cormain} solving the slow
divergence problem.

Section \ref{op} contains some open problems.

The Appendix written by M. Kapovich and B. Kleiner contains the
proofs of Theorem \ref{main} (that a finitely presented lacunary
hyperbolic group is hyperbolic) and other useful results about
asymptotic cones of finitely presented groups.

{\bf Acknowledgement.} The authors are grateful to Cornelia Dru\c
tu and Bruce Kleiner for very fruitful conversations. We are also
grateful to Michael Kapovich and Bruce Kleiner for adding their
unpublished results as an Appendix to our paper.


\section{Preliminaries}
\label{prel}



\subsection{Cayley graphs and van Kampen diagrams}
\label{ctgs}


Given a word $W$ in an alphabet $S$, we denote by $|W| $ its length.
We also write $W\equiv V$ to express the letter--for--letter
equality of words $W$ and $V$.

Let $G$ be a group generated by a set $S$. Recall that the {\it
Cayley graph} $\Gamma (G,S)$ of a group $G$ with respect to the set
of generators $S$ is an oriented labeled 1--complex with the vertex
set $V(\Gamma (G,S))=G$ and the edge set $E(\Gamma (G,S) )=G\times
S^{\pm 1}$. An edge $e=(g,a)$ goes from the vertex $g$ to the vertex
$ga$ and has label $\Lab (e)\equiv a$. As usual, we denote the
initial and the terminal vertices of the edge $e$ by $e_-$ and $e_+$
respectively. Given a combinatorial path $p=e_1\ldots e_k$ in the
Cayley graph $\Gamma (G,S) $, where $e_1,  \ldots , e_k\in E(\Gamma
(G,S) )$, we denote by $\Lab(p)$ its label. By definition, $\Lab
(p)\equiv \Lab (e_1)\ldots \Lab (e_k).$ We also denote by
$p_-=(e_1)_-$ and $p_+=(e_k)_+$ the initial and terminal vertices of
$p$ respectively. The length $|p|$ of $p$ is the number of edges in
$p$.

The {\it (word) length} $|g|$ of an element $g\in G$ with respect to
the generating set $S$ is defined to be the length of a shortest
word in $S$ representing $g$ in $G$. The formula
$dist(f,g)=|f^{-1}g|$ defines a metric on $G$. We also denote by
$dist $ the natural extension of the metric to $\Gamma (G,S)$.

Recall that a {\it van Kampen diagram} $\Delta $ over a presentation
\begin{equation}
G=\langle S\; | \; \mathcal R\rangle \label{ZP}
\end{equation}
is a finite oriented connected planar 2--complex endowed with a
labeling function $\Lab : E(\Delta )\to S^{\pm 1}$, where $E(\Delta
) $ denotes the set of oriented edges of $\Delta $, such that $\Lab
(e^{-1})\equiv (\Lab (e))^{-1}$. Given a cell $\Pi $ of $\Delta $,
we denote by $\partial \Pi$ the boundary of $\Pi $; similarly,
$\partial \Delta $ denotes the boundary of $\Delta $. The labels of
$\partial \Pi $ and $\partial \Delta $ are defined up to cyclic
permutations. An additional requirement is that the label of any
cell $\Pi $ of $\Delta $ is equal to (a cyclic permutation of) a
word $R^{\pm 1}$, where $R\in \mathcal R$. Labels and lengths of
paths are defined as in the case of Cayley graphs.

The van Kampen Lemma states that a word $W$ over an alphabet $S$
represents the identity in the group given by (\ref{ZP}) if and only
if there exists a connected simply--connected planar diagram $\Delta
$ over (\ref{ZP}) such that $\Lab (\partial \Delta )\equiv W$
\cite[Ch. 5, Theorem 1.1]{LS}.


\subsection{Tree-graded spaces}
\label{tgs}


Here we collect all the necessary definitions and basic properties
of tree-graded spaces from \cite{DS} needed in this paper.

\begin{definition}\label{tgspace}
Let $\free$ be a complete geodesic metric space and let $\pp$ be a
collection of closed geodesic non-empty subsets (called
{\it{pieces}}). Suppose that the following two properties are
satisfied:

\begin{enumerate}

\item[($T_1$)] Every two different pieces have at most one common
point.

\item[($T_2$)] Every non--trivial simple geodesic triangle (a
simple loop composed of three geodesics) in $\free$ is contained
in one piece.
\end{enumerate}

Then we say that the space $\free$ is {\em tree-graded with
respect to }$\pp$.
\end{definition}

\noindent  For technical reasons it is convenient to allow $\mathcal
P$ to be empty. Clearly $\mathbb F$ is tree--graded with respect to
the empty collections of pieces only if $\mathbb F$ is a tree.

By \cite[Proposition 2.17]{DS}, property $(T_2)$ in this
definition can be replaced by each of the following two
properties.

\begin{enumerate}
\item[{\rm ($T_2'$)}]
 For every topological arc $\cf:[0,d]\to \free$, where $\cf(0)\ne \cf(d)$, and any $t\in
[0,d]$, let $\cf[t-a,t+b]$ be a maximal sub-arc of $\cf$
containing $\cf (t)$ and contained in one piece. Then every other
topological arc with the same endpoints as $\cf$ must contain the
points $\cf (t-a)$ and $\cf (t+b)$.
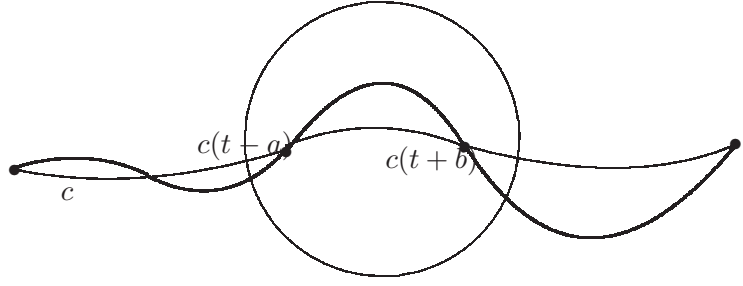
\begin{figure}[!ht]
\centering
\unitlength .85mm 
\linethickness{0.4pt}
\ifx\plotpoint\undefined\newsavebox{\plotpoint}\fi 
\begin{picture}(128.04,45.33)(0,0)
\put(93.47,23.86){\line(0,1){.97}}
\put(93.44,24.83){\line(0,1){.968}}
\put(93.38,25.8){\line(0,1){.964}}
\multiput(93.27,26.77)(-.0306,.19171){5}{\line(0,1){.19171}}
\multiput(93.11,27.73)(-.0327,.15844){6}{\line(0,1){.15844}}
\multiput(92.92,28.68)(-.02987,.1176){8}{\line(0,1){.1176}}
\multiput(92.68,29.62)(-.03125,.10323){9}{\line(0,1){.10323}}
\multiput(92.4,30.55)(-.0323,.09154){10}{\line(0,1){.09154}}
\multiput(92.08,31.46)(-.03309,.0818){11}{\line(0,1){.0818}}
\multiput(91.71,32.36)(-.03369,.07354){12}{\line(0,1){.07354}}
\multiput(91.31,33.24)(-.031701,.061663){14}{\line(0,1){.061663}}
\multiput(90.86,34.11)(-.032159,.056156){15}{\line(0,1){.056156}}
\multiput(90.38,34.95)(-.032499,.051229){16}{\line(0,1){.051229}}
\multiput(89.86,35.77)(-.032736,.046783){17}{\line(0,1){.046783}}
\multiput(89.3,36.56)(-.032883,.042741){18}{\line(0,1){.042741}}
\multiput(88.71,37.33)(-.032951,.039042){19}{\line(0,1){.039042}}
\multiput(88.09,38.08)(-.032948,.035637){20}{\line(0,1){.035637}}
\multiput(87.43,38.79)(-.032882,.032486){21}{\line(-1,0){.032882}}
\multiput(86.74,39.47)(-.036033,.032515){20}{\line(-1,0){.036033}}
\multiput(86.02,40.12)(-.039438,.032476){19}{\line(-1,0){.039438}}
\multiput(85.27,40.74)(-.043136,.032364){18}{\line(-1,0){.043136}}
\multiput(84.49,41.32)(-.047176,.032167){17}{\line(-1,0){.047176}}
\multiput(83.69,41.87)(-.051618,.031877){16}{\line(-1,0){.051618}}
\multiput(82.86,42.38)(-.060579,.033726){14}{\line(-1,0){.060579}}
\multiput(82.01,42.85)(-.066814,.033334){13}{\line(-1,0){.066814}}
\multiput(81.15,43.28)(-.07394,.0328){12}{\line(-1,0){.07394}}
\multiput(80.26,43.68)(-.0822,.0321){11}{\line(-1,0){.0822}}
\multiput(79.35,44.03)(-.09192,.03119){10}{\line(-1,0){.09192}}
\multiput(78.44,44.34)(-.1036,.03){9}{\line(-1,0){.1036}}
\multiput(77.5,44.61)(-.1348,.03251){7}{\line(-1,0){.1348}}
\multiput(76.56,44.84)(-.15883,.03078){6}{\line(-1,0){.15883}}
\multiput(75.61,45.02)(-.19207,.02828){5}{\line(-1,0){.19207}}
\put(74.65,45.17){\line(-1,0){.966}}
\put(73.68,45.26){\line(-1,0){.969}}
\put(72.71,45.32){\line(-1,0){.971}}
\put(71.74,45.33){\line(-1,0){.97}}
\put(70.77,45.29){\line(-1,0){.968}}
\multiput(69.8,45.22)(-.2408,-.0303){4}{\line(-1,0){.2408}}
\multiput(68.84,45.09)(-.19133,-.03292){5}{\line(-1,0){.19133}}
\multiput(67.88,44.93)(-.13546,-.02967){7}{\line(-1,0){.13546}}
\multiput(66.93,44.72)(-.11723,-.03129){8}{\line(-1,0){.11723}}
\multiput(66,44.47)(-.10284,-.0325){9}{\line(-1,0){.10284}}
\multiput(65.07,44.18)(-.09114,-.0334){10}{\line(-1,0){.09114}}
\multiput(64.16,43.85)(-.07461,-.03124){12}{\line(-1,0){.07461}}
\multiput(63.26,43.47)(-.067501,-.031922){13}{\line(-1,0){.067501}}
\multiput(62.39,43.06)(-.061275,-.032445){14}{\line(-1,0){.061275}}
\multiput(61.53,42.6)(-.055762,-.032836){15}{\line(-1,0){.055762}}
\multiput(60.69,42.11)(-.050832,-.033116){16}{\line(-1,0){.050832}}
\multiput(59.88,41.58)(-.046384,-.033299){17}{\line(-1,0){.046384}}
\multiput(59.09,41.01)(-.04234,-.033398){18}{\line(-1,0){.04234}}
\multiput(58.33,40.41)(-.03864,-.033421){19}{\line(-1,0){.03864}}
\multiput(57.59,39.78)(-.035235,-.033377){20}{\line(-1,0){.035235}}
\multiput(56.89,39.11)(-.03369,-.034936){20}{\line(0,-1){.034936}}
\multiput(56.22,38.41)(-.032077,-.036423){20}{\line(0,-1){.036423}}
\multiput(55.57,37.68)(-.031997,-.039828){19}{\line(0,-1){.039828}}
\multiput(54.97,36.93)(-.033712,-.046085){17}{\line(0,-1){.046085}}
\multiput(54.39,36.14)(-.033569,-.050534){16}{\line(0,-1){.050534}}
\multiput(53.86,35.33)(-.033333,-.055467){15}{\line(0,-1){.055467}}
\multiput(53.36,34.5)(-.032991,-.060983){14}{\line(0,-1){.060983}}
\multiput(52.89,33.65)(-.032523,-.067213){13}{\line(0,-1){.067213}}
\multiput(52.47,32.77)(-.03191,-.07433){12}{\line(0,-1){.07433}}
\multiput(52.09,31.88)(-.0311,-.08258){11}{\line(0,-1){.08258}}
\multiput(51.75,30.97)(-.03341,-.10255){9}{\line(0,-1){.10255}}
\multiput(51.45,30.05)(-.03234,-.11695){8}{\line(0,-1){.11695}}
\multiput(51.19,29.12)(-.03088,-.13519){7}{\line(0,-1){.13519}}
\multiput(50.97,28.17)(-.02886,-.15919){6}{\line(0,-1){.15919}}
\multiput(50.8,27.21)(-.0324,-.2405){4}{\line(0,-1){.2405}}
\put(50.67,26.25){\line(0,-1){.967}}
\put(50.58,25.28){\line(0,-1){3.877}}
\multiput(50.68,21.41)(.0332,-.2404){4}{\line(0,-1){.2404}}
\multiput(50.81,20.45)(.02936,-.15909){6}{\line(0,-1){.15909}}
\multiput(50.98,19.49)(.0313,-.13509){7}{\line(0,-1){.13509}}
\multiput(51.2,18.55)(.03271,-.11684){8}{\line(0,-1){.11684}}
\multiput(51.47,17.61)(.03374,-.10244){9}{\line(0,-1){.10244}}
\multiput(51.77,16.69)(.03137,-.08248){11}{\line(0,-1){.08248}}
\multiput(52.11,15.78)(.03214,-.07423){12}{\line(0,-1){.07423}}
\multiput(52.5,14.89)(.032736,-.067109){13}{\line(0,-1){.067109}}
\multiput(52.93,14.02)(.033184,-.060878){14}{\line(0,-1){.060878}}
\multiput(53.39,13.17)(.033509,-.055361){15}{\line(0,-1){.055361}}
\multiput(53.89,12.34)(.033729,-.050428){16}{\line(0,-1){.050428}}
\multiput(54.43,11.53)(.031977,-.043423){18}{\line(0,-1){.043423}}
\multiput(55.01,10.75)(.032123,-.039726){19}{\line(0,-1){.039726}}
\multiput(55.62,9.99)(.032192,-.036322){20}{\line(0,-1){.036322}}
\multiput(56.26,9.27)(.032191,-.033171){21}{\line(0,-1){.033171}}
\multiput(56.94,8.57)(.035341,-.033265){20}{\line(1,0){.035341}}
\multiput(57.65,7.91)(.038746,-.033298){19}{\line(1,0){.038746}}
\multiput(58.38,7.27)(.042446,-.033263){18}{\line(1,0){.042446}}
\multiput(59.15,6.67)(.046489,-.033152){17}{\line(1,0){.046489}}
\multiput(59.94,6.11)(.050937,-.032955){16}{\line(1,0){.050937}}
\multiput(60.75,5.58)(.055866,-.03266){15}{\line(1,0){.055866}}
\multiput(61.59,5.09)(.061377,-.032251){14}{\line(1,0){.061377}}
\multiput(62.45,4.64)(.067601,-.031708){13}{\line(1,0){.067601}}
\multiput(63.33,4.23)(.07471,-.031){12}{\line(1,0){.07471}}
\multiput(64.22,3.86)(.09125,-.03311){10}{\line(1,0){.09125}}
\multiput(65.14,3.53)(.10294,-.03217){9}{\line(1,0){.10294}}
\multiput(66.06,3.24)(.11733,-.03092){8}{\line(1,0){.11733}}
\multiput(67,2.99)(.13555,-.02924){7}{\line(1,0){.13555}}
\multiput(67.95,2.78)(.19143,-.03231){5}{\line(1,0){.19143}}
\put(68.91,2.62){\line(1,0){.963}}
\put(69.87,2.5){\line(1,0){.968}}
\put(70.84,2.43){\line(1,0){.97}}
\put(71.81,2.4){\line(1,0){.971}}
\put(72.78,2.41){\line(1,0){.969}}
\put(73.75,2.47){\line(1,0){.965}}
\multiput(74.71,2.57)(.19198,.02889){5}{\line(1,0){.19198}}
\multiput(75.67,2.72)(.15873,.03128){6}{\line(1,0){.15873}}
\multiput(76.63,2.9)(.1347,.03294){7}{\line(1,0){.1347}}
\multiput(77.57,3.13)(.1035,.03033){9}{\line(1,0){.1035}}
\multiput(78.5,3.41)(.09182,.03148){10}{\line(1,0){.09182}}
\multiput(79.42,3.72)(.0821,.03236){11}{\line(1,0){.0821}}
\multiput(80.32,4.08)(.07384,.03304){12}{\line(1,0){.07384}}
\multiput(81.21,4.47)(.066708,.033545){13}{\line(1,0){.066708}}
\multiput(82.07,4.91)(.056441,.031657){15}{\line(1,0){.056441}}
\multiput(82.92,5.38)(.051517,.03204){16}{\line(1,0){.051517}}
\multiput(83.75,5.9)(.047074,.032317){17}{\line(1,0){.047074}}
\multiput(84.55,6.45)(.043033,.0325){18}{\line(1,0){.043033}}
\multiput(85.32,7.03)(.039335,.032601){19}{\line(1,0){.039335}}
\multiput(86.07,7.65)(.03593,.032629){20}{\line(1,0){.03593}}
\multiput(86.79,8.3)(.032779,.03259){21}{\line(1,0){.032779}}
\multiput(87.47,8.99)(.032835,.035741){20}{\line(0,1){.035741}}
\multiput(88.13,9.7)(.032827,.039146){19}{\line(0,1){.039146}}
\multiput(88.75,10.45)(.032748,.042845){18}{\line(0,1){.042845}}
\multiput(89.34,11.22)(.032587,.046887){17}{\line(0,1){.046887}}
\multiput(89.9,12.02)(.032336,.051332){16}{\line(0,1){.051332}}
\multiput(90.42,12.84)(.031982,.056257){15}{\line(0,1){.056257}}
\multiput(90.9,13.68)(.031506,.061763){14}{\line(0,1){.061763}}
\multiput(91.34,14.54)(.03346,.07364){12}{\line(0,1){.07364}}
\multiput(91.74,15.43)(.03283,.08191){11}{\line(0,1){.08191}}
\multiput(92.1,16.33)(.03201,.09164){10}{\line(0,1){.09164}}
\multiput(92.42,17.25)(.03092,.10333){9}{\line(0,1){.10333}}
\multiput(92.7,18.18)(.03371,.13451){7}{\line(0,1){.13451}}
\multiput(92.93,19.12)(.0322,.15854){6}{\line(0,1){.15854}}
\multiput(93.13,20.07)(.02999,.19181){5}{\line(0,1){.19181}}
\put(93.28,21.03){\line(0,1){.965}}
\put(93.38,21.99){\line(0,1){1.871}}
\put(57,21.86){\circle*{1.8}} \put(84.75,22.11){\circle*{.71}}
\put(85,22.36){\circle*{1}}
\qbezier(14.5,19.11)(33.88,15.11)(56.75,22.11)
\qbezier(57,22.86)(70.75,28.36)(84.5,22.86)
\qbezier(84.75,22.86)(110.88,15.74)(127.5,23.11)
\put(14.5,19.11){\circle*{1.58}}
\put(127.25,23.11){\circle*{1.58}} \thicklines
\qbezier(14.25,19.36)(25.75,22.86)(35.25,18.36)
\qbezier(57,21.86)(72.88,43.24)(85.25,22.11)
\put(22.75,15.36){\makebox(0,0)[cc]{$\cf$}}
\put(50.61,22.76){\makebox(0,0)[cc]{$\cf(t-a)$}}
\put(79.81,20.55){\makebox(0,0)[cc]{$\cf(t+b)$}}
\put(84.93,22.6){\circle*{1.69}}
\qbezier(35.32,18.08)(47.35,11.98)(56.66,21.86)
\qbezier(85.04,22.7)(103.43,-5.99)(127.29,22.91)
\end{picture}
\caption{Property ($T_2'$).} \label{fig1.1}
\end{figure}
\item[{\rm ($T_2''$)}] Every simple loop in $\free$ is contained
in one piece.
\end{enumerate}

In order to avoid problems with pieces that are singletons, we
shall always assume that pieces in a tree-graded space cannot
contain each other.

Let us define a partial order relation on the set of tree-graded
structures of a space. If $\pp$ and $\pp'$ are collections of
subsets of $X$, and a space $X$ is tree-graded with respect to
both $\pp$ and $\pp'$, we write $\pp \prec \pp'$ if for every set
$M \in \pp$ there exists $M'\in \pp'$ such that $M\subset M'$. The
relation $\prec$ is a partial order because by our convention
pieces of $\pp$ (resp. $\pp'$) cannot contain each other.

\begin{lemma}(\cite[Lemma 2.31]{DS}) \label{cutting}
Let $X$ be a complete geodesic metric space containing at least
two points and let $\calc$ be a non-empty set of global cut points
in $X$.

\begin{enumerate}
\item[(a)] There exists the largest in the sense of $\prec$
collection $\calp$ of subsets of $X$ such that
\begin{itemize}
\item $X$ is tree-graded with respect to $\pp$; \item any piece in
$\pp$ is either a singleton or a set with no global cut-point from
$\calc$.
\end{itemize}
Moreover the intersection of any two distinct pieces from $\calp$
is either empty or a point from $\calc$. \item[(b)] Let $X$ be a
homogeneous space with a cut-point. Then every point in $X$ is a
cut-point, so let $\calc=X$.  Let $\pp$ be the set of pieces
defined in part (a).  Then for every $M\in \calp$ every $x\in M$
is the {\em projection} of a point $y\in X\setminus M$ onto $M$
(i.e. the closest to $y$ point in $M$).
\end{enumerate}
\end{lemma}

\begin{lemma}(\cite[Lemma 2.15]{DS})\label{cut}
Let $\free$ be a tree-graded metric space. Let $A$ be a path
connected subset of $\free$ without cut points. Then $A$ is
contained in a piece. In particular every simple loop in $\free$
is contained in a piece of $\free $.
\end{lemma}

\begin{lemma}(\cite[Lemma 2.28]{DS})\label{sir}
Let $\g=\g_1\g_2\dots\g_{2m}$ be a curve in a tree-graded space
$\free$ which is a composition of geodesics. Suppose that all
geodesics $\g_{2k}$ with $k\in \{ 1,\dots ,m-1 \}$ are non-trivial
and for every $k\in \{ 1,\dots ,m \}$ the geodesic $g_{2k}$ is
contained in a piece $M_k$ while for every $k\in \{ 0,1,\dots ,m-1
\}$ the geodesic $\g_{2k+1}$ intersects $M_k$ and $M_{k+1}$ only
in its respective endpoints. In addition assume that if
$\g_{2k+1}$ is empty then $M_k\ne M_{k+1}$. Then $\g$ is a
geodesic.
\end{lemma}

\begin{lemma}(\cite[Corollary 2.10]{DS}) \label{strconv}
\begin{itemize}
  \item[(1)] Every simple path in $\free$ joining two points in a piece is contained in the piece.
  \item[(2)] Every non-empty intersection between a simple path in $\free$ and a
             piece is a subpath.
\end{itemize}
\end{lemma}

\begin{lemma}(\cite[Corollary 2.11]{DS}) \label{projA} Let $A$ be a connected
subset (possibly a   point) in a tree-graded space $\free $ which
intersects a piece $M$ in at most one point.
\begin{itemize}
\item[(1)] The subset $A$ projects onto $M$ in a unique point $x$.
\item[(2)] Every path joining a point in $A$ with a point in $M$
contains $x$.
\end{itemize}
\end{lemma}


\subsection{Asymptotic cones}
\label{asco}


 Let us recall the definition of asymptotic cones. A
non-principal ultrafilter $\omega$ is a finitely additive measure
defined on all subsets $S$ of ${\mathbb N}$, such that $\omega(S)
\in \{0,1\}$, $\omega (\mathbb N)=1$, and $\omega(S)=0$ if $S$ is
a finite subset. For a bounded sequence of numbers $x_n$, $n\in
\N$, the limit $\lo x_n$ with respect to $\omega$ is the unique
real number $a$ such that $\omega(\{i\in {\mathbb N}:
|x_i-a|<\epsilon\})=1$ for every $\epsilon>0$. Similarly, $\lo
x_n=\infty $ if $\omega(\{i\in {\mathbb N}: x_i>M \})=1$ for every
$M>0$.

Given two infinite sequences of real numbers $(a_n)$ and $(b_n)$
we write $a_n=\oom (b_n)$ if $\lo a_n/b_n =0$. Similarly $a_n=\To
(b_n)$ (respectively $a_n=\Oo (b_n)$) means that $0<\lo
(a_n/b_n)<\infty $ (respectively $\lo (a_n/b_n)<\infty $).

Let $(X_n,\dist_n)$, $n\in\N$, be a metric space. Fix an arbitrary
sequence $e=(e_n)$ of points $e_n\in X_n$. Consider the set ${\cal
F}$ of sequences $g=(g_n)$, $g_n\in X_n$, such that
$\dist_n(g_n,e_n) \le c $ for some constant $c=c(g)$. Two sequences
$(f_n)$ and $(g_n)$ of this set ${\cal F}$ are said to be {\em
equivalent} if $\lio\dist_n(f_n,g_n) =0$. The equivalence class of
$(g_n)$ is denoted by $(g_n)^\omega$. The $\omega$-{\em limit} $\lio
(X_n)_e$ is the quotient space of equivalence classes where the
distance between $(f_n)^\omega$ and $(g_n)^\omega$ is defined as
$\lio\dist(f_n,g_n)$.

An {\em asymptotic cone} $\Con^\omega(X,e, d)$ of a metric space
$(X,\dist)$ where $e=(e_n)$, $e_n\in X$, and $d=(d_n)$ is an unbounded
non-decreasing {\it scaling sequence} of positive real numbers,
is the $\omega$-limit of spaces
$X_n=(X,\dist/d_n)$. The asymptotic cone is a complete space; it is
a geodesic metric space if $X$ is a geodesic metric space
(\cite{Gr3,Drutu}). Note that $\Con^\omega(X, e,d)$ does not depend
on the choice of $e$ if $X$ is homogeneous (say, if $X$ is a
finitely generated group with a word metric), so in that case, we
shall omit $e$ in the notation of an asymptotic cone.

If $(Y_n)$ is a sequence of subsets of $X$ endowed with the induced
metric, we define $\lio(Y_n)_e$ to be the subset of $\Con^\omega(X,
e,d)$ consisting of $x\in \Con^\omega(X, e,d)$ that can be
represented by sequences $(x_n)$, where $x_n\in Y_n$.

An {\em asymptotic cone} of a finitely generated group $G$ with a
word metric is the asymptotic cone of its Cayley graph (considered
as the discrete space of vertices with the word metric). Asymptotic
cones corresponding to two different finite generating sets of $G$
(and the same ultrafilters and scaling constants) are bi-Lipschitz
equivalent. The asymptotic cone $\Con^\omega(G,d)$ of a group $G$ is
a homogeneous geodesic metric space with transitive group of
isometries $\God$ consisting of sequences $(g_n)$, ${g_n\in G}$ such
that $|g_n|\le Cd_n$ for some constant $C$ depending on the sequence
(here $|g_n|$ is the word length of $g_n$). The action is by
multiplication on the left: $(g_n)\circ (h_n)^\omega=(g_n
h_n)^\omega$.

Recall that a geodesic $p$ in $\CX $ is called a {\it limit
geodesic} if $p =\lo p _n$, where for every $n\in
\mathbb N$,  $p_n$ is a geodesic in $X$. The lemma below was
proved in \cite[Corollary 4.18]{D}.

\begin{lemma} \label{p3} Assume that in an asymptotic
cone $\CX $, a collection of closed subsets $\mathcal P$ satisfies
$(T_1)$ and every non--trivial simple triangle in $\CX $ whose
sides are limit geodesics is contained in a subset from $\mathcal
P$. Then $\mathcal P$ satisfies $(T_2)$, i.e., $\CX $ is
tree--graded with respect to $\mathcal P$.
\end{lemma}


\subsection{Hyperbolic groups}


Recall that a geodesic space $X$ is {\it $\delta $--hyperbolic} (or
simply {\it hyperbolic}, for brevity) if for any geodesic triangle
$\Delta $ in $X$, each side of $\Delta $ is contained in the closed
$\delta$--neighborhood of the union of the other two sides. This
$\delta$ is called the {\em hyperbolicity constant} of $X$. A group
$H$ is {\it $\delta $--hyperbolic} (or simply {\it hyperbolic}) if
it is generated by a finite set $S$ and its Cayley graph $\Gamma(H,
S)$ endowed with the combinatorial metric is a hyperbolic metric
space.

Recall that a path $p$ is called $(\lambda , c)$--quasi--geodesic
for some $\lambda \in (0,1]$, $c\ge 0$ if for any subpath $q$ of
$p$, we have $$\d (q_-, q_+)\ge \lambda |q| -c.$$ The property of
hyperbolic spaces stated below is well--known although it is
usually formulated in a slightly different manner (see, for
example, \cite{GH}).

\begin{lemma}\label{qg}
For any $\lambda \in (0, 1]$, $c>0$ there exists $\theta( \lambda,
c)$ such that any two $(\lambda , c)$--quasi--geodesic paths $p$,
$q$ in a $\delta $--hyperbolic metric space such that $p_-=q_-$
and $p_+=q_+$ belong to the closed $\theta (\lambda
,c)$--neighborhoods of each other.
\end{lemma}

The next property can easily be derived from the definition of a
hyperbolic space by cutting the $n$--gon into triangles.

\begin{lemma}\label{quad}
For any $n\ge 3$, any side of a geodesic $n$--gon in a $\delta
$--hyperbolic space belongs to the closed $(n-2)\delta
$--neighborhood of the other $(n-1)$ sides.
\end{lemma}


\section{Lacunary hyperbolic groups: characterization and examples}


\label{chgcaae}


\subsection{A characterization of lacunary hyperbolic groups}

\label{acochg}

We say that a metric space $X$ is {\it {\ch}} if one of the
asymptotic cones of $X$ is an $\mathbb R$--tree. In particular,
every hyperbolic metric space is {\ch}. A group $G$ is {\it {\ch}}
if it is finitely generated and the corresponding Cayley graph is
{\ch}. Clearly this notion is independent of the choice of the
finite generating set. We also say that a metric space $X$ is {\it
almost homogeneous} if there is a homogeneous subspace $Y\subseteq
X$ such that $\dist _{Hau} (X, Y)<\infty $, where $\dist_{Hau}$ is
the Hausdorff distance between $X$ and $Y$. That is, there exists
$\e>0$ such that for every $x\in X$, we have $\dist (x, Y)<\e $.
Given a group $G$, any Cayley graph of $G$ endowed with the
combinatorial metric is almost homogeneous. Note also that if $X$
is almost homogeneous, then every asymptotic cone of $X$ is
homogeneous.

\begin{lemma}\label{tree-like}
Let $X$ be a metric space. Then the following properties are
equivalent.

\begin{enumerate}
\item[1)] Some asymptotic cone of $X$ is an $\mathbb R$--tree.

\item[2)] Some asymptotic cone of $X$ is {\ch}.

\end{enumerate}
If, in addition, $X$ is almost homogeneous, these properties are
equivalent to
\begin{enumerate}
\item[3)] Some asymptotic cone of $X$ is locally isometric to an
$\mathbb R$--tree.
\end{enumerate}
\end{lemma}

\begin{proof}
It suffices to show that $2)\Rightarrow 1)$ and $3) \Rightarrow 1)$.
Recall that the set $\mathcal C(X)$ of asymptotic cones of $X$ is
closed under taking ultralimits \cite[Corollary 3.24]{DS}. Further
let $\dist $ denote the standard metric on $\CX$. It is
straightforward to check that for every $k>0$, $(\CX, \frac1k \dist
)$ is isometric to the asymptotic cone ${\rm Con} (X, e, (kd_n))$.
In particular, if $Y\in \mathcal C(X)$, then every asymptotic cone
of $Y$ belongs to $C(X)$. This yields $2)\Rightarrow 1)$.

Further if $\CX $ is homogeneous and locally isometric to an
$\mathbb R$--tree, there is $c>0$ such that $\CX $ contains no
simple nontrivial loops of length at most $c$. We consider the
sequence of cones $C_k={\rm Con} (X, e, (kd_n))$ for $k\to 0$.
Clearly $C_k$ has no simple nontrivial loops of length at most
$c/k$. Then for any non--principal ultrafilter $\omega$,
$\lio(C_k)_e$ is an $\mathbb R$--tree.
\end{proof}

By Theorem \ref{main} of Kapovich and Kleiner, if a finitely
presented group $G$ is {\ch}, then, in fact, it is hyperbolic.
Combining this with Lemma \ref{tree-like}, we obtain the following.

\begin{proposition}\label{fp}
The following conditions are equivalent for any finitely presented
group $G$.
\begin{enumerate}
\item[1)] Some asymptotic cone of $G$ is {\ch}.

\item[2)] Some asymptotic cone of $G$ is locally isometric to an
$\mathbb R$--tree.

\item [3)] All asymptotic cones of $G$ are $\mathbb R$--trees,
i.e., $G$ is hyperbolic.
\end{enumerate}
\end{proposition}

In contrast, in Sections 4, 5, we construct a non--hyperbolic
finitely generated group all of whose asymptotic cones are
quasi--isometric and locally isometric to an $\mathbb R$--tree.

The next theorem describes the structure of {\ch} groups. Given a
group homomorphism $\al \colon G\to H$ and a generating set $S$ of
$G$, we denote by $r_S(\alpha )$ the {\it injectivity radius of
$\alpha $ with respect to $S$}, i.e., the radius of the largest ball
$B$ in $G$ such that $\alpha $ is injective on $B$.

\begin{theorem}\label{dirlim}
Let $G$ be a finitely generated group. Then the following conditions
are equivalent.

\begin{enumerate}

\item[1)] $G$ is {\ch}.

\item[2)] There exists a scaling sequence $d=(d_n)$ such that $\CG$ is an $\mathbb R$--tree
for any non--principal ultrafilter $\omega $.

\item[3)] $G$ is the direct limit of a sequence of finitely generated
groups and epimorphisms
$$
G_1\stackrel{\al_1}\longrightarrow
G_2\stackrel{\al_2}\longrightarrow \ldots $$ such that $G_i$ is
generated by a finite set $S_i$, $\alpha_i(S_i)=S_{i+1}$, and each
$G_i$ is $\delta _i$--hyperbolic, where
$\delta _i=o(r_{S_i}(\alpha _i)).$
\end{enumerate}
\end{theorem}

\begin{proof}
We are going to show that 2) $\Rightarrow $ 1) $\Rightarrow $ 3)
$\Rightarrow $ 2). The first implication is trivial. Let us prove
the second one. Suppose that $\CG $ is an $\mathbb R$--tree for some
$\omega $ and non-decreasing $d=(d_n)$. For every $n$, we denote by
$H_n$ the group given by the presentation $\langle S\, |\mathcal
R_n\rangle $, where $\mathcal R_n$ consists of labels of all cycles
in the ball of radius $d_n$ around the identity in $\Gamma (G, S)$.
Note that $G$ is a quotient of $H_n$ and the canonical map $H_n\to
G$ is injective on the ball of radius $d_n$. It follows that the
natural epimorphisms $H_n\to H_m$ are also injective on the balls of
radius $d_n$ for arbitrary $m\ge n$.

Observe that by Kapovich-Kleiner's Theorem \ref{localglobal} and
Lemma \ref{sc} from the Appendix, there are constants $C_1$, $C_2$,
$C_3$ with the following property. Let $H$ be a group having a
finite presentation $\langle S\, |\, \mathcal R\rangle $,
$D=\max\limits_{R\in \mathcal R} |R|$. Assume that for some $\delta
$ and $d\ge \max \{ C_1\delta , D\} $, every ball of radius $C_2d$
in $\Gamma(H,S)$ is $\delta $--hyperbolic. Then $\Gamma (H, S)$ is
$C_3d$--hyperbolic.

Since $\CG $ is an $\mathbb R$--tree, balls of radius $d_n$ in
$\Gamma (G,S)$ (and hence in $\Gamma (H_n,S)$) are $o_\omega
(d_n)$--hyperbolic. Recall that any $\delta $--hyperbolic graph
endowed with the combinatorial metric becomes $1$--connected after
gluing $2$--cells along all combinatorial loops of length at most
$16\delta$ (see, for example, \cite[Ch. III.H, Lemma 2.6]{BH}).
Hence $H_n$ admits a finite presentation with generating set $S$
and relations of lengths $o_\omega (d_n)$. In particular, for
every positive integer $k$, there exists $n(k)$ such that the
sequence $(n(k))$ is strictly increasing and the following
conditions hold.
\begin{enumerate}
\item[(H1)] $H_{n(k)} $ admits a finite presentation with
generating set $S$ and relations of lengths at most
$\frac{d_{n(k)}}{C_2k}$.

\item[(H2)] Every ball of radius $\frac{d_{n(k)}}{k}$ in the
$\Gamma (H_{n(k)},S)$ is $\frac{d_{n(k)}}{C_1C_2k}$--hyperbolic.
\end{enumerate}
These conditions allow us to apply the above observation to
$H=H_{n(k)} $, $d=D=\frac{d_{n(k)}}{C_2k}$, and $\delta
=\frac{d_{n(k)}}{C_1C_2k}$. Thus $H_{n(k)} $ is
$\frac{C_3d_{n(k)}}{C_2k}$--hyperbolic. Now setting $G_k=H_{n(k)}$
we obtain a sequence of groups and homomorphisms
\begin{equation}\label{seqGk}
G_1\stackrel{\al_1}\longrightarrow
G_2\stackrel{\al_2}\longrightarrow \ldots ,
\end{equation}
where $G_k$ is $\delta _k$--hyperbolic for $\delta _k=
\frac{C_3d_{n(k)}}{C_2k}=o(d_{n(k)})$ and $\al _k$  is injective on
the ball of radius $d_{n(k)}$ as desired.

To prove 3) $\Rightarrow $ 2) we fix any sequence $d=(d_n)$ such
that
\begin{equation}\label{dchoice}
\delta _n=o(d_n),\;\;\; d_n=o(r_{S}(\al_n)).
\end{equation}
Let $\omega $ be an arbitrary non--principal ultrafilter.
According to Lemma \ref{p3} applied to the collection $\mathcal P$
of all one--element subsets of $\CG $, to show that $\CG $ is a
tree it suffices to prove that it contains no simple non--trivial
limit geodesic triangles.

Suppose that $pqs$ is a non--trivial simple triangle in $\CG $
whose sides are limit geodesics. Clearly $pqs=\lo H_n$, where
$H_n=p_na_nq_nb_ns_nc_n$ is a geodesic hexagon in $\Gamma (G, S)$
such that $p=\lo p_n$, $q=\lo q_n$, $s=\lo s_n$,
\begin{equation}\label{abcn}
|a_n|=\oom(d_{n}),\;\;\; |b_n|=\oom(d_{n}),\;\;\;
|c_n|=\oom(d_{n}),
\end{equation}
and perimeter of $H_n$ satisfies $|H_n|=\Oo(d_n)$. By
(\ref{dchoice}) we have $|H_n|=o(r_{S}(\al_n))$, hence the label of
$H_n$ represents $1$ in $G_n$ \oas. Thus $H_n$ may be considered as
a configuration in the Cayley graph of $G_n$ \oas.

Let $p$ be a non--trivial side of $pqs$. Lemma \ref{quad} implies
that $p_n$ belongs to the closed $4\delta _n$--neighborhood of the
other sides \oas. Combining this with (\ref{abcn}) and
(\ref{dchoice}), we obtain that the $\Oo(\delta_n)+\oom(d_n)$-neighborhood of $q_n\cup r_n$ contains $p_n$. Since $\Oo(\delta_n)=\oom(d_n)$, $p$ belongs
to the union of $q$ and $r$. This contradicts the assumption that the triangle
$pqr$ is simple.
\end{proof}

\begin{remark}\label{dlrem}
It is easy to see from the proof of Theorem \ref{dirlim} that (3)
can be replaced with the following (a priori stronger) condition,
which will be useful for some applications.

\bigskip

{\it (4) $G$ is the direct limit of a sequence of finitely generated
groups and epimorphisms
$$
G_1\stackrel{\al_1}\longrightarrow
G_2\stackrel{\al_2}\longrightarrow \ldots $$ such that $G_i$ is
generated by a finite set $S_i$, $\alpha_i(S_i)=S_{i+1}$, each $G_i$
is $\delta _i$--hyperbolic, $\alpha_i$ is injective on the ball of
radius $r_i$ of the group $G_i=\la S_i\ra$, where $\delta _i=o(r_i)$
and
\begin{enumerate}
 \item[(4a)] the sequence of the numbers $r_i$ is non-decreasing;
 \item[(4b)] the group $G_i$ has a presentation $\la S_i\mid
 \mathcal P_i\ra $, where $\max\limits_{P\in \mathcal P_i}|P|=o(r_i)$.
\end{enumerate}}

Indeed it is easy to see that the proof of the implication
1)$\Rightarrow $ 3) ensures (d1) for $r_i=d_{n(i)}$ and (d2) follows
from condition (H1).
\end{remark}

\begin{remark}
The third condition from the theorem implies that every lacunary
hyperbolic group (as well as any other limit of hyperbolic groups)
embeds into an ultraproduct of hyperbolic groups. Indeed let
$\phi\colon G_1\to \prod ^\omega G_i $ be the homomorphism defined
by the rule $\phi (g)=(\al _1(g), \al _2\circ\al_1 (g), \ldots )$.
If $\omega $ is non-principal, it is straightforward to see that
${\rm Ker } (\phi )= \bigcup\limits_{i=1}^\infty {\rm Ker} ({\al
_i\circ\cdots \circ \al_1}) $ and hence $\phi (G_1)\cong G$. In
particular, if a universal sentence holds in the first order group
language holds in all hyperbolic groups, then it holds in all
lacunary hyperbolic groups.

This observations is similar to Maltsev's Local Theorems
\cite{Mal}. It provides us with a uniform way of proving
(nontrivial) universal theorems for lacunary hyperbolic groups. As
an example, the reader may verify that the sentence $$\forall
x\forall y\, (x^{-1}y^2x=y^3\, \Rightarrow \, [y,x^{-1}yx]=1 )$$
is a theorem in the class of all hyperbolic group and hence in the
class of all lacunary hyperbolic groups, but it is not a theorem
in the class of all groups.
\end{remark}


\subsection{Lacunary hyperbolic groups and the classical small cancellation condition}
\label{lgatcscc}

We begin with examples of \ch groups constructed by means of the
classical small cancellation theory. Recall the small cancellation
condition $C'(\mu )$. Given a set $\mathcal R$ of words in a certain
alphabet, one says that $\mathcal R$ is {\it symmetrized} if for any
$R\in \mathcal R$, all cyclic shifts of $R^{\pm 1}$ are contained in
$\mathcal R$.

\begin{definition}\label{cl}
Let
\begin{equation}\label{pres32}
G=\langle S\mid\mathcal R\rangle
\end{equation}
be a group presentation, where $\mathcal R$ is a symmetrized set of
reduced words in a finite alphabet $S$.  A common initial subword of
any two distinct words in $\mathcal R$ is called a {\em piece}. We
say that $\mathcal R$ satisfies the $C^\prime (\mu )$ if any piece
contained (as a subword) in a word $R\in \mathcal R$ has length
smaller than $\mu |R|$.
\end{definition}

The main property of groups with $C'(\mu )$-presentations is given
by the Greendlinger Lemma below.

\begin{lemma}[{\cite[Theorem V.4.4.]{LS}}] \label{greendlin}
Let $\mathcal R$ be a symmetrized set of  words in a finite alphabet
$S$ satisfying a $C'(\mu )$ condition with $\mu \le 1/6$, $\mathcal
P=\langle S\mid\mathcal R\rangle $. Assume that a reduced \vk
diagram $\Delta$ over $\cal P$ with cyclically reduced boundary path
$q$ has at least one cell. Then $q$ and the boundary path of some
cell $\Pi$ in $\Delta$ have a common subpath $t$ with $|t|>(1-3\mu
)|\partial\Pi|$.
\end{lemma}

Given a van Kampen diagram $\Delta $ over (\ref{pres32}), we denote
by $\mathcal A (\Delta )$ the sum of the perimeters of all cells in
$\Delta$. The next lemma easily follows from the Greendlinger Lemma
by induction on the number of cells in the diagram.

\begin{lemma}\label{sum-of-per}
Suppose that a group presentation (\ref{pres32}) satisfies the
$C^\prime (\mu )$--condition for some $\mu \le 1/6$. Then for any
reduced diagram $\Delta $ over (\ref{pres32}), we have:
\begin{enumerate}
\item[(a)] $|\partial \Delta |>(1-3\mu) |\partial \Pi |> |\partial \Pi |/2$ for any cell $\Pi
$ in $\Delta $.

\item[(b)] $|\partial \Delta |>(1-6\mu )\mathcal A (\Delta )$.
\end{enumerate}
\end{lemma}

The lemma below was actually proved in \cite{Short} although it was
not stated explicitly there. The explicit statement is due to
Ollivier \cite{Oll}.

\begin{lemma}
Suppose that there exists $C>0$ such that for every minimal van
Kampen diagram $\Delta $ over (\ref{pres32}), we have $|\partial
\Delta |\ge C\mathcal A(\Delta )$. Then provided $G$ is finitely presented,
it is $\delta $--hyperbolic, where $\delta \le 12\max\limits_{R\in \mathcal R} |R|
/C^2$.
\end{lemma}

\begin{cor}\label{hypconst}
Suppose that a finite group presentation (\ref{pres32}) satisfies
the $C^\prime (\mu )$--condition for some $\mu <1/6$. Then $G$ is
$\delta $--hyperbolic, where $\delta \le 12\max\limits_{R\in
\mathcal R} |R|/(1-6\mu )^2$.
\end{cor}

\begin{definition}
We say that a subset $L\subset \mathbb N$ is {\it sparse}, if for
any $\lambda >0$, there exists a segment $I=[a,b]\subset [1, +\infty
)$ such that $I\cap L=\emptyset $ and $a/b<\lambda $.
\end{definition}

Given a (not necessary finite) presentation (\ref{pres32}), we
denote by $L(\mathcal R)$ the set $\{ |R| \mid R\in \mathcal R\}
$. The following result provides us with a rich source of examples
of {\ch} groups.

\begin{proposition}\label{sparse}
Let (\ref{pres32}) be a group presentation with finite alphabet $S$,
satisfying the $C^\prime
(\mu )$ small cancellation condition for some $\mu<1/6$. Then the
group $G$ is {\ch} if and only if the set $L(\mathcal R)$ is sparse.
\end{proposition}

\begin{proof}
Suppose that $L(\mathcal R)$ is sparse. Then for every $n\in \mathbb
N$, there exists a segment $I_n=[a_n, b_n]\subset \mathbb R$ such
that
\begin{equation}\label{In}
I_n\cap L( \mathcal R)=\emptyset,
\end{equation}
$b_n/a_n>n$, and $a_{n+1}>b_n$. We set $\mathcal R_n=\{ R\in
\mathcal R\mid |R|\le a_n\} $ and $G_n=\la S\mid\mathcal R_n\ra $.
Then $G$ is the limit of the sequence of the groups $G_n$ and the
obvious homomorphisms $\al _n\colon G_n\to G_{n+1}$. By Corollary
\ref{hypconst}, $G_n$ is $\delta _n$--hyperbolic, where $\delta
_n=O(a_n)=o(b_n)$. On the other hand, by Lemma \ref{sum-of-per} (a)
and (\ref{In}) we have $r_S(\al _n)\ge b_n/2$. Hence $G$ is lacunary
hyperbolic by Theorem \ref{dirlim}.

Now assume that $G$ is lacunary hyperbolic. Let $G=\lim G_i$ as in
Remark \ref{dlrem}. Then, in the notation of Remark \ref{dlrem}, we
have $r_i\ge \lambda (2M_i+1)$ for given $\lambda>0$ and all
sufficienly large $i$, where $M_i=\max\limits_{P\in \mathcal
P_i}|P|$. To prove that $L(\mathcal R)$ is sparse, it suffices to
show that there is no  $R\in \mathcal R$ with $r_i\ge |R| \ge
2M_i+1$. By condition (d1) in Remark \ref{dlrem}, the natural
homomorphism $G_i\to G$ is also injective on balls of radius $r_i$.
Hence $R=1$ in $G_i$. However, by Lemma \ref{sum-of-per} (a), all
words from $\mathcal P_i$ represent $1$ in the group $G_i^\prime
=\la S\mid \mathcal T_i\ra $, where $\mathcal T_i=\{ R\in \mathcal
R\mid |R|\le 2M_i\} $. Hence $R=1$ in $G_i^\prime $. As
$R\notin\mathcal T_i$, this contradicts Lemma \ref{greendlin} and
the $C^\prime (\mu )$--condition.
\end{proof}

\subsection{Relative and lacunar hyperbolicity}
\label{rach}

Recall a definition of relatively hyperbolic groups. There are at
least six equivalent definitions (the first one is due to Gromov
\cite{Gr3}). We use the definition whose equivalence to the other
definitions is proved in \cite{DS}.

\begin{definition}(\cite[Theorem 8.5]{DS})\label{defsub} Let $G$ be a
finitely generated group, $H_1,...,H_n$ subgroups of $G$. Then $G$
is called {\em (strongly) hyperbolic relative to peripheral
subgroups} $H_1,...,H_n$ if every asymptotic cone $\Con^\omega(G,d)$
of $G$ is tree-graded with respect to the collection of nonempty
ultralimits of sequences of left cosets $\lio (g_jH_i)$, $g_j\in G$,
where different sequences of cosets $(g_jH_i), (g'_jH_i)$ define the
same piece if they coincide $\omega$-almost surely.
\end{definition}

\begin{remark}\label{remsub}
Recall that in a finitely generated relatively hyperbolic group $G$,
every peripheral subgroup is finitely generated and
quasi--isometrically embedded in $G$ \cite{Os06}. This implies that
each ultralimit $\lio(g_jH_i)$ is either empty or bi--Lipschitz
equivalent to the asymptotic cone $\Con^\omega(H_i,d)$ with respect
to a finite generating set of $H_i$ (see \cite{DS} for details).
\end{remark}

The first two claims of the following proposition provides us with a
way of constructing new \ch groups from given ones. The third claim
is related to the following problem. It is well--known that if a
group $G$ is hyperbolic relative to a finitely presented subgroup
$H$, then $G$ is finitely presented itself. However it is still
unknown (see \cite[Problem 5.1]{Os06}) whether finite presentability
of $G$ implies finite presentability of $H$ (although $H$ is
finitely generated whenever $G$ is \cite{Os06}). An ultimate
negation of this implication would be the following statement:
\begin{quote}
Any finitely generated recursively presented group $H$ embeds into a
finitely presented relatively hyperbolic group $G$ as a peripheral
subgroup.
\end{quote}

Propositions \ref{sparse} and \ref{thsub2} (c) imply that this
statement does not hold.

\begin{proposition} \label{thsub2}
\begin{itemize}
\item[(a)] If a finitely generated group $G$ is hyperbolic
relative to a \ch subgroup $H$, then $G$ is itself \ch.

\item[(b)] Every \ch group $H$ embeds into a 2--generated \ch group
$G$. Moreover one can assume that $G$ is hyperbolic relative to $H$.

\item[(c)] If $G$ is hyperbolic relative to a \ch subgroup $H$ and $H$ is
not finitely presented, then $G$ is not finitely presented.
\end{itemize}
\end{proposition}

\proof (a) Suppose that $G$ is hyperbolic relative to $H$ and
$\Con^\omega(H,d)$ is an $\R$-tree. By Definition \ref{defsub} and
Remark \ref{remsub}, $\Con^\omega(G,d)$ is tree-graded relative to a
collection of $\mathbb R$--trees. In particular, $\Con^\omega(G,d)$
has no nontrivial simple loops, i.e., it is an $\R $--tree.

(b) By \cite[Theorem 1.1]{AMO} applied to the free group of rank
$2$, there exists a 2-generated group $G$ such that $H$ embeds in
$G$ and $G$ is hyperbolic relative to $H$. It remains to use part
(a).

(c) The group $G$ is \ch by (a). If $G$ was finitely presented, it
would be hyperbolic by Kapovich-Kleiner's Theorem \ref{main} from
the Appendix. Since peripheral subgroups are quasi--isometrically
embedded into relatively hyperbolic groups \cite{Os06} and
quasi--isometrically embedded subgroups of hyperbolic groups are
hyperbolic \cite{Short}, $H$ is hyperbolic. Hence $H$ is finitely
presented that contradicts our assumption.
\endproof

Observe that the first assertion of Proposition \ref{thsub2} cannot
be generalized to the case of several peripheral subgroups.
Moreover, we have the following.

\begin{example}\label{315}
The free product $H_1\ast H_2$ of \ch groups is not necessarily
\ch. Indeed it is not hard to construct a set of words $\mathcal
R=\{ R_i, i\in \mathbb N\} $ in a finite alphabet $S$ such that
$|R_i|=i$ and $\mathcal R$ satisfies the $C^\prime (1/7)$
condition. It is also easy to find two subsets $N_1, N_2\subset
\mathbb N$ such that both $N_1, N_2$ are sparse and $N_1\cup N_2
=\mathbb N$. Set
$$
H_1=\la S\mid R_i=1,\, i\in N_1\ra ,\;\;\; H_2=\la S\mid R_i=1,\,
i\in N_2\ra .
$$
Then by Proposition \ref{sparse} $H_1, H_2$ are \ch while $H_1\ast
H_2$ is not.
\end{example}

\subsection{Subgroups of \ch groups}
\label{twor}

The next theorem shows that subgroups of \ch groups share common
properties with subgroups of hyperbolic groups.

\begin{definition}Let $H$ be a subgroup of a finitely generated group $G=\la S\ra$.
Then the {\em growth function} of $H$ (relative to $S$) is the
function
$$f_{H,G}(n)=\#(\Ball_G(n)\cap H)$$
where $\Ball_G(n)$ is the ball of radius $n$ around $1$ in the group
$G$ (in the word metric related to $S$. We say that $H$ has {\em
exponential growth in $G$} if its growth function $f_{H,G}(n)$ is
bounded from below by an exponent $d^n$ for some $d>1$.
\end{definition}

\begin{theorem}\label{sub} Let $G$ be a \ch group. Then

\begin{itemize}
\item[(a)] Every finitely presented subgroup of $G$ is a subgroup of
a hyperbolic group.

\item[(b)] Every undistorted (i.e. quasi-isometrically embedded)
subgroup of $G$ is \ch.

\item[(c)] Let $H$ be a (not necessarily finitely generated) bounded torsion
subgroup of $G$. Then the growth of $H$ in $G$ is not exponential.
\end{itemize}

\end{theorem}

\proof (a) Indeed, every finitely presented subgroup of a lacunary
hyperbolic group $G$ that is a direct limit of hyperbolic groups
$G_i$ as in Theorem \ref{dirlim} is isomorphic to a subgroup of one
of the $G_i$'s.

\begin{remark} It is worth noting that for the same reason, every finitely
presented subgroup of the free Burnside group $B$ of any
sufficiently large odd exponent is cyclic because $B$ is a direct
limit of hyperbolic groups \cite{Ivanov}, periodic subgroups of
hyperbolic groups are finite \cite{Gr}, and finite subgroups of $B$
are cyclic \cite{Ad}.
\end{remark}

\medskip

(b) If $H$ is a finitely generated undistorted subgroup of a \ch
group $G$ then every asymptotic cone $\Con^\omega(H,(d_n))$ of $H$
is bi-Lipschitz homeomorphic to a subspace of the corresponding
asymptotic cone $\Con^\omega(G, (d_n))$. Since a connected subspace
of an $\R$-tree is an $\R$-tree, every undistorted subgroup of a \ch
group is \ch itself.

\medskip

(c) Let $G=\la S\ra$ be a direct limit of hyperbolic groups $G_i=\la
S_i\ra$ and homomorphisms $\alpha_i$ as in Remark \ref{dlrem},
$S=S\iv, S_i=S_i\iv$ Note that the volume of the ball
$\Ball_{G_i}(r)$ of radius $r$ in $G_i$ is at most $a^r$ for $a=\#
S+1$.

Suppose $H\le G$ has exponential growth in $G$ and bounded
torsion: $h^n=1$ for all $h\in H$. Denote by $H_i$ the preimage of
$H$ in $G_i$. Since $H$ has exponential growth in $G$, there is
$d>1$ such that the number of elements from $H_i$ of length at most
$t$ in the generators $S_i$ is at least $d^t$ for every integer $t$.
We denote this subset by $B_i(t)$.

Suppose for some $i$ and $t_0$, all the elements of $B_i(2t_0)$ have
finite orders. Then all elements of both $B_i(t_0)$ and
$B_i(t_0)B_i(t_0)$ have finite orders. It is proved in \cite[Lemma
17]{IvOl} that under these assumptions the subset $B_i(t_0)$ must be conjugate
in the hyperbolic group $G_i$ to
a subset of $\Ball_{G_i}(205\delta_i)$. Therefore $d^{t_0}\le \card
B_i(t_0)\le a^{205\delta_i}$, and so $t_0\le C\delta_i$, for $C= 205
\log_d a$. Hence there is a constant  $D\ge 2C$ such that in every
intersection $H_i\cap \Ball_{G_i}(D\delta_i)$, there exists an
element $g_i$ of infinite order.

Since $H$ has exponent $\le n$, the epimorphism $G_i\to G$ is not
injective on the ball of radius $|g_i^{n}| \le nD\delta_i$ for every
$i$. Hence in the notation of Remark \ref{dlrem},  we have
$r_i<nD\delta_i$. Therefore $
 \frac{\delta_i}{r_i}>\frac{1}{nD}
$ for all $i$ and (c) is proved by contradiction.
\endproof

\begin{remark} Note that both the growth condition and the bounded
torsion condition are essential in Theorem \ref{sub}(c). Indeed,
Theorem \ref{cext1} gives examples of \ch groups with infinite
(central) subgroups of any given finite exponent $\ge 2$; Theorem
\ref{am} below provides an example of \ch group with infinite
locally finite (torsion) normal subgroup, and \ch groups from
Theorem \ref{ExoticQuotients}, part 1, and Theorem \ref{th2} are
torsion themselves.
\end{remark}

\begin{cor} \label{corsub} (a) A \ch group $G$ cannot contain copies of $\Z^2$ or
Baumslag-Solitar groups.

(b) The free Burnside group $B(m,n)$ with sufficiently large
exponent $n$ and $m\ge 2$ cannot be a subgroup of a \ch group.

(c) The lamplighter group $(\Z/2\Z) \mathrm{wr} \Z$ cannot be a
subgroup of a \ch group.
\end{cor}

\begin{proof}
(a) It follows from Theorem \ref{sub}, part (a). Indeed $\Z^2$ and
Baumslag-Solitar groups are finitely presented and are not subgroups
of hyperbolic groups

\medskip
(b) Indeed, $B(m,n)$ has exponential growth for any $m\ge 2$, $n\gg
1$ \cite{Ivanov}. Hence it has exponential growth in any group $G$
containing $B(m,n)$ and we can use Theorem \ref{sub} (c).

\medskip

(c) Consider any short exact sequence $$1\to N\to G\to Q\to 1$$
where $G=\la S\ra$, $Q=\la SN\ra$. Let $f_{N,G}$ be the growth
function of $N$ in $G$, $f_G$ and $f_Q$ be the growth functions of
$G$ and $Q$ relative to the generating sets $S$ and $SN$
respectively. Then the following inequality obviously holds:
$$f_G(n)\le f_{N,G}(2n)f_Q(n)$$
for every $n\ge 1$. Applying this inequality to the short exact
sequence $$1\to N\to (\Z/2\Z) \mathrm{wr} \Z\to \Z\to 1$$ where
$N=(\Z/2\Z)^{\Z}$ is the base of the wreath product, we deduce
that $N$ has exponential growth in the lamplighter group, and,
consequently, in any finitely generated group containing the
lamplighter group. Since $N$ has exponent $2$, it cannot be a
subgroup of a \ch group by Theorem \ref{sub} (c).

\begin{remark} It is clear that the same argument shows that a \ch group has
no finitely generated subgroup $H$ of exponential growth which is an
extension of a bounded torsion subgroup by a nilpotent group.
\end{remark}
\end{proof}


\subsection{Lacunary hyperbolic amenable groups}
\label{chag}

Recall that the class of elementary amenable groups is defined to be
the smallest class containing all finite and Abelian groups and
closed under taking directed unions, extensions, quotients, and
subgroups. In this section, we will construct an elementary amenable
group some of whose asymptotic cones are $\R$-trees.

Pick a prime number $p$ and a non-decreasing sequence ${\bf c}$ of
positive integers $c_1\le c_2\le \ldots $. Consider the group
$A=A(p,{\bf c})$ generated by $a_i, i\in \Z$ subject to the
following relations:
$$a_i^p=1, i\in\Z,$$ $$[...[a_{i_0},a_{i_1}],...,a_{i_{c_n}}]=1$$ for every
$n$ and all commutators with $\max\limits_{j,k} |i_j-i_k|\le n$. The
group $A=A(p,{\bf c})$ is locally nilpotent since arbitrary $a_j,
a_{j+1},\dots, a_{j+n}$ generate a nilpotent subgroup of nilpotency
class at most $c_n$. Since the locally nilpotent group $A$ is
generated by elements of order $p$, it is a $p$-group \cite{Hall}.
Notice that for $l\ge 0$, there is a retraction $\pi_l$ of $A$ onto
the finite subgroup $A(l)$ generated by $a_0,\dots,a_l$
($\pi_l(a_j)=1$ if $j\ne 0,\dots l$).

The group $A$ admits the automorphism $a_i\to a_{i+1}$ ($i\in\mathbb
Z$). Denote by $G=G(p,{\bf c})$ the extension of $A$ by this
automorphism, i.e., $G$ is generated by the normal subgroup $A$ and
an element $t$ of infinite order such that $ta_it^{-1}=a_{i+1}$ for
every integer $i$.

\begin{lemma}\label{newgr} The group $G=G(p,{\bf c})$ satisfies
the following properties:
\begin{itemize}
\item[(a)] $G$ is 2-generated,
\item[(b)] $G$ is (locally nilpotent $p$-group)-by-(infinite cyclic),
and so it is elementary amenable,

\item[(c)] $G$ admits an epimorphism onto the wreath product
$(\Z/p\Z) wr \Z$ and so $G$ is not virtually cyclic.
\item[(d)] $G$ is residually (finite $p$-group).
\end{itemize}
\end{lemma}

\proof (a) It is clear that $G$ is generated by $a_0$ and $t$.

(b) and (c) follow form the construction.

(d) Fix $m=p^s$ for some $s>0$ and consider a ``circular version"
$B$ of the group $A$, namely the (finite) group $B$ generated by $m$
elements $b_{[i]}$, where $[i]$ is a residue class modulo $m$, with
defining relations $b_{[i]}^p=1$ and
$[...[b_{[i_0]},b_{[i_1]}],...,b_{[i_{c_n}]}]=1$ for every $n$ and
all commutators with $\max\limits_{j,k} |[i_j-i_k]|\le n$, where
$|[i]|$ is the smallest non-negative integer $s$ such that either
$i-s$ or $i+s$ is $0$ modulo $m$. The group $B$ is a finite
$p$-group.

Notice that if $0\le l\le m/2$, then $B$ has a retraction on the
subgroup $B(l)$ generated by $b_{[0]},\dots, b_{[l]}$ which are
subject to all the defining relations of $B$ involving these
generators. Since $l\le m/2$, we have $|[i_j-i_k]|=|i_j-i_k|$ for
$i_j, i_k\in \{0,\dots,l$, and so the group $B(l)$ is naturally
isomorphic with $A(l)$.

We have an epimorphism $\alpha_m: A\to B=B_m$ such that
$\alpha_m(a_i)=b_{[i]}$ for all $a_i$-s, and, by the previous
observation, this homomorphism maps the subgroup $A(l)\le A$
isomorphically onto the subgroup $B(l)\le B$ provided $l\le m/2$.

Let $H_m$ be the extension of $B=B_m$ by the automorphism of order
$m$ : $t_m^{-1}b_{[i]}t_m=b_{[i+1]}, t_m^m=1$. We have
$|H_m|=m|B_m|$, so $H_m$ is a $p$-group. Then the epimorphism
$\alpha_m$ extends to the epimorphism $\beta_m$ of $G$ onto $H_m$
such that $\beta_m(t)=t_m$. The intersection of kernels of arbitrary
infinite family of homomorphisms $\beta_m$ contains no non-trivial
elements from $A(l)$, $l=0,1, 2,\dots$, and so it is trivial. Hence
$G$ is a residually (finite $p$-group).
\endproof

\begin{lemma}\label{amenab}
The groups $G(p, {\bf c})$ are limits of hyperbolic (in fact
virtually free) groups satisfying all assumptions of Theorem
\ref{dirlim}, provided the sequence ${\bf c}$ grows fast enough.
\end{lemma}
\begin{proof}

We chose a sequence ${\bf c}= (c_1, c_2,\dots )$ by induction.

Let $C_0=\ast_{i=-\infty}^\infty \, \langle a_i\, |\, a_i^p=1\rangle
$, $c_1=1$. Suppose that the integers $c_1\le \dots\le c_n$ are
already chosen. For every $l\le n$, we denote by $U_l$ the normal
subgroup of $C_{0}$ generated (as a normal subgroup) by all
commutators of the form
$$[...[a_{i_0},a_{i_1}],...,a_{i_{c_l}}],$$ where $\max\limits_{j,k}
|i_j-i_k|\le l$. Then we set $V_n=\prod_{l=1}^{n}U_l$ and
$C_n=C_0/V_n$.

The map $a_i\to a_{i+1}$, $i=0,1,\ldots $, extends to an
automorphism of $C_n$. Let
$$G_n=\langle C_n, t\, |\, a_i^t=a_{i+1}, \, i\in \Z\rangle
$$
be the corresponding extension by the automorphism. Clearly $G_n$ is
generated by $\{ a_0, t\} $, and the set of defining relations of
$G_{n-1}$ is a subset of the set of defining relations of $G_n$ (and
the set of defining relations of $G$). Thus the identity map on this
set induces epimorphisms $$G_1\to G_2\to...\to G_{n-1}\to G_n\to
G.$$ Observe that $G_n$ splits as an HNN--extension of its subgroup
generated by $a_0, \ldots , a_n$. This subgroup is nilpotent of
class at most $c_n$ and generated by elements of order $p$, hence it
is a finite $p$-group \cite{Hall}. This implies that $G_n$ is
virtually free and so the Cayley graph of $G_n$ corresponding to the
generators $t, a_0$ is $\delta _n$-hyperbolic for some $\delta_n$.

Now we are going to explain that choosing $c_{n+1}\gg c_n$ we can
always ensure the condition $\delta _{n}=o(t_{n})$ from Theorem
\ref{dirlim}. By the results of Higman \cite{Hig},  $C_n$ is
residually (finite $p$-group), and so it is residually nilpotent.

Now observe that for every $c=c_{n+1}$, the image of
$U_{n+1}=U_{n+1}(c)$ in $C_n$ belongs to the $(c+1)$st term $\gamma
_{c+1}(C_{n})$ of the lower central series of $C_{n}$ and since
$C_{n}$ is residually nilpotent we can make the natural homomorphism
$C_{n}\to C_{n}/\gamma _{c+1}(C_{n})$ (and, hence, the homomorphism
$C_{n}\to C_{n+1}$) injective on any given finite subset by choosing
big enough $c$. Hence the homomorphism $G_{n}\to G_{n+1}$ can be
made injective on the ball of radius $t_n=\exp (\delta _n)$, for
example.

Finally we note that the set of relations of the direct limit group
$G$ coincides with the set of relation of the group $G(p,{\bf c})$
which is not virtually cyclic by Lemma \ref{newgr}, part (c). This
completes the proof.
\end{proof}

Combining Theorem \ref{dirlim} and Lemma \ref{amenab}, we obtain the
following.

\begin{theorem}\label{am}
There exists a finitely generated elementary amenable group $G$ and
a scaling sequence $d=(d_n)$ such that $G$ is not virtually cyclic
and for any ultrafilter $\omega $, the asymptotic cone $\CG $ is an
$\mathbb R$--tree.
\end{theorem}

Note that $G(p,{\bf c})$ from Lemma \ref{newgr} is clearly not
finitely presented (because it is not virtually free and is a direct
limit of virtually free groups). Hence it has a non-simply connected
asymptotic cone \cite{Drutu}. Thus we obtain the following
corollary.

\begin{cor}
There is a finitely generated elementary amenable group having at
least two non-homeomorphic asymptotic cones.
\end{cor}


\section{Graded small cancellation and circle-tree asymptotic cones}
\label{cscactac}

Recall the definition of circle-trees.

\begin{definition}\label{ct}
We say that a metric space $X$ is a {\it circle--tree}, if $X$ is
tree graded with respect to a collection of circles (with the
standard length metric) whose radii are uniformly bounded from below
and from above by positive constants. In particular, every
circle--tree is locally isometric to an $\mathbb R$--tree.
\end{definition}

Note that by Lemma \ref{tree-like}, any group having a circle-tree
asymptotic cone is \ch. In this Section, we shall show that the
class of groups all of whose asymptotic cones are trees or
circle-trees is very large and contains all groups given by
presentations satisfying certain small cancellation conditions. As
a consequence, this class contains groups all of whose proper
subgroups are cyclic or finite.



\subsection{The Greendlinger
Lemma for small cancellation presentations over hyperbolic
groups}\label{CSC1}


Let $H$ be a group generated by a finite set $S$. We will consider
quotient groups of $H$ as groups given by presentations over $H$
(i.e. presentations including all relations of $H$ plus some extra
relations). Our goal is to generalize Definition \ref{cl} for such
presentations.

We start with a definition of a piece. In what follows we write
$U\equiv V$ for two words $U$ and $V$ is some alphabet to express
letter--by--letter equality.

\begin{definition}\label{piece}
Let $H$ be a group generated by a set $S$. Let $\mathcal R$ be a
symmetrized set of reduced words in $S^{\pm 1}$. For $\e
>0$, a subword $U$ of a word $R\in \mathcal R$ is called a {\it
$\e $--piece}  if there exists a word $R^\prime \in \mathcal R$
such that:

\begin{enumerate}
\item[(1)] $R\equiv UV$, $R^\prime \equiv U^\prime V^\prime $, for
some $V, U^\prime , V^\prime $; \item[(2)] $U^\prime = YUZ$ in $H$
for some words $Y,Z$ such that $\max \{ |Y|, \,|Z|\} \le \e $;
\item[(3)] $YRY^{-1}\ne R^\prime $ in the group $H$.
\end{enumerate}
Note that if $U$ is an $\e$-piece, then $U'$ is an $\e$-piece as
well.
\end{definition}

Recall that a word $W$ in the alphabet $S^{\pm 1}$ is called {\it
$(\lambda , c)$--quasi--geodesic} (respectively {\it geodesic}) in
$H$ if any path in $\Gamma (H, S)$ labeled by $W$ is $(\lambda ,
c)$--quasi--geodesic (respectively geodesic).

\begin{definition}\label{SC} Let $\e \ge 0$, $\mu\in (0,1)$, and $\rho
>0$.
We say that a symmetrized set $\mathcal R$ of words over the
alphabet $S^{\pm 1}$ satisfies the {\em condition $C(\e , \mu ,\rho
)$} for the group $H$, if
\begin{enumerate}
\item[($C_1$)] All words from $\mathcal R$ are geodesic in $H$;
\item[($C_2$)] $|R|\ge \rho $ for any $R\in \mathcal R$;
\item[($C_3$)] The length of any $\e$-piece contained in any word $R\in \mathcal R$ is smaller than $\mu |R|$.
\end{enumerate}
\end{definition}

Suppose now that $H$ is a group defined by
\begin{equation}\label{H}
H=\langle S\, |\, \mathcal O\rangle ,
\end{equation}
where $\mathcal O$ is the set of all relators (not only defining)
of $H$. Given a symmetrized set of words $\mathcal R$, we consider
the quotient group
\begin{equation}\label{quot}
H_1=\langle H\, |\, \mathcal R\rangle = \langle S\, |\, \mathcal
O\cup \mathcal R\rangle .
\end{equation}
A cell in a van Kampen diagram over (\ref{quot}) is called an {\it
$\mathcal R$--cell} (respectively, an {\it $\mathcal O$--cell}) if
its boundary label is a word from $\mathcal R$ (respectively,
$\mathcal O$). We always consider van Kampen diagrams over
(\ref{quot}) up to some natural elementary transformations. For
example we do not distinguish diagrams if one can be obtained from
the other by joining two distinct $\mathcal O$--cells having a
common edge or by the inverse transformation, etc. (see
\cite[Section 5]{Ols} for details).

\begin{figure}
\vspace{1mm}
    \hspace{15mm}\includegraphics[width=95mm]{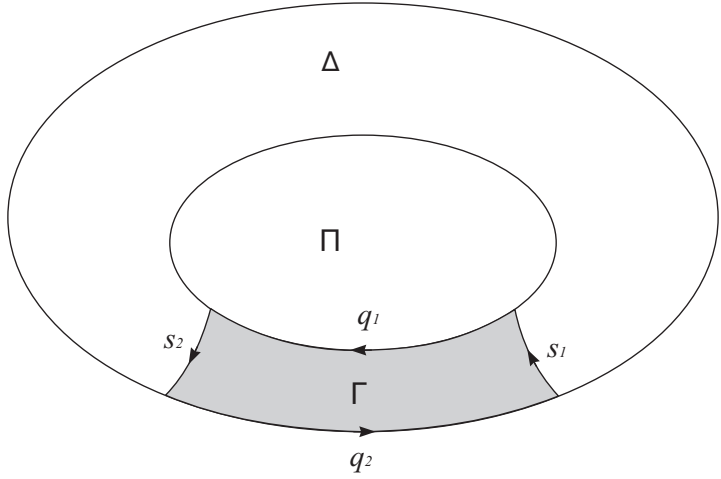}\\
    \vspace{-5mm}
  \caption{A contiguity subdiagram}\label{cont-fig}
\end{figure}

Let $\Delta $ be a van Kampen diagram over (\ref{quot}), $q$ a
subpath of its boundary $\partial \Delta $, $\Pi $, $\Pi ^\prime $ some $\mathcal
R$--cells of $\Delta $. Suppose that there is a simple closed path
$p=s_1q_1s_2q_2$ in $\Delta $, where $q_1$ (respectively $q_2$) is a
subpath of the boundary $\partial \Pi $ (respectively $q$ or $\partial \Pi
^\prime $) and $\max \{ |s_1|,\, |s_2|\} \le \e $ for some constant
$\e $. By $\Gamma $ we denote the subdiagram of $\Delta $ bounded by
$p$. If $\Gamma $ contains no $\mathcal R$--cells, we say that
$\Gamma $ is an {\it $\e $--contiguity subdiagram} of $\Pi $ to the
part $q$ of $\partial \Delta $ or $\Pi ^\prime $ respectively (Fig.
\ref{cont-fig}). The subpaths $q_1$ and $q_2$ are called {\it
contiguity arcs} of $\Gamma $ and the ratio $|q_1|/|\partial \Pi |$
is called the {\it contiguity degree} of $\Pi $ to $\partial \Delta
$ (or or $\Pi ^\prime $) and is denoted by $(\Pi , \Gamma ,
\partial \Delta )$ (or $(\Pi , \Gamma ,
\Pi^\prime )$).

The following easy observation will often be useful.

\begin{lemma}\label{geodpart}
Suppose that the group $H$ is hyperbolic. Let $\mathcal R$ be a set
of geodesic in $H$ words, $\Delta $ a diagram over (\ref{quot}), and
$q$ a subpath of $\partial \Delta $ whose label is geodesic in
$H_1$. Then for any $\e\ge 0$, no $\mathcal R$--cell $\Pi $ in
$\Delta $ have an $\e $--contiguity subdiagram $\Gamma $ to $q$ such
that $(\Pi , \Gamma , q)>1/2 +2\e/|\partial \Pi |$.
\end{lemma}

\begin{proof}
Let $\Gamma $ be an $\e$--contiguity subdiagram of an $\mathcal
R$--cell $\Pi $ to $q$, $\partial \Gamma =s_1q_1s_2q_2$, where $q_2$
is a subpath of $q$ and $\partial \Pi =q_1r$. Let also $\gamma $
denote the contiguity degree $(\Pi, \Gamma , q)$. Then we have
$|q_2|\ge |q_1|-|s_1|-|s_2|= \gamma |\partial \Pi | - 2\e $ since
$\Lab (q_1)$ is geodesic in $H$. On the other hand, $|q_2| \le
|r|+|s_1|+|s_2|\le (1-\gamma )|\partial \Pi | +2\e$ as $\Lab (q_2)$
is geodesic in $H$. These two inequalities yield $\gamma \le 1/2
+2\e/|\partial \Pi |$.
\end{proof}

Given a van Kampen diagram $\Delta $ over (\ref{quot}), we call a
combinatorial map from the $1$-skeleton $\Sk^{(1)} (\Delta )$ to
the Cayley graph $\Gamma (H_1, S)$ {\it natural} if it preserves
labels and orientation of edges. The following easy observation
will be useful.

\begin{lemma}\label{easyobs} Suppose that $H$ is hyperbolic.
Let $\mathcal R$ be a symmetrized set of words in $S^{\pm 1}$
satisfying the condition $C(\e , \mu, \rho )$ for some $\e \ge 0$,
$\mu\in (0,1)$, $\rho >0$. Suppose that $\Pi $, $\Pi ^\prime $ are
two $\mathcal R$--cells in a diagram $\Delta $ over (\ref{quot}) and
$\Gamma $ is an $\e $--contiguity subdiagram of $\Pi ^\prime $ to
$\Pi $ such that $\gamma =(\Pi ^\prime , \Gamma , \Pi)\ge \mu $.
Then
\begin{equation}\label{raznitsa}
||\partial \Pi | - |\partial \Pi ^\prime ||\le 2 \e
\end{equation}
and for any natural map $\phi\colon \Sk^{(1)} (\Delta )\to \Gamma
(H_1, S)$, the Hausdorff distance between $\phi (\partial \Pi )$
and $\phi (\partial \Pi ^\prime )$ does not exceed $2\e +2\delta $.
\end{lemma}

\begin{proof}
Let $\partial \Pi =uv$, $(\partial \Pi ^\prime )^{-1}=u^\prime
v^\prime $, and $\partial\Gamma =yuz(u^\prime)^{-1}$. Since $ \Lab
(yuz)=\Lab (u^{\prime })$ in $H$, we have
\begin{equation}\label{yuv}
\Lab (y)\Lab (uv)\Lab (y)^{-1}=\Lab (u^\prime v^\prime )
\end{equation}
in $H$ by the $C(\e , \mu, \rho )$--condition. As labels of
$\partial \Pi $ and $\partial \Pi ^\prime $ are geodesic in $H$,
(\ref{yuv}) implies (\ref{raznitsa}).

Further let $Q=abcd$ be a quadrangle in $\Gamma (H, S)$ such that
$a$, $b$, $c$, and $d$ are labeled by $\Lab (y)$, $\Lab (uv)$, $\Lab
(y)^{-1}$, and $(\Lab (u^\prime v^\prime ))^{-1}$, respectively. By
the first assertion of Lemma \ref{quad}, $b$ and $d$ belong to the
closed $(\e +2\delta)$--neighborhoods of each other. To finish the
proof it remains to note that the map $\Gamma (H, S)\to \Gamma (H_1,
S)$ induced by the homomorphism $H\to H_1$ does not increase the
distance.
\end{proof}

We call a (disc) \vk diagram over (\ref{quot}) {\em minimal} if it
has minimal number of $\mathcal R$--cells among all disc diagrams
with the same boundary label.  The first part of the following
result is an analog of the Greendlinger Lemma \ref{greendlin} for
presentations over hyperbolic groups.

\begin{lemma}\label{gamma-cell}
Suppose that $H$ is a $\delta$-hyperbolic group having presentation
$\la S\, |\, \mathcal O\rangle$ as in (\ref{H}), $\varepsilon\ge
2\delta$, $0<\mu\le 0.01$, and $\rho$ is large enough (it suffices
to choose $\rho > 10^6\varepsilon /\mu$). Let $H_1$ be given by a
presentation
$$
H_1=\langle H\, |\, \mathcal R\rangle = \langle S\, |\, \mathcal
O\cup \mathcal R\rangle $$ as in (\ref{quot}) where $\mathcal R$ is
a finite symmetrized set of words in $S^{\pm 1}$ satisfying the
$C(\e , \mu, \rho )$--condition. Then the following statements hold.

\begin{enumerate}
\item  Let
$\Delta$ be a minimal disc diagram over (\ref{quot}). Suppose that
$\partial\Delta=q^1\cdots q^t$, where the labels of $q^1, \ldots ,
q^t$ are geodesic in $H$ and $t\le 12$. Then, provided $\Delta$ has
an $\mathcal R$-cell, there exists an $\mathcal R$-cell $\Pi$ in
$\Delta$ and disjoint $\varepsilon$-contiguity subdiagrams
$\Gamma_1, \dots, \Gamma_t$ (some of them may be absent) of $\Pi$ to
$q^1,\ldots, q^t$ respectively such that
$$(\Pi,\Gamma_1,q^1)+\dots+(\Pi,\Gamma_t,q^t)>1-23\mu .$$

\item $H_1$ is a $\delta_1$-hyperbolic group with $\delta_1\le 4r$
where $r=\max\{|R| \mid R\in {\cal R}\}$.
\end{enumerate}
\end{lemma}

\proof The first part of the lemma is essentially a special case of
\cite[Lemma 6.6]{Ols} where the parameters $\lambda$ and $c$ of the
quasi-geodesity of the defining words are equal to $1$ and $0$,
respectively because the words in $\cal R$ are geodesic in $H$. The
minor corrections in the argument of \cite{Ols} leading to this
special case are the following.

We replace $t\le 4$ by $t\le 12$. Note that the proof from
\cite{Ols} works even in the case $t\le k$, for any fixed $k$, but
then we should replace 23 by $C=C(k)$ in the formulation of the
lemma.

As in the proof from \cite{Ols}, we need to consider geodesic
quadrangles $s_1q_1s_2q_2$ with ``short" sides $s_1, s_2$ and
``long" sides $q_1, q_2$. Then if a point $o\in q_1$ is far from the
ends of $q_1$, say, $\min (\dist(o,q_-), \dist (o,q_+))\ge
\max(|s_1|,|s_2|)+2\delta$ then in \cite{Ols}, it is proved that the
distance from $o$ to $q_2$ is bounded from above by $c_1=13\delta$.
In our case, by Lemma \ref{quad}, we can take $c_1=2\delta$.

Finally, in \cite[Lemma 6.2]{Ols}, we can replace the upper estimate
$n\sqrt\rho$ by $n\mu\rho$ (this is possible because $\rho$ is large
enough).

The proof of the second statement of the lemma is divided into two
steps.

{\em Step 1.} First we consider a minimal diagram $\Gamma$ over
$H_1$ whose contour contains a subpath $pq$ where segments $p$ and
$q$ are geodesic in $\Gamma$. Assume that there is an $\cal R$-cell
$\Pi$ with two contiguity subdiagrams $\Gamma_p$ and $\Gamma_q$ to
$p$ and $q$, respectively, and with contiguity arcs $v_p\subseteq p$
and $v_q\subseteq q$, such that $(v_p)_-=p_-$ and $(v_q)_+=q_+$. We
claim that the Hausdorff distance between the images $\bar p$ and
$\bar q$ of $p$ and $q$ in the Cayley graph of $H_1$ does not exceed
$2r$. We may assume that the boundary path of $\Gamma$ is the
product of $pq$, an arc of $\Pi$, and two side arcs of $\Gamma_p$
and $\Gamma_q$ (see Fig. \ref{fig1})

\begin{figure}
 \hspace{10mm}
 \includegraphics{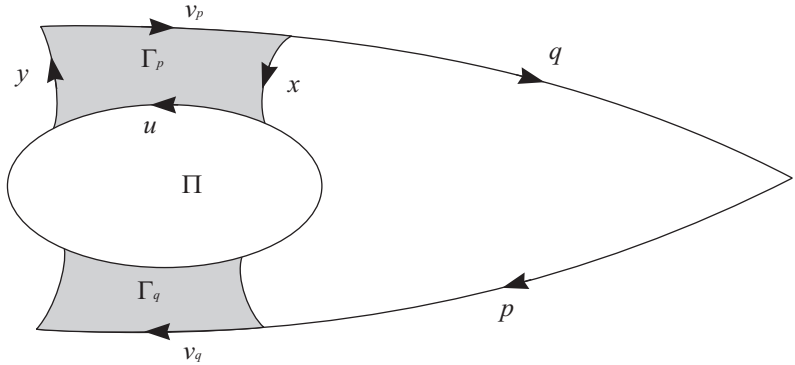}\\
 \caption{The diagram $\Gamma $.}
\label{fig1}
\end{figure}

 Let $v_pxuy$ be the boundary path of $\Gamma_p$ where $u$ is the
 contiguity arc of $\Gamma_p$ lying on $\partial\Pi$. Since $v_p$ is geodesic we have
 $|v_p|\le |x|+|y|+|u|< 2\varepsilon +r$. Therefore, using $\Pi$ and the
 contiguity subdiagrams, we can connect every point of $v_p$ (of $v_q$)
 to $q$ (to $p$) with a path of length $\le |v_p|/2 + 2\e +r <2r$
 since, according to the first statement of the lemma, $\rho$ is chosen
 much greater than $\varepsilon$.

Then we consider a maximal set of $\cal R$-cells $\Pi,
\Pi',\Pi'',\dots$ in $\Gamma$ having disjoint contiguity subdiagrams
$\Gamma_p, \Gamma_q, \Gamma'_p, \Gamma'_q,\dots$ to both $p$ and
$q$. After suitable enumeration, they provide us with decompositions
$p=v_pw_pv'_pw'_p\dots $ and $q=\dots w'_qv'_qw_qv_q$ where $v_p,
v'_p, v''_p,\dots$ and $v_q,v'_q, v''_q,\dots$ are the contiguity
arcs of the above contiguity subdiagrams.  As in the previous
paragraph, we have that the distance between every point of $v'_p,
v''_p,\dots$ (of $v'_q, v''_q,\dots$) and $q$ (and $p$) is less than
$2r$. Thus it suffices to obtain the same estimate for the distance
between a point of one of $w_p, w'_p,\dots$ and $q$. (More
precisely, it suffices to do this with the images $\bar w_p,\dots,
\bar q$ of these paths in the Cayley graph of $H_1$.)

For example, $w'_p$ is a section of a loop $w'_p abc w'_qdef$
where $b$ and $e$ are arcs of $\Pi'$ and $\Pi''$, respectively,
and $a,c,d,f$ are geodesics of length at most $\varepsilon$.
Denote by $\Xi$ the subdiagram bounded by this octagon. If $\Xi$
is a diagram over $H$,  then by Lemma \ref{quad}, every point of
$\bar w'_p$ is at distance at most $6\delta$ from the union of the
remaining $7$ sides. Since $\max(|b|,|e|)\le r$, we have that the
distance between a point of $w'_p$ and $w'_q$ is at most
$6\delta+2\e+r<2r$ since the parameter $\rho$ is chosen so that
$\rho >6\delta+2\e$.

Thus to complete Step 1, it suffices to show that $\Xi$ contains no
$\cal R$-cells. Arguing by contradiction, we have an $\cal R$-cell
$\pi$ and its contiguity subdiagrams $\Gamma_1,\dots,\Gamma_8$ (some
of them may be absent) to $w'_p, a,\dots, f$, respectively, with
$$(\pi,\Gamma_1,w'_p)+\dots+(\pi,\Gamma_8,f)>1-23\mu $$ by the first
assertion of the lemma.  Here the contiguity degree to $b$ and $e$
are less than $\mu$ by the $C(\e,\mu,\rho)$-condition. Since the
lengths of $a,c,d,f$ are less than $\e\ll\rho$, the contiguity
degree of $\pi$ to each  of these four boundary sections of $\Xi$ is
less than $\mu/2$. (The accurate proof of the latter inequality is
given in \cite[Lemma 6.5(a)]{Ols}.) Hence
$(\pi,\Gamma_1,w'_p)+(\pi,\Gamma_5,w'_q)>1-27\mu .$

If, for example, $\Gamma_5$ is absent, then $w'_p$ is homotopic in
$\Gamma$ to the path $x'zy'$, where $x'$ and $y'$ are side arcs of
contiguity subdiagram and $z$ an arc on the boundary $\partial\pi$
with length $<27\mu|\partial\pi|$. Since $w'_p$ is geodesic, we have
$|w'_p|<2\e+27\mu|\partial\pi|$. On the other hand, $w'_p$ is
homotopic to $x'u'y'$ in the diagram $\Gamma_p$ over $H$ where $u'$
is the arc of $\pi$ of length at least $(1-27\mu)|\partial\pi|$.
Therefore $|w'_p|+2\e > (1-27\mu)|\partial\pi|$. Since
$|\partial\pi|\ge\rho$ we obtain
 $(1-54\mu)\rho<4\e$ that contradicts $\rho >
10^6\mu^{-1}\varepsilon$.

Thus both $\Gamma_1$ and $\Gamma_5$ are present, and so the cell
$\pi$ can be added to the set $\Pi,\Pi',\dots$ contrary the
maximality. This contradiction completes Step 1.

 {\em Step 2.} Assume that $\delta_1> 4r$.
Then, by Rips' definitions of hyperbolicity (see \cite{Gr}, 6.6)
there exists a geodesic triangle $xyz$ in the Cayley graph of $H_1$
such that

\begin{quote}
($\star$) for arbitrary three points $o_1,o_2,o_3$ chosen on the
sides $x$, $y$, and $z$, respectively, we have
$$\max(\dist(o_1,o_2),\dist(o_2,o_3) ,\dist(o_3,o_1))\ge\delta_1 >
4r.$$
\end{quote}

Let $\Delta$ be a minimal diagram over $H_1$ corresponding to the
triangle. We preserve the notation $xyz$ for the boundary of
$\Delta$. A subdiagram $\Gamma$ is said to be an $xy$-{\em corner}
of $\Delta$ if for some subpaths $p$ of $x$ and $q$
 of $y$ such that $p_+=q_-$ (i.e., $pq$ is a subpath of $xy$), it contains a cell $\Pi$ with two
 contiguity subdiagrams $\Gamma_p$ and $\Gamma_q$ satisfying the conditions of Step 1,
 and it is bounded by $p$, $q$, $\Pi$, $\Gamma_p$ and $\Gamma_q$ as in the first paragraph
 of Step 1. The $xy$-corner is called maximal, if the sum $|p|+|q|$ is maximal. By
 definition, it consists of the single vertex $x_+=y_-$ if there exists no $\Pi$ as above.
 Similarly we define $yz$- and $zx$-corners.

    In this notation, we suppose  that the $xy$-corner $\Gamma=\Gamma_{xy}$ is maximal, $x$ contains the
    subsegment $p$, and let $p'$ be a similar segment of a maximal $zx$-corner $\Gamma_{zx}$ where $p'$         also lies on $x$.
    If $p$ and $p'$ together cover $x$, then every point of $x$ belongs to the $2r$-neighborhood
    of the union of two other sides of the triangle $xyz$. This implies that there is a point
    $o_1$ on $x$ and two points $o_2$ and $o_3$ on $y$ and $z$, respectively, such that
    $\max (\dist(o_1,o_2),\dist(o_1,o_3))\le 2r$, and this contradicts Condition $(\star)$.
    Hence $p$ and $p'$ must be disjoint. It follows that all three maximal corners
    $\Gamma_{xy},\Gamma_{yz}$ and $\Gamma_{zx}$ can be chosen pairwise
    disjoint(see Fig. \ref{fig2}).

\begin{figure}
\hspace{0mm} \includegraphics[width=124mm]{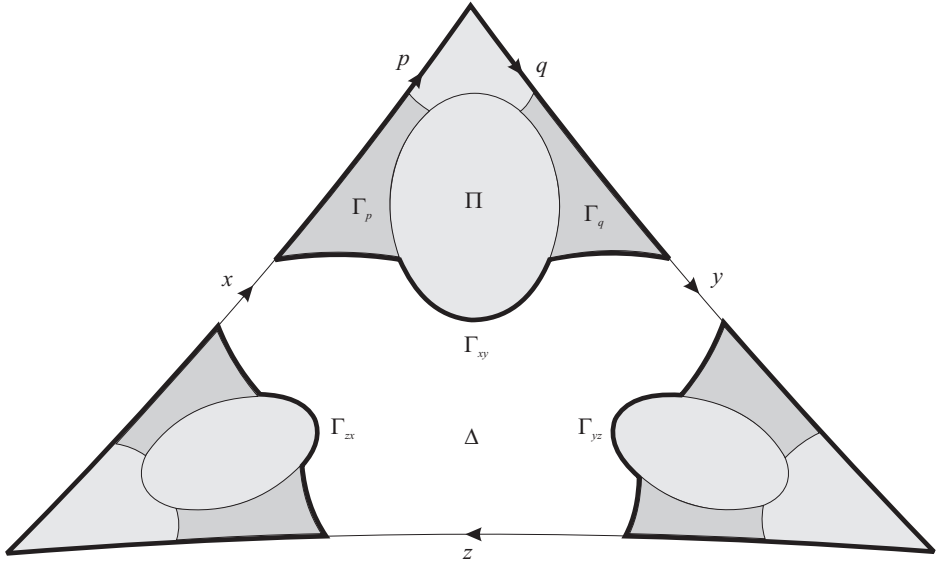}\\
 \caption{The corners $\Gamma_{xy},\Gamma_{yz}$, and $\Gamma_{zx}$ in $\Delta $.}
\label{fig2}
\end{figure}

    Now consider the diagram $\Delta'$ obtained from $\Delta$ by cutting off these
    three corners. The contour of $\Delta'$ is of the form $x'ay'bz'c$ where $x',y'$, and $z'$
    belong in $x,y$, and $z$, respectively. In turn, $a= a_1a_2a_3$ where $a_2$ is an
    arc of the cell $\Pi$ from the definition of corner, and $a_1$, $a_3$ are side arcs
    of $\Gamma_p$ and $\Gamma_q$. Similarly, $b=b_1b_2b_3$ and $c=c_1c_2c_3$. Thus
    we have a decomposition of $\partial\Delta'$ in 12 sections. Note that
    $\sum_{i=1}^3 (|a_i|+|b_i|+|c_i|)\le 6\e + 3r$.

    If $\Delta'$ is a diagram over $H$, then by Lemma \ref{quad}, $x'$ belongs to the
    $10\delta$-neighborhood of the union of the remaining 11 sides of the $12$-gon,
    and so $x'$ belongs to the $(10\delta+6\e+3r)$-neighborhood of $y\cup z$.
    Since $10\delta+6\e<\rho\le r$, every point of $x'$ belongs to the
    $4r$-neighborhood of $y\cup z$. Recall that every point of $x\backslash x'$
    is in $p\cup p'$ and belongs in the $2r$-neighborhood of $y\cup z$. Similarly,
    every point of $y$ or $z$ belongs in $4r$-neighborhood of the union of two
    other sides, which contradicts to the choice of the triangle $xyz$.

    Then we may assume that $\Delta'$ has an $\cal R$-cell, and so, by the
    first assertion of the lemma, it has a cell $\pi$ with contiguity subdiagrams
    $\Gamma_{x'},\Gamma_{a_1},\dots$ to the 12 sections, such that the sum of
    their contiguity degrees is greater than $1-23\mu$. As at the Step 1, the sum
    of contiguity degrees for all $a_i, b_i,c_i$ ($i=1,2,3$) is less than $6\mu$.
    So the sum of the contiguity degrees to $x', y'$, and $z'$ is greater than
    $1-29\mu$. It is impossible that two of these three subdiagrams be absent
    because the degree of contiguity of $\pi$ to a geodesic section cannot
    be as high as $1-29\mu$ by the argument of Step 1 with  coefficient $27$
    replaced by $29$.

    Therefore we may assume that there are contiguity subdiagrams of $\pi$ to
    both $x'$ and $y'$. But this contradicts the maximality of the corner
    $\Gamma_{xy}$, and the lemma is completely proved.
\endproof


\subsection{Comparing various small cancellation conditions}
\label{cvscc}

The purpose of this section is to compare the small cancellation
condition introduced in Section \ref{CSC1} with the classical small
cancellation condition.

\begin{lemma} \label{S1S2}
Let $R_1, R_2, \ldots $ be an infinite set of distinct cyclically
reduced words in a finite alphabet $S$, $\mathcal R_i$ the set of
all cyclic shifts of $R_i^{\pm 1}$, and ${\cal
R}=\bigcup_{i=1}^\infty \mathcal R_i$. Assume that
\begin{enumerate}
\item[($SC_1$)] The set $\cal R$  satisfies $C^\prime (\lambda)$
for some $\lambda\le 1/10$.

\item[($SC_2$)] The set $\mathcal R_n$ ($n=1,2,\dots$) satisfies
$C^\prime (\lambda_n)$ where $\lambda_n\to 0$.
\end{enumerate}
Let also $H = \langle{S} | R_{i_1},\dots,R_{i_j}\rangle$ for some
$i_1<\dots<i_j$. Then, for every $\mu>0, \varepsilon \ge 0$, and
$\rho>0$, there is $n>i_j$ such that the set $\mathcal R_n$
satisfies the $C(\varepsilon,\mu, \rho)$-condition over $H$.
\end{lemma}

\proof We first prove that any subword $V$ of a cyclic shift of
$R_n^{\pm 1}$ is geodesic in the group $H$.  Let $U$ be a geodesic
word equal to $V$ in $H$.  Notice that $V$ has no subword $W$
which is also a subword of
 one of the words $R$
from $\mathcal T=\mathcal R_{i_1}\cup  \dots \cup \mathcal R_{i_j}$
 with $|V|\ge \lambda|R|$ by condition ($SC_1$).  Then by Lemma
 \ref{greendlin}, either $U\equiv V$, and we are done, or there
 is a subword $W$ of $U$ which is also a prefix
 of  a word $R$ from $\mathcal T$ with $|W|>
 (1-3\lambda-2\lambda)|R|\ge |R|/2$. Therefore the subword $W$ is
 equal in $H$ to a shorter word $T^{-1}$, where $R\equiv WT=1$ in
 $H$. But then $U$ is not geodesic, a contradiction.

 Now, it remains to show
that a word from $\mathcal R_n $ has no $\varepsilon$-pieces of
length at least $\mu|R_n|$ if $n$ is large enough.  We may further
assume by ($SC_2$) that  $\lambda_n<\mu/2$ and $|R_n|>
\mu^{-1}\max(|R_{i_j}|,4\varepsilon)$. (Observe that ($SC_2$)
implies $|R_n|\to \infty $ and $n\to \infty $.)

 Assume that $R, R'\in \mathcal R_n$, $R\equiv UV, R'\equiv U'V'$,
$|U|\ge\mu|R_n|$,
 and $U'=YUZ$ in $H$ for some geodesic in $H$ words $Y, Z$ with
 $\max(|Y|,|Z|)\le\varepsilon$. We must prove that
 $YRY^{-1}=R'$ in $H$.

 Consider a diagram $\Delta$ over $H$ with boundary path
 $p_1q_1p_2q_2$, where $\Lab (p_1)\equiv Y$, $\Lab (q_1)\equiv U$, $\Lab
 (q_2)\equiv Z$, and $\Lab (p_2)\equiv (U')^{-1}$.

 If $\Delta$ has no cells then there is a common subpath $t$ of $q_1$
 and $q_2^{-1}$ with $|t|\ge\lambda_n|R_n|$ because $$|p_1|+|p_2|\le
 2\varepsilon<\mu|R_n|/2< (\mu-\lambda_n)|R_n|,$$ and every edge
 of $q_1$ belong to the path $p_2q_2p_1$.
 From the small cancellation condition ($SC_1$), we see that
 $T\equiv \Lab(t)$ is a prefix of a unique word $R_0$ from $\mathcal R_n$.
 Since $U\equiv X_1TX_2$ and $U'\equiv X_3TX_4$ for some words
 $X_1,X_2, X_3, X_4$, we obtain that the conjugates $X_1^{-1}RX_1$
 and $X_3^{-1}R'X_3$ are both freely equal to $R_0$, i.e.
 $R=X_1X_3^{-1}R'(X_1X_3^{-1})^{-1}$. Since $YX_1X_3^{-1}$ is a
 label of a closed path in $\Delta$, we have
 $YRY^{-1}=(YX_1X_3^{-1})R'(YX_1X_3^{-1})^{-1} = R'$ in $H$, as
 required.

 Now assume by contradiction that $\Delta$ has a cell. No cell $\Pi$
 has a boundary arc with length $>|\partial\Pi|/2$  contained in
 $p_1$ or in $p_2$ because these segments are geodesic. Since $q_1$ is a
 geodesic path by (1), no cell $\Pi$ can have vertices in both
 $p_1$ and $p_2$ because otherwise
 $$\mu|R_n|\le |q_1|\le|p_1|+|\partial\Pi|/2+|p_2|\le 2\varepsilon
 +|R_{i_j}|/2\le\max(|R_{i_j}|,4\varepsilon)$$ contrary the choice of $n$.
Then it follows from
 Lemma \ref{greendlin}, that there is a cell $\Pi$ in $\Delta$,
 whose boundary has a common subpath $t$ with either $q_1$ or
 $q_2$, where $|t|>\frac12 (1-3\lambda-1/2)|\partial\Pi|\ge
 \lambda|\partial\Pi|$. But this contradicts condition ($SC_1$) because
 $n\notin \{i_1,\dots,i_j\}$.
\endproof

\begin{remark}\label{exist}
Infinite set{s} of words in the alphabet $\{ a, b\} $ satisfying
($SC_1$) and ($SC_2$) were constructed in various places (see, for
example, \cite{EO}). Moreover, one can find such sets satisfying
($SC_1$) for arbitrary small $\lambda >0$.
\end{remark}

Our next goal is to relate $C(\e, \mu , \rho)$ to the following
small cancellation condition $C(\e , \mu , \lambda , c, \rho )$
introduced in \cite{Ols}. Essentially, this condition is obtained
from $C(\e, \mu, \rho)$ by replacing the word ``geodesic" by
``quasi-geodesic"; the additional parameters are the quasi-geodesic
parameters. We shall show that any $C(\e , \mu , \lambda , c, \rho )$--presentation
can be transformed into a $C(\e', \mu', \rho')$--presentation of the same group for suitable $\e', \mu', \rho'$.

\begin{definition}
Let $H$ be a hyperbolic group generated by a finite set $S$. A
symmetrized set $\mathcal R$ of words in $S^{\pm 1}$ satisfies the
{\it condition $C(\e , \mu , \lambda , c, \rho )$} for some
constants $\e >0$, $\mu\in (0,1)$, $\lambda \in (0,1]$, $c\ge 0$,
$\rho >0$, if all words in $\mathcal R$ are $(\lambda ,
c)$--quasi--geodesic and conditions ($C_2$), ($C_3$) from
Definition \ref{SC} hold. Thus $C(\e, \mu, \rho )$ is equivalent
to $C(\e, 1,0,\mu, \rho)$.
\end{definition}

In the two lemmas below the following notation is used. Let
$\mathcal R_0 $ {be} a set of words in $S^{\pm 1}$, $\mathcal R_0
^\prime $ the set obtained from $\mathcal R_0$ by replacing each
$R\in \mathcal R_0$ with a shortest word $R^\prime $ such that $R$
and $R^\prime $ are conjugate in $H$. Denote by $\mathcal R$
(respectively $\mathcal R^\prime$) the set of all cyclic shifts of
words from $\mathcal R_0 ^{\pm 1}$ (respectively $(\mathcal
R_0^\prime) ^{\pm 1})$. Note that all words in $\mathcal R^\prime $
are geodesic in $H$.

Given $\lambda >0, c\ge 0$, we define $$\kappa =8\delta +4\theta ,$$
where $\delta $ is the hyperbolicity constant of the {Cayley graph}
$\Gamma (H, S)$, and $\theta =\theta (\lambda ,c)$ is the constant
from Lemma \ref{qg}.

\begin{lemma}\label{conj}
Suppose that all words in $\mathcal R$ are $(\lambda ,
c)$--quasi-geodesic in $H$ for some $\lambda , c$ and have lengths
at least $(\kappa+c)/\lambda $.  Then for every $W\in \mathcal
R^\prime $, there is $R_W\in \mathcal R$ such that $W$ and $R_W$ are
conjugate by an element of length at most $\kappa /2$ in $H$.
\end{lemma}

\begin{proof}
Any word $W\in \mathcal R^\prime $ is conjugate to some word $R\in
\mathcal R^{\pm 1}$ in $H$. Let $T_{W,R}$ denote a shortest word conjugating
$W$ to $R$ in $H$. Let also $U$ and $S$ be cyclic shifts of $W$ and $R$ respectively such that
\begin{equation}\label{TUS}
|T(U,S)|\le |T(U^\prime , S^\prime )|
\end{equation}
for any cyclic shifts $U^\prime $ and $S^\prime $ of $W$ and $R$.

There is a 4-gon $asbu^{-1}$ in $\Gamma (H,S)$ such that $\Lab (u)\equiv U$, $\Lab (s)\equiv  S$, $\Lab
(a)\equiv \Lab (b^{-1})\equiv T(U,S)$. Clearly $s$ is $(\lambda , c)$--quasi--geodesic, and
$a,u,b$ are geodesic. Moreover
\begin{equation}\label{dus}
\dist (u,s)\ge |a|
\end{equation}
by (\ref{TUS}). There are two cases to consider.

{\it Case 1.} If $|a|\le 2\delta +\theta $, then by Lemma \ref{qg} and Lemma
\ref{quad}, $u$ and $s$ belong to the closed $(4\delta
+2\theta)$--neighborhoods of each other. Therefore, $W$ is
conjugate to a cyclic shift of $R$ by a word of length at most
$4\delta +2\theta $.

{\it Case 2.} Now assume that $|a|> 2\delta +\theta $.
In particular, $\dist (u,s)> 2\delta +\theta$  by (\ref{dus}).
Consider the middle point $m$ of $s$. Lemmas \ref{qg} and \ref{quad}
imply that $\dist (m, a\cup b)\le 2\delta +\theta
$. For definiteness, assume that $\dist (m, a)\le 2\delta +\theta
$. Let $z$ be the point on $a$ such that $\dist (m,a)=\dist
(m,z)$. Then
$$
\dist (z, a_-)\ge \dist (m, a_-)-2\delta -\theta \ge \dist (s,
u)-2\delta -\theta \ge |a|-2\delta -\theta
$$
by (\ref{dus}). Therefore, $$\dist (z, s_-)\le |a|-\dist (z, a_-)\le 2\delta
+\theta .$$ Consequently, $$\dist (m, s_-) \le \dist (m, z)+\dist
(z, s_-)\le 2(2\delta +\theta ).$$ This means that a cyclic shift
of $R$ represents an element of length at most $4(2\delta +\theta
)$ in $H$. Hence $|R|\le \lambda^{-1} (8\delta +4\theta+c)$ that
contradicts our assumption. Hence this case in impossible.
\end{proof}

\begin{lemma}\label{Chyp}
Suppose that $\mathcal R$ satisfies $C(\e , \mu , \lambda , c,
\rho )$ for some $\e >0$, $\mu\in (0,1)$, $\lambda \in (0,1]$,
$c>0$, and
\begin{equation}\label{ro}
\rho>2(c+3\kappa )/\lambda .
\end{equation}
Then $\mathcal R^\prime $ satisfies $C(\e^\prime, \mu^\prime,
\rho^\prime )$ for
\begin{equation}\label{newpar}
\e^\prime = \e -2\kappa ,\;\;\; \mu^\prime = 2\mu/\lambda,\;\;\;
\rho ^\prime =\lambda \rho -c -\kappa .
\end{equation}
\end{lemma}

\begin{proof}
By (\ref{ro}) and Lemma \ref{conj} for any word $W \in \mathcal
R^\prime$, there is a word $R_W\in \mathcal R$ that is conjugate
to $W$ by a word of length at most $\kappa /2$ in $H$. In
particular, this yields the last inequality in (\ref{newpar}).

Suppose now that for some words $W_1, W_2$ of $\mathcal R^\prime
$, we have $W_1\equiv U_1V_1$, $W_2\equiv U_2V_2$, and $U_1=YU_2Z$
in $H$, where $|Y|, |Z|\le \e^\prime $ and $|U_1|\ge \mu ^\prime
|W_1|$. Let $A_1,A_2$ be words of lengths at most $\kappa /2$ such
that
\begin{equation}\label{rwi}
R_{W_i}=A_iW_iA_i^{-1},\; i=1,2
\end{equation}
in $H$. Using Lemma \ref{qg} and Lemma \ref{quad}, we can find
initial subwords $C_i$ of $R_{W_i}$, $i=1,2$, such that
\begin{equation}\label{ci}
C_i=A_iU_iB_i,\; i=1,2
\end{equation}
in $H$, where
\begin{equation}\label{bi}
|B_i|\le \theta +2\delta +\kappa/2<\kappa.
\end{equation}
Thus $C_1=A_1YA_2^{-1}C_2B_2^{-1}ZB_1$ in $H$. Note that $$\max\{
|B_2^{-1}ZB_1|, |A_1YA_2^{-1}|\} \le \e^\prime +2\kappa =\e .$$
Now using subsequently (\ref{ci}), (\ref{bi}), (\ref{rwi}), and
(\ref{ro}) we obtain
$$
\begin{array}{rl}
|C_1|\ge & |U_1|-|A_1|-|B_1|\ge \mu^\prime |W_1|-2\kappa\ge
\mu^\prime (\lambda |R_{W_1}|-c-\kappa) -2\kappa \\ &
\\{\ge} & 2\mu
|R_{W_1}|-2\mu (c+\kappa )/\lambda -2\kappa\ge \mu |R_{W_1}|.
\end{array}
$$
Hence
$$
R_{W_1}= (A_1YA_2^{-1}) R_{W_2} (A_1YA_2^{-1})^{-1}
$$
in $H$ by the $C(\e, \lambda, c,\mu,\rho)$--condition. Combining
this with (\ref{rwi}), we obtain $W_1=YW_2Y^{-1}$.
\end{proof}


\subsection{Groups with circle-tree asymptotic cones}
\label{gwctac}

In this section we define the graded small cancellation condition
and prove that all asymptotic cones of any group given by a {graded
small cancellation} presentation are circle-trees.

In the next two sections we show that many classical small
cancellation groups as well as some `monsters' obtained by methods
from \cite{Ols} admit {graded small cancellation} presentations.

\begin{definition}\label{classQ} Let $\alpha, K$ be positive
numbers.  We say that the presentation
\begin{equation}\label{qpres}
\langle S\mid \mathcal R\rangle =\left\langle S\,\left|\,
\bigcup\limits_{i=0}^\infty \mathcal R_i\right.\right\ra
\end{equation}
of a group $G$ is a $\qq(\alpha,K)$-{\em presentation} if the
following conditions hold for some sequences $\e =(\e_n)$, $\mu
=(\mu _n)$, and $\rho =(\rho _n)$ of positive real numbers
($n=1,2\dots $).
\begin{enumerate}

\item[($ \mathbf Q_0$)] The group $G_0=\langle S\mid \mathcal R_0\rangle $
is $\delta_0$-hyperbolic for some $\delta_0$.

\item[($\mathbf Q_1$)] For every $n\ge 1$, $\calr_{n}$ satisfies $C(\e_n,
\mu _n, \rho _n)$ over $G_{n-1}=\left\langle S\, \left|\,\right.
\bigcup\limits_{i=0}^{n-1} \mathcal R_i\right\rangle .$

\item[($\mathbf Q_2$)] $\mu_n=o(1)$, $\mu _n\le
\alpha$, and $\mu_n\rho_n>K\e_n$ for any $n\ge 1$.

\item[($\mathbf Q_3$)] $\e_{n+1}>8\max\{|R|, R\in \calr_n\}=O(\rho_n)$.
\end{enumerate}
\end{definition}

The following lemma shows that if $\alpha$ is small enough and $K$
is big enough, $\qq(\alpha,K)$-presentations have properties
resembling the properties of ordinary small cancellation
presentations.

\begin{lemma}\label{asc}
Let (\ref{qpres}) be a $\qq(.01, 10^6)$--presentation. Then the
following conditions hold.
\begin{enumerate}

\item[(a)] For every $n\ge 1$, Lemma \ref{gamma-cell} applies to
$H=G_{n-1}$ and $H_1=G_{n}=\langle H\mid \mathcal R_{n}\rangle $.
In particular, $G_n$ is $\delta _n$--hyperbolic, where $\delta
_n\le 4\max\limits_{R\in \mathcal R_{n}}|R|$.

\item[(b)] $\e _n=o(\rho _n)$.

\item[(c)] $\rho _{n} =o(\rho _{n+1})$; in particular, $\rho _n\to
\infty $ as $n\to \infty $ and $\delta _n=o(\rho _{n+1})$.

\item[(d)] $\rho _n=o(r_S(G_{n}\to G_{n+1}))$, where $r_S$ is the
injectivity radius.
\end{enumerate}
\end{lemma}
\begin{proof}
The first assertion easily follows from ($\mathbf Q_0$)--($\mathbf Q_2$) by
induction. Assertions (b) and (c) follow immediately from ($\mathbf Q_2$)
and ($\mathbf Q_3$). Finally given $g\in {\mathrm{Ker}}\, (G_n\to G_{n+1})$,
$g\ne 1$, we consider a geodesic (in $G_n$) word $W$ representing
$g$ and a minimal van Kampen diagram over $G_{n+1}$ with boundary
label $W$. Applying Lemma \ref{gamma-cell} for $r=1$ and taking into
account that words in $\mathcal R_{n+1}$ are geodesic in $G_n$, we
obtain $|W|>(1-23\mu _{n+1})|\partial \Pi |-2\e_n$, where $\Pi $ is
an $\mathcal R_{n+1}$--cell provided by Lemma \ref{gamma-cell}.
Hence $|g|\ge (1-o(1))\rho_{n+1}$. Combining this with (c) we obtain
(d).
\end{proof}

\begin{definition}
From now {on} the condition $\qq=\qq(.01, 10^6)$ will be called the
{\em graded small cancellation} condition.
\end{definition}

The following statement is an immediate corollary of Lemma
\ref{S1S2}. It shows, in particular, that the class of groups
admitting graded small cancellation presentations is large.

\begin{cor}\label{subset}
Let $R_1, R_2, \ldots $ be an infinite set of distinct cyclically
reduced words in a finite alphabet $S$, $\mathcal R_i$ the set of
all cyclic shifts of $R_i^{\pm 1}$, and ${\cal
R}=\bigcup_{i=1}^\infty \mathcal R_i$. Assume that conditions
$(SC_1)$ for $\lambda<\frac1{100}$ and $(SC_2)$ of Lemma \ref{S1S2} are satisfied.  Then there
is an infinite sequence $i_1<i_2<\ldots $ such that $\left\langle
S\,\left|\, \bigcup_{j=1}^\infty \mathcal
R_{i_j}\right.\right\rangle $ satisfies the graded small
cancellation condition.
\end{cor}

\begin{definition}\label{od-vis}
Given an ultrafilter $\omega $ and a scaling sequence $d=(d_n)$, we
say that a sequence of real numbers $f=(f_n)$ is {\it $(\omega ,
d)$--visible} if there exists a subsequence $(f_{n_i})$ of $f$ such
that $f_{n_i}=\To (d_i)$.
\end{definition}

\begin{theorem}\label{classG7}
For any group $G$ having a \gsc pre\-senta\-tion, any ultrafilter
$\omega$, and any sequence of scaling constants $d=(d_n)$, the
asymptotic cone $\CG$ is a circle--tree. $\CG $ is an $\mathbb
R$--tree if and only if the sequence $(\rho_n)$ from Definition
\ref{classQ} is not $(\omega , d)$--visible.
\end{theorem}

The proof of the theorem is divided into a sequence of lemmas.
Throughout the rest of the section we fix arbitrary scaling sequence
$d=(d_n)$ and ultrafilter $\omega $.

Let us fix a group $G$ having a graded small cancellation
presentation (\ref{qpres}). In what follows, we denote by $\ds $
(respectively $\dh $) the distance (respectively the Hausdorff
distance) in the Cayley graph {$\Gamma (G, S)$.}


\begin{lemma}\label{Rn2}
Any subword of any word $R_n\in \mathcal R_n$ of length at most
$|R_n|/2$ is $(1-\oom(1),0)$--quasi--geodesic.
\end{lemma}
\begin{proof}
Suppose $U_n=V_n$ in $G$, where $U_n$ is a subword of $R_n\in
\mathcal R_n$, $|U_n|\le |R_n|/2$, and $V_n$ is a geodesic word in
$G$. By Lemma \ref{asc}, $U_n=V_n$ in $G_{n}$ \oas. Let $\Delta _n$
be a minimal diagram over $G_{n}$ with boundary $p_nq_n$, where
$\Lab (p_n)=U_n$, $\Lab (q_n)=V^{-1}_n$. If $\Delta _n$ has no
$\mathcal R_n$--cells, then $|U_n|=|V_n|$ since $R_n$ is geodesic in
$G_{n-1}$. Thus we may assume that $\Delta _n$ has at least one
$\mathcal R_n$--cell. Let $\Pi _n$ be the $\mathcal R_n$--cell,
$\Gamma ^1_n$, $\Gamma ^2_n$ the contiguity diagrams to $p$, $q$,
respectively, provided by Lemma \ref{gamma-cell}.

Observe that $\lio (\Pi_n, \Gamma ^2_n, q_n)\le 1/2$. (Indeed
otherwise $(\Pi_n, \Gamma ^2_n, q_n)>1/2 +2\e_n/|\Pi_n|$ {\oas }
as $2\e_n/|\Pi_n|=o(1)$ and $V_n$ could not be geodesic by Lemma
\ref{geodpart}.) Thus $\lio (\Pi_n, \Gamma ^1_n, p_n)\ge 1/2>\mu
_n$. Hence $||\partial \Pi _n|-|R_n||\le 2\e_n$ by Lemma
\ref{easyobs}. Now there are two cases to consider.

{\it Case 1.} If $\lio (\Pi _n, \Gamma _n^2, q_n)=1/2$, then
$$
\begin{array}{rl}
|V_n|\ge & (1/2-\oom (1))|\partial \Pi _n|-2\e_n
\\& \\ {\ge} & (1/2-\oom (1))(|R_n|-2\e_n)-2\e_n\ge (1
-\oom(1))|U_n|
\end{array}
$$
and we are done.

{\it Case 2.} If $\lio (\Pi _n, \Gamma _n^2, q_n)<\theta<1/2$,
then $\lio (\Pi _n, \Gamma _n^1, p_n) >1-\theta $ and we get a
contradiction as
$$
|U_n|>(1-\theta )|\partial \Pi _n|-2\e_n \ge (1-\theta )
(|R_n|-2\e_n)-2\e_n> \frac12|R_n|
$$
\oas.
\end{proof}

\begin{definition} \label{defC}
Let $\mathcal A$ denote the set of
all loops in the Cayley graph $\Gamma (G, S)$ labeled by words from
the set of relators $\calr$ (see (\ref{qpres})). We say that a sequence $(p_n)$ of elements of $\mathcal A$ is {\it
asymptotically visible} (relative to the scaling sequence $d$ and ultrafilter $\omega $) if
\begin{equation}\label{vis}
|p_n|=\To (d_n), \;\;\; {\rm and }\;\;\; \ds (1, p_n)=\Oo (d_n).
\end{equation}
\noindent By $\C =\C (d, \omega )$ we denote the collection of all
distinct limits $\lo p_n$, where $(p_n)$ ranges in the set of all
asymptotically visible sequences of elements of the set $\mathcal A$.
\end{definition}

\begin{lemma}\label{circles}
Every piece from $C=\lo p_n\in \C $ is isometric to a circle of
length $\lio \frac{|p_n|}{d_n}$.
\end{lemma}

\begin{proof}
By Lemma \ref{Rn2}, $p_n$ equipped with the metric induced from
$\Gamma (G,S) $ is $(1-\oom(1),0)$--quasi--isometric to a circle
of length $|p_n|$ and our lemma follows. Checking details is
straightforward and we leave this to the reader.
\end{proof}

The next observation is quite trivial and follows immediately from
Property (${\bf Q}_3$) of the \gsc condition (Definition
\ref{classQ}) and Lemma \ref{asc} (c).

\begin{lemma}\label{Rn}
Suppose that the sequence $(\rho _n)$ is $(\omega , d)$--visible.
Let $({\cal R}_{i_n})$ be a subsequence of $({\cal R}_n)$, such that
$\rho _{i_n}=\To (d_n)$. Then for any asymptotically visible
sequence $(p_n)$, $\Lab(p_n)\in {\cal R}$, we have $\Lab (p_n)\in{\cal R}_{i_n}$ $\omega
$--almost surely.
\end{lemma}

\begin{lemma}\label{2points}
Suppose that for some asymptotically visible sequences $(p_n)$ and
$(q_n)$, the intersection $\lo p_n \cap \lo q_n$ contains at least
two distinct points. Then $\lo p_n=\lo q_n$.
\end{lemma}

\begin{proof}
The assumptions of the lemma imply that for every $n\in \mathbb
N$, there exist quadrangles $Q_n=p_n^\prime s_n q_n^\prime t_n$ in
$\Gamma (G,S)$ such that $p_n^\prime $ and $q_n^\prime $ are
subpaths of of $p_n$ and $q_n$, respectively, and $s_n$, $t_n$ are
geodesics such that
\begin{equation}\label{lpq}
|p_n|=\To (d_n), \;\;\; |q_n|=\To (d_n),
\end{equation}
\begin{equation}\label{lst}
|s_n|=\oom(d_n),\;\;\; |t_n|=\oom(d_n).
\end{equation}
Without loss of generality we may also assume that
\begin{equation}\label{half}
|p_n^\prime |\le \frac12 |p_n|,\;\;\; |q_n^\prime |\le \frac12
|q_n|.
\end{equation}

Let $({\cal R}_{i_n})$ be the subsequence of $({\cal R}_n)$
provided by Lemma \ref{Rn}. Thus $\Lab (p_n), \Lab (q_n)\in {\cal
R}_{i_n} $ \oas. Note that (\ref{lpq}) and (\ref{lst}) imply
$$
|Q_n|=\To (d_n)=\To (\rho _{i_n}).
$$
Therefore, by Lemma \ref{asc} (d) there is a van Kampen diagram
$\Xi _n$ with the boundary label $\Lab(Q_n)$ having no $\mathcal
R_j$--cells for $j>i_n$ \oas. For simplicity we keep the notation
$p_n^\prime $, $s_n$, $q_n^\prime $, $t_n$ for the corresponding
parts of $\partial \Xi _n$.  Let also $\Delta _n$ be the diagram
obtained from $\Xi _n$ by attaching two $\mathcal R_{i_n}$--cells
$\Sigma ^1 _n$ and $\Sigma ^2_n$ along $p_n^\prime $ and $q_n
^\prime $ respectively so that the natural map $\Sk ^{(1)}(\Delta
_n)\to \Gamma (G, S)$ sending parts of $\partial \Xi _n$ to the
corresponding sides of $Q_n$, maps $\partial \Sigma ^1 _n $ to
$p_n$ and $\partial \Sigma ^2_n$ to $q_n$. There are two cases to
consider.

{\it Case 1.} $\Xi _n$ has no $\mathcal R_{i_n}$--cells {\oas }
and hence $\Lab (Q_n)=1$ in $G_{i_n-1}$ \oas . We recall that
$\Gamma (G_{i_n-1}, S)$ is $\delta _{i_n-1}$--hyperbolic. Let
$p_n^{\prime\prime }$ be the subpath of $p_n^\prime $ such that
$\dist ((p_n^{\prime\prime })_{\pm }, (p_n^\prime )_{\pm})= \max
\{ |s_n|, \,|t_n|\} +2\delta _{i_n-1}$. (Note that $|p_n^\prime |>
2(\max \{ |s_n|, \,|t_n|\} +2\delta _{i_n-1})$ \oas.) By Lemma
\ref{quad}, $\dist _{Hau} (p_n^{\prime\prime} , q_n^\prime )\le
2\delta _{i_n-1} $. Thus we may assume that $\Delta _n$ contains a
$2\delta _{i_n-1}$--contiguity subdiagram $\Gamma _n$ of $\Sigma
^1 _n$ to $\Sigma ^2_n$ such that
$$
(\Sigma ^1 _n, \Gamma _n, \Sigma ^2 _n) = \frac{|p^\prime _n
|-2(\max \{ |s_n|, \,|t_n|\} +2\delta _{i_n-1})}{|\partial \Sigma
^1_n|}=\frac{\To (d_n)-\oom (d_n)}{\To (d_n)}>\mu _{i_n}
$$
\oas (see Fig. \ref{2p-fig}). Now Lemma \ref{easyobs} and Lemma
\ref{asc} imply that
$$\dh (p_n, q_n)\le \e_n +2\delta _{i_n-1} =\oom(d_n) .$$ Therefore
$\lo p_n=\lo q_n $.

\begin{figure}
\vspace{-1cm}\hspace{10mm}
  \includegraphics{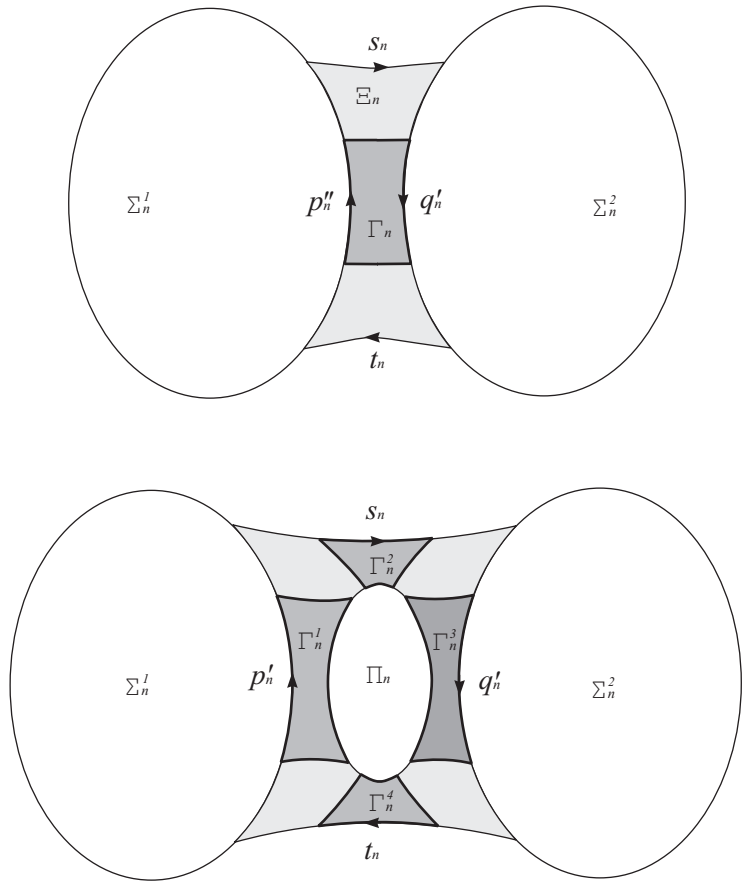}\\
  \vspace{-3mm}
  \caption{Two cases in the proof of Lemma \ref{2points}}\label{2p-fig}
\end{figure}

{\it Case 2.} Suppose now that $\Delta _n$ has at least one
$\mathcal R_{i_n}$--cell \oas . Then by Lemma \ref{gamma-cell}, we
may assume that there is an $\mathcal R_{i_n}$--cell $\Pi _n$ in
$\Xi _n$ and disjoint $\varepsilon _{i_n}$-contiguity subdiagrams
$\Gamma^1_n, \dots, \Gamma^4_n$ of $\Pi _n$ to $p_n^\prime, s_n,
q_n^\prime ,t_n$, respectively, such that
$$
(\Pi _n,\Gamma^1_n,p_n^\prime )+(\Pi _n,\Gamma^2_n,s_n)+ (\Pi
_n,\Gamma^3_n,q^\prime_n) + (\Pi _n,\Gamma^4_n,t_n) >1-23\mu
_{i_n}
$$
\oas (Fig. \ref{2p-fig}).  Using (\ref{lst}) and Lemma \ref{asc}
one can easily show that
$$
(\Pi _n,\Gamma^2_n,s_n)+ (\Pi _n,\Gamma^4_n,t_n) =\oom(1).
$$
Hence
$$
(\Pi _n, \Gamma^1_n, p_n^\prime )+(\Pi _n,\Gamma^3_n,
q^\prime_n)=1-\oom(1) .
$$
Note that $\Gamma ^3_n$ may also be considered as a contiguity
subdiagram of $\Pi_n $ to $\Sigma ^2 _n$ in $\Delta _n$.

Suppose first that $(\Pi _n,\Gamma^1_n,p_n^\prime
)=(\Pi_n,\Gamma^1_n,\Sigma ^1_n )\le\mu _{i_n}$ \oas. Then
$(\Pi_n,\Gamma^3_n,\Sigma ^2_n )=1-\oom(1)$. Applying Lemma
\ref{easyobs} we obtain
$$
||\partial\Pi_n |- |\partial \Sigma ^2_n||<2\e_n.
$$
Hence the length of the contiguity arc $v_n$ of $\Gamma _n^3$ to
$\Sigma _n^2$ satisfies
$$
|v_n|\ge (1-\oom (1))|\partial \Pi _n| -2\e_n= (1-\oom
(1))|\partial \Sigma^2 _n| .
$$
However this contradicts (\ref{half}).

Therefore $(\Pi _n,\Gamma^1_n,p_n^\prime )>\mu _{i_n}$ and
similarly $(\Pi _n,\Gamma^3_n,q_n^\prime )>\mu _{i_n}$. Let $w_n$
be the image of $\partial \Pi _n$ under the natural map $\Sk
^{(1)}(\Delta _n)\to \Gamma (G, S)$ sending $\partial \Sigma ^1_n$
and $\partial \Sigma ^2_n$ to $p_n$ and $q_n$ respectively. Then
using Lemma \ref{easyobs} and Lemma (\ref{asc}) we obtain
$$\dh (p_n,q_n)\le \dh (p_n, w_n)+\dh (q_n, w_n)=\oom(d_n).$$ Thus
$\lo (p_n)=\lo (q_n)$ again.
\end{proof}

\begin{lemma}\label{simpletr}
Every simple triangle in $\CG $ whose sides are limit geodesics is
contained in a subset from $\C $.
\end{lemma}

\begin{proof}
Suppose that $pqs$ is a simple triangle in $\CG $ whose sides are
limit geodesics. As in the proof of Theorem \ref{dirlim},
$pqs=\lo H_n$, where $H_n$ is a geodesic hexagon in $\Gamma (G,
S)$. Let $\Delta _n$ be a van Kampen diagram over (\ref{qpres})
with boundary label $\Lab (H_n)$.

First of all we assume that $(\rho _n)$ is $(\omega , d)$--visible
and denote by $({\cal R}_{i_n})$ be the subsequence of $({\cal
R}_n)$ provided by Lemma \ref{Rn}. Note that $|H_n|=\To (d_{i_n})$.
Arguing as in the proof of the previous lemma, one can show that
$\Delta _n$ has no $\mathcal R_{j}$--cells for $j>i_n$ \oas. To
simplify our notation we identify the $1$--skeleton of $\Delta _n$
with its natural image in $\Gamma (G, S)$.

If $\Delta _n$ has no $\mathcal R_{i_n}$--cells \oas, we obtain a
contradiction as in the third paragraph of the proof of
Theorem \ref{dirlim}. Thus we may assume that $\Delta _n$ has
at least one $\mathcal R_{i_n}$--cell \oas\, and $\Delta _n$ is
minimal over $G_{i_n}$. By Lemma \ref{gamma-cell} there is an
$\mathcal R_{i_n}$--cell $\Pi _n$ in $\Delta _n$ and $\e
_{i_n}$--contiguity subdiagrams $\Gamma ^1_n,\ldots , \Gamma ^6_n$
of $\Pi _n$ to the sides of $H_n$ such that $\partial \Pi _n$
belongs to the closed $(23\mu _{i_n}|\partial \Pi _n|)
$--neighborhood of the union $\bigcup_{i=1}^6 \Gamma _n^i $ in
$\Delta _n$ \oas. Let $w_n$ denote the natural image of $\partial
\Pi _n$ in $\Gamma (G, S)$. By Lemma \ref{quad} $w_n$ belongs to
the closed $(23\mu _{i_n}|\partial \Pi _n| +\e _{i_n} +2\delta
_{i_n-1})$--neighborhood of $H_n$ in $\Gamma (G,S)$. Further Lemma
\ref{asc} implies that $23\mu _{i_n}|\partial \Pi _n| +\e _{i_n}
+2\delta _{i_n-1}= \oom (d_n)$. Hence $\lio w_n\subseteq pqs$.
Since $\lio w_n $ is a circle by Lemma \ref{circles} and $pqs$ is
simple, we have $pqs=\lio w_n\in \mathcal C$.

Finally we assume that $(\rho _n)$ is not $(\omega , d)$--visible.
Let $j_n$ be the maximal number such that $\Delta _n$ has at least
one $\mathcal R_{j_n}$--cell. Let $\Pi _n$ and $\Gamma _n^1, \ldots
, \Gamma _n^6$ be the $\mathcal R_{j_n}$--cell and the $\e
_{j_n}$--contiguity subdiagrams of $\Pi _n$ to sides of $H_n$
provided by Lemma \ref{gamma-cell}. Again we can easily show that
the total length of the contiguity arcs of $\Gamma _n^1, \ldots ,
\Gamma _n^6$ to the sides of $H_n$ is $\To (\rho _{j_n})$. Since
$(\rho _n)$ is not $(\omega , d)$--visible, $\rho _{j_n}=\To (d_n)$
is impossible. Therefore $\rho _{j_n}=\oom(d_n)$. Hence
$\delta_{j_n}=\oom(d_n)$ by Lemma \ref{asc}. This again leads to a
contradiction as in the third paragraph of the proof of Theorem
\ref{dirlim}. Thus $\CG $ has no nontrivial simple triangles whose
sides are limit geodesics, i.e., it is tree--graded with respect to
the empty collection of pieces (i.e., it is an $\mathbb R$-tree) by
Lemma \ref{p3}.
\end{proof}

\begin{proof}[Proof of Theorem \ref{classG7}]
Let $G$ have a \gsc presentation (\ref{qpres}). For any fixed
scaling sequence $d=(d_n)$ and any ultrafilter $\omega $, let $\C
=\C (d, \omega )$ be the collection of subsets of $\CG $ described
in Definition \ref{defC}. If $(\rho _n)$ is $(\omega , d)$--visible,
all elements of $\C $ are circles whose radii are uniformly bounded
from below and from above by positive constants by Lemmas
\ref{circles} and \ref{Rn}. Further by Lemma \ref{2points}, $\C $
satisfies $(T_1)$ (see Definition \ref{tgspace}). Applying now Lemma
\ref{simpletr} and Lemma \ref{p3} we conclude that $\CG$ is
tree--graded with respect to $\C $, i.e., is a circle--tree. If
$(\rho _n)$ is not $(\omega , d)$--visible, the same arguments show
that $\CG$ is tree--graded with respect to $\C $. But $\C $ is empty
in this case, hence $\CG $ is an $\mathbb R$--tree.
\end{proof}


\subsection{Groups without free subgroups}
\label{gwfs}


Recall that any torsion--free non--elementary hyperbolic group has
an infinite quotient group with finite (or cyclic) proper subgroups
\cite{Ols}. We show in this section that some of these groups have
\gsc presentations. In what follows, $H$ denotes a hyperbolic group
generated by a finite set $S$.

The following Lemma \ref{old-gamma-cell} is an analog of Lemma
\ref{gamma-cell} for such presentations (and a particular case of
Lemma 6.6 \cite{Ols}), while Lemma \ref{primyk} is a
quasi-geodesic analog of Lemma \ref{geodpart}.

\begin{lemma}\label{old-gamma-cell}
For any hyperbolic group $ H$ and any $\lambda>0$, there is
$\mu_0>0$ such that for any $\mu\in (0,\mu_0]$ and any $c\ge 0$,
there are $\varepsilon\ge 0$ and $\rho>0$ with the following
property:

Let a finite symmetrized presentation $H_1 = \langle H|{\cal
R}\rangle$ satisfy the condition
$C(\varepsilon,\mu,\lambda,c,\rho)$, and $\Delta$ a minimal diagram
over $H_1$ whose boundary is a product of
$(\lambda,c)$-quasi-geodesic paths $p$ and $q$. Then provided
$\Delta$ has an $\cal R$-cell, there exists an $\cal R$-cell $\Pi$
in $\Delta$ and disjoint $\varepsilon$-contiguity subdiagrams
$\Gamma_1$ and $\Gamma_2$ (one of them may be absent) of $\Pi$ to
$p$ and $q$, respectively, such that
$(\Pi,\Gamma_1,p)+(\Pi,\Gamma_2, q)>1-23\mu$.
\end{lemma}
\endproof

\begin{lemma}\label{primyk} Let a presentation $H_1=\langle H|{\cal R}\rangle$
satisfy a $C(\varepsilon, \mu,\lambda,c,\rho)$-condition with $\mu
<\lambda^2/100 $ and $\rho> 2\mu^{-1}(c+2\varepsilon) $. Let
$\Delta$ be a minimal diagram over $H_1$ with a
$(\lambda,c)$-quasi-geodesic subpath $q$ of the boundary, and
$\Gamma$ a contiguity subdiagram of an $\cal R$-cell $\Pi$ to $q$.
Then $\psi=(\Pi,\Gamma,q)<1-24\mu$.

\proof Let $p_1q_1p_2q_2$ be the boundary of the contiguity
subdiagram of $\Pi $ to $q$, where $q_1p=\partial \Pi$ and $q_2$ is
a subpath of $q$. Since $q_2$ is a $(\lambda,c)$-quasi-geodesic
subpath of $q$, we have $\lambda|q_2|-c\le |p_2|+|p|+|p_1|\le
2\varepsilon+ (1-\psi)|\partial\Pi|$. On the other hand, $|q_2|\ge
\lambda |q_1|-c -|p_1|-|p_2|$ since the path $q_1$ is
$(\lambda,c)$-quasi-geodesic. Hence $|q_2|\ge
\lambda\psi|\partial\Pi|-c -2\varepsilon$. These two estimates for
$|q_2|$ give us the inequality
$$\lambda^{-1}(|\partial\Pi|(1-\psi)+2\varepsilon+c) \ge \psi\lambda
|\partial\Pi| -c - 2\varepsilon,$$ that is $$\psi\le
(\lambda^{-1}|\partial\Pi| +(1+\lambda^{-1})
(c+2\varepsilon))(|\partial\Pi|(\lambda^{-1}+\lambda))^{-1}.$$ Since
$|\partial\Pi|\ge \rho$, we obtain from the assumptions of the
lemma:
$$\psi \le \frac{\lambda^{-1}\rho +(1+\lambda^{-1})
(c+2\varepsilon)}{\rho(\lambda^{-1}+\lambda)}<
\frac{1+\mu}{1+100\mu}< 1-24\mu.$$
\endproof

\end{lemma}

Below we say that a bi--infinite path $p$ in the Cayley graph of a
group generated by a finite set $S$ is {\it $V$--periodic} (or just
{\it periodic}), if $p$ is labeled by the bi--infinite power of some
word $V$ in $S^{\pm 1}$.

\begin{theorem}\label{ExoticQuotients}

\begin{enumerate}
\item Let $G$ be an arbitrary non-elementary hyperbolic group.
Then there exists an infinite torsion quotient group $Q_1$ of $G$
admitting a \gsc presentation.

\item Let $G$ be an arbitrary torsion free non--cyclic hyperbolic
group with a finite set of generators $S$. Then there exists an
infinite non-Abelian torsion free quotient group $Q_2$ of $G$
admitting a \gsc presentation and such that all proper subgroups of
$Q_2$ are cyclic. Moreover, every periodic bi--infinite path in the
Cayley graph $\Gamma (Q_2, S)$ is a Morse quasi--geodesic.
\end{enumerate}
\end{theorem}

\begin{proof}
$1$. Infinite torsion quotient group $Q_1$ of an arbitrary
non-elementary hyperbolic group $G$ was constructed in \cite{Ols}
(Corollary 2) as a direct limit of a sequence of hyperbolic groups
$G=G(0)\to G(1)\to\dots\to G(i-1)\to G(i)\to\dots$, where each
step is a transition from $H=G(i-1)$ to $H_1=G(i)= \langle G(i-1)
| V^m\rangle$ for a word $V=V_i$ and a sufficiently large $m=m_i$.
Lemmas 4.1 and 6.7 of \cite{Ols} claim that the set of the cyclic
shifts of the words $V^{\pm 1}$ satisfy a $C(\varepsilon, \mu,
\lambda, c,\rho)$-condition, where $\lambda=\lambda_i = \lambda(V,
i-1)>0$, $c=c_i=c(V, i-1)\ge 0$, the positive $\mu=\mu_i$ can be
selected arbitrary small, then $\varepsilon=\varepsilon_i$ can be
chosen arbitrary large, and afterwards $m$, and therefore
$\rho=\rho_i$, can be chosen arbitrary large.

It follows from lemma \ref{Chyp} that one can replace the defining
word $V_i^m$ by a conjugate in $H$ word $R_i$ having minimal
length in its conjugacy class, so that the set of cyclic shifts of
$R_i$ satisfy $C(\varepsilon'_i,\mu'_i,\rho'_i)$-condition with
parameters $\varepsilon'_n, \mu'_n, \rho'_n$ ($n=1,2,\dots$)
satisfying the definition \ref{classQ}. This proves the first
statement of the proposition.

$2.$ To construct $Q_2$ we denote by $\mathcal F$ the set of all
2-generated subgroups of $G$ and enumerate all elements $p_1, p_2,
\ldots $ of the set $\mathcal P = S \times \mathcal F$. We set
$G(0)=G$ and proceed by induction. Suppose that a (torsion--free)
hyperbolic group $G(i-1)$ and relators $R_1, \ldots , R_{i-1}$ are
already constructed. Then we consider the first pair, say
$p_k=(s,K)\in \mathcal P$, such that the image $K' $ of $K$ in
$G(i-1)$ is non--elementary and the image $s'$ of $s$ in $G(i-1)$
does not belong to $K^\prime $.  As in the proof of \cite[Corollary
1]{Ols}, we can choose a word $R_i$ of the form $R_i\equiv
X_0U^mX_1U^m\ldots X_lU^m$, where $X_0$ represents an element of
$sK$ and $U, X_1, \ldots , X_l$ represent elements of $K$, such that
the set of all cyclic shifts of $R_i^{\pm 1}$ satisfies a
$C(\varepsilon, \mu, \lambda, c, \rho)$-condition. Here
$\lambda=\lambda_i >0$, the positive $\mu=\mu_i$ can be selected
arbitrary small, and then $c=c_i\ge 0$ arbitrary large, then
$\varepsilon=\varepsilon_i$ can be chosen arbitrary large, and
afterwards $m$, and therefore $\rho=\rho_i$, can be chosen arbitrary
large. Such a choice of the parameters is guaranteed by lemmas 4.2
and 6.7 of \cite{Ols}. Then $G(i)=\langle G(i-1)|R_i\rangle$, and
the group $Q_2$ is defined to be the limit of the sequence
$G=G(0)\to G(1)\to \dots \to G(i)\to \dots$. Hence, as in the first
part of the proof, one can choose the parameters so that $Q_2$ has a
\gsc presentation. As in \cite[Corollary 1]{Ols}, $Q_2$ is a
non-Abelian torsion free group with cyclic proper subgroups. (The
only difference is that now we are adding only one relation $R_i$
when passing from $G(i-1)$ to $G(i)$, while in \cite{Ols}, the set
${\cal F}$ was enumerated, and, for given $K\in {\cal F}$, one
imposed finitely many relations to obtain $G(i)$, namely, one
relation for every $s\in S$.)

 To ensure the Morse property for bi--infinite periodic paths in
$\Gamma (Q_2, S)$ we have to make the following additional changes
in the scheme from \cite{Ols}. For every $i=1,2,\ldots $, when
passing from $G(i-1)$ to $G(i)$ we fix a set of words $\mathcal
V_i$ in $S^{\pm 1}$ of lengths at most $i$ having infinite order
in $G_{i-1}$. Arbitrary power of a word from ${\cal V}_i$ is
$(\lambda(i), c(i))$-quasi-geodesic in the hyperbolic group
$G(i-1)$ for some $\lambda(i)>0$ and $c(i)\ge 0$, and one can
chose the constants $\lambda_i$ and $c_i$ in the previous
paragraph so that $\lambda_i\le\lambda(i)$ and $c_i\ge c(i)$. Then
\cite[Lemma 2.5]{Ols} allows us to chose the words $R_i$ and the
parameters in the previous paragraph so that the following is
true.

$(*)$ {\em Let $\Delta $ be a minimal diagram over $G(i)$, $q$ a
part of $\partial \Delta $ such that $\Lab(q)=V_i^n$ for some $1\le
i\le s$, $n\in \mathbb N$. Then $\Delta $ contains no $\e_i
$--contiguity subdiagrams of $R_i$--cells to $q$ with contiguity
degree at least $\mu_i $.}

Assume that $V$ is a word representing a non--trivial element in
$Q_2$. Then $V\in \mathcal V_i$ for some $i$. Assume that for some
$n>0$, $V^n=U$ in $Q_2$, where $U$ is a geodesic word in $Q_2$.
Let $j=j(n)$ be the smallest positive integer such that $V^n=U$ in
$G(j)$. Suppose that $j\ge i$ (hence $V\in \mathcal V_j$).
Consider a minimal diagram $\Delta$ over $G(j)$ with $\partial
\Delta =pq$, where $\Lab (p)\equiv V^n$  and $\Lab (q^{-1})\equiv
U$. Let $\Pi $, $\Gamma _1$, $\Gamma _2$ be the $\mathcal
R_j$--cell of $\Delta $ and the $\e_j$--contiguity diagrams of
$\Pi $ to $p$ and $q$, respectively, provided by Lemma
\ref{old-gamma-cell}. Then $(\Pi,\Gamma_1, p)<\mu _j$ by  $(*)$,
and so $(\Pi,\Gamma _2,q)>1-24\mu _j$ that contradicts Lemma
\ref{primyk}. Hence $j<i$, i.e., $V^n=U$ in the hyperbolic group
$G(i)$, where $i$ is independent of $n$. This implies that any
bi--infinite $V$--periodic path in $\Gamma (Q_2, S)$ is
quasi--geodesic.

It remains to prove that any bi--infinite $V$--periodic path in
$\Gamma (Q_2, S)$ is Morse. Let us fix arbitrary $L, C>0$ and
consider any $(L,C)$--quasi--geodesic word $W$ such that $V^n=W$ in
$Q_2$ for some $n$. As above let $j$ be the smallest positive
integer such that $V^n=W$ in $G(j)$.

Since in a hyperbolic group every bi-infinite periodic geodesic is
Morse \cite{Short}, we would finish the proof if we show that $j$
can be bounded from above by some constant $J=J(L, C)$ independent
of $n$. Let us choose $J$ so that $\lambda_J<L$, $C<c_J$, and $V\in
\mathcal V_J$. Again we consider a minimal diagram $\Delta$ over
$Q_2$ with $\partial \Delta =pq$, where $\Lab (p)\equiv V^n$, and
$\Lab (q^{-1})\equiv W$.

By contradiction, assume that $\rank(\Delta)=k>J$.  Let $\Pi $,
$\Gamma _1$, $\Gamma _2$ be the $\mathcal R_k$--cell of $\Delta $
and the $\e_j$--contiguity diagrams of $\Pi $ to $p$ and $q$,
respectively, provided by Lemma \ref{old-gamma-cell}. Then
$(\Pi,\Gamma _2,q)>1-24\mu _j$ as above, contrary to Lemma
\ref{primyk}.
\end{proof}

\begin{remark}\label{fsub}
In a similar way, one can use methods of \cite{Ols} to construct an
infinite group $Q$ admitting a \gsc presentation, and such that all
proper subgroups of $Q$ are finite. In that construction, one would
have to use Theorem 4 from \cite{Ols}. Note that there exists a
slight error in the formulation of that theorem. Let $E_0$ be the
elementary group and $C$ its infinite cyclic normal subgroup from
the formulation of Theorem 4 \cite{Ols}. Since the group $H$
satisfies the quasi-identity $$x^2y=yx^2 \rightarrow xy=yx,$$ the
center $Z$ of $E_0$ has a finite odd index in $E_0$ by Proposition 2
\cite{Ols}. Hence $Z$ contains the Sylow 2-subgroup $P_2$ of $E_0$
and an infinite cyclic subgroup $C$ such that the product $CP_2$ is
of odd index in $E_0$. To make the formulation of Theorem 4
\cite{Ols} correct, one needs to add the condition that the cyclic
subgroup $C$ in $E_0$ is chosen with this additional property,
namely, the order of $E_0/(CP_2)$ is odd. (This condition was used
in the proof of Theorem 4 \cite{Ols}.)

Also the direct limits of hyperbolic groups in \cite{MO} can be
chosen satisfying Condition 3) of Theorem \ref{dirlim}. Hence there
exist torsion and torsion free examples of divisible (and even
verbally complete) lacunary hyperbolic groups.
\end{remark}

\subsection{Floyd boundary}
\label{fb}

Finally we note a relation between cut points in asymptotic cones
and the Floyd boundary. Recall that the {\it Floyd boundary}
$\partial G$ of a finitely generated group $G=\langle S\rangle $ is
defined as follows (see \cite{Flo}). Let $\dist _F$ be the metric on
$\Gamma =\Gamma (G,S)$ obtained by setting the lengths of each edge
$e$ to be equal to $(1+\dist(e,1))^{-2}$. It is easy to see that
$\Gamma $ is bounded with respect to $\dist _F$. Let
$\overline{\Gamma }$ be the metric completion of $(\Gamma , \dist
_F)$. Then $\partial G=\overline{\Gamma }\setminus\Gamma $.

The Floyd boundary of a group $G$ is a quasi--isometry invariant.
If $\partial G$ consists of $0$ (respectively $2$) points, $G$ is
finite (respectively virtually cyclic). If $\partial G$ consists
of $0$, $1$ or $2$ points it is said to be trivial. Otherwise it
is uncountable (and, moreover, $\partial G$ is a boundary in the
sense of Furstenberg). If $\partial G$ is nontrivial, $G$ contains
a free non--Abelian subgroup. In particular, $G $ is trivial for
any amenable group. The class of groups with nontrivial Floyd
boundary includes non--elementary hyperbolic groups,
non--elementary geometrically finite Kleinian groups, groups with
infinitely many ends, and many other examples. (For more details
we refer to \cite{Kar}.)

\begin{proposition}\label{Floyd}
Let $G$ be a finitely generated group whose Floyd boundary
consists of at least 2 points. Then all asymptotic cones of $G$
have cut points.
\end{proposition}

\begin{proof}
Let us fix a scaling sequence $d=(d_n)$ and an ultrafilter $\omega
$. Let also $(x_n)$, $(y_n)$ be sequences of elements of $G$ that
converge to distinct points $x,y\in \partial G$. For each $x_n$ we
fix a geodesic $\gamma _{n}$ in $\Gamma (G,S)$ connecting $x_n$ to
$1$. Since $\Gamma (G,S)$ is locally finite, there is an infinite
ray $\gamma $ such that $\gamma _-=1$ and the combinatorial length
of the common part of $\gamma $ and $\gamma _{n}$ tends to $\infty
$ as $n\to \infty$. This means, in particular, that the sequence
of vertices of $\gamma $ converge to $x$ as $n\to \infty$. Thus we
may assume that $|x_n|=n$ and $|y_n|=n$. Consider the subsequences
$a_n=x_{[d_n]}$ and $b_n=y_{[d_n]}$ and set $a=(a_n)^\omega$,
$b=(b_n)^\omega $. Note that $\dist (a, (1)^\omega )=\dist (b,
(1)^\omega )=1$.

Suppose that $\CG $ has no cut points. Then for some $\e>0$, there
is a path $p$ in $\CG \setminus {\rm Ball}\,((1)^\omega , \e  )$
of some length $L$ connecting $a$ to $b$. Now applying standard
methods it is easy to show that $a_n$ and $b_n$ can be connected
by a path $p_n$ in $\Gamma (G, S)$ such that $|p_n|<2Ld_n$ and
$p_n$ avoids the ball of radius $\e d_n/2$ centered at $1$ in
$\Gamma (G, S)$ (with respect to the combinatorial metric) {\oas}.
Hence we have
$$
\dist _F(a_n, b_n)\ge \frac{2Ld_n}{(\e d_n/2)^2} = o(1).
$$
Thus $(a_n)$ and $(b_n)$ converge to the same points of the Floyd
boundary and we get a contradiction.
\end{proof}

We note that the converse to Theorem \ref{Floyd} does not hold.
Indeed all asymptotic cones of groups constructed in this section
have cut points by Theorem \ref{classG7}. On the other hand their
Floyd boundary consists of a single point since they contain no
non--Abelian free subgroups and are not virtually cyclic.

\begin{cor}
There exists a lacunary hyperbolic group $G$ such that the Floyd
boundary $\partial G$ consists of a single point and all asymptotic
cones of $G$ are circle-trees.
\end{cor}


\section{Central extensions of lacunary hyperbolic groups}
\label{ceochg}


\subsection{Asymptotic cones of group extensions}
\label{acoge}

In this section, we obtain some results about asymptotic cones of
group extensions. These results are used in the next two sections.

\begin{lemma}[{\bf Asymptotic cones of isometry
groups}]\label{isometry} Suppose that a finitely generated group
$G$ acts isometrically on a metric space $X$. Fix an arbitrary
point $x\in X$. Then the map $\al\colon G\to X$ defined by $\al
(g)=gx$ for any $g\in G$ induces a continuous map $\hat
\al\colon\Con^\omega(G,d)\to \Con ^\omega (X, d)$.
\end{lemma}

\proof Note that the map $\al$ is C--Lipschitz for $$C= \max\{
\dist (x, sx)\, |\, s\in S^{\pm 1} \} ,$$ where $S$ is a finite
generating set of $G$. Indeed if $g=s_1\ldots s_n$ for some $g\in
G$ and $s_1, \ldots , s_n\in S^{\pm 1}$, then
$$
\dist (x, gx)\le \sum\limits_{i=1}^{n} \dist (g_{i-1}x,
g_{i-1}s_{i}x)\le \sum\limits_{i=1}^{n} \dist (x,s_{i}x)\le Cn.
$$
where $g_0=1$ and $g_i=s_1\cdots s_i$ for $1\le i\le n$. Thus for
any scaling sequence $d=(d_n)$ and any non--principal ultrafilter
$\omega $, the map $(g_n)^\omega\mapsto (g_nx)^\omega$ from $\Con
^\omega (G, d)$ to $\Con^\omega(X,d)$ is continuous.
\endproof

Given a group $G$ generated by a finite set $S$ and a normal
subgroup $N$ of $G$, we endow the group $G$ and the quotient group
$G/N$ by the word metric with respect to the set $S$ and its image
in $G/N$, respectively. We also assume that $N$ is endowed with the
metric induced from $G$. Thus for any $d$ and $\omega $, $\CN $ may
be considered as a subset of $\CG$. Set $x=1\in G/N$. Then Lemma
\ref{isometry} applied to the natural action of $G$ on the quotient
group $G/N$ by left multiplications gives us the map $\hat \alpha
\colon \CG \to \Con^\omega (G/N, d).$ Note that $\alpha $ is the
natural homomorphism $G\to G/N $ in this case. Given $b\in
\Con^\omega(G/N,d)$, we call the subset $\mathcal F_b=\hat\alpha
^{-1}(b)\subseteq \CG$ {\it fiber}.

Recall that the group
$$
\God =\left\{ \left. (g_n)\in \prod ^\omega G\;\right| \; |g_n|=\Oo
(d_n)\right\}
$$
acts transitively by isometries on $\CG$ by left multiplication. Let
$\Nod$ be the subgroup of $\God$ defined as follows:
$$
\Nod =\left\{ \left. (g_n)\in \prod ^\omega N\;\right| \; |g_n|=\Oo
(d_n)\right\}.
$$

\begin{theorem}[{\bf Asymptotic cones of quotient groups}]\label{factors}
Let $G=\la S\ra$ be a finitely generated group $G$, $N$ a normal
subgroup of $G$ endowed with the metric induced from $G$. Then the
following conditions hold.
\begin{enumerate}
\item[(a)] The map $\hat \al$ is surjective.

\item[(b)] For any $b\in \CGN$, we have $\mathcal F_b=\gamma\,\CN
$ for some element $\gamma \in \God$.

\item[(c)] The action of $\God $ permutes fibers, that is, for any $b\in \CGN $ and $\gamma\in \God$, $\gamma \mathcal F_b$ is a fibre.

\item[(d)] The action of $\Nod $ stabilizes each fiber (as a set) and acts on each fiber transitively.

\item[(e)] If $\Con ^\omega (N,d)$ is discrete, the map $\hat\alpha \colon \CG\to \CGN $ is locally isometric. If
$\Con ^\omega (N,d)$ consists of a single point, then $\CG $ and
$\Con ^\omega (G/N,d)$ are isometric.
\end{enumerate}
\end{theorem}

\proof Let $\sigma \colon G/N\to G$ be a section that assigns to
each element $x\in G/N$ a shortest preimage of $x$ in $G$. If
$h=(h_nN)^\omega$ is a point in $\Con^\omega(G/N,d)$, then
$g=(\sigma (h_nN))^\omega )$ belongs to $\CG$ and $\hat \alpha
(g)=h$. Thus $\hat \alpha $ is surjective.

Further for any $b\in \Con^\omega(G/N,d) $, we have
\begin{equation}\label{fib}
\mathcal F_b=\{ (g_n)^\omega \in \CG \mid (\alpha (g_n))^\omega
=b\} .
\end{equation}
Let us fix any element $(f_n)^\omega \in \mathcal F_b$ and set
$\gamma =(f_n)$. Clearly $\gamma \, \CN\subseteq \mathcal F_b$. If
$g=(g_n)^\omega \in \mathcal F_b$, then $(\alpha (f_n))^\omega
=(\alpha (g_n))^\omega $, i.e., $|\alpha (g_n^{-1})\alpha
(f_n)|=\oom (d_n)$. Let $s_n$ be a shortest preimage of $\alpha
(g_n^{-1}f_n)$ in $G$. Then $|s_n|=\oom (d_n)$ and
$u_n=f_n^{-1}g_ns_n\in N$. Note that $|u_n|\le |f_n|+|g_n|+|s_n|=\Oo
(d_n)$. Hence $u=(u_n)^\omega \in \CN $. Clearly $\gamma u=
(g_ns_n)^\omega =(g_n)^\omega$.  Thus $\mathcal F_b\subseteq \gamma
\, \CN $ and the second assertion is proved.

Similarly it is easy to show that each subset of the form $\gamma
\CN $, where $\gamma \in \God $, is a fiber. This and the second
assertion imply the third one.

Given two elements $(g_n)^\omega $ and $(h_n)^\omega $ of $\mathcal
F_b$, the element $(h_ng_n^{-1})^\omega $ belongs to $\Nod $ and
takes $(g_n)^\omega $ to $(h_n)^\omega $. This proves the fourth
assertion.

Finally, assume that a ball of radius $\e \in (0, \infty]$ in $\CN
$ consists of a single point. Let $h=(h_n)^\omega $ and
$g=(g_n)^\omega $ be two elements of $\CG $ such that $\dist (g,
h)< \e/2$. Note that
$$|\sigma (g_n^{-1}h_nN)|\le |g_n^{-1}h_n|< \e d_n/2 $$ \oas. Then
for $u_n=g_n^{-1}h_n \big(\sigma (g_n^{-1}h_nN)\big)^{-1}$, we
have $|u_n|<\e d_n $ \oas. Since $u_n\in N$, $(u_n)^\omega $
belongs to the ball of radius $\e $ in $\Con ^\omega (N, d)$
around $(1)^\omega $. Hence $(u_n)^\omega =(1)^\omega$, i.e.,
$|u_n|=o_\omega (d_n)$. Finally we obtain
$$
\begin{array}{rl}
\dist (\hat \alpha(g), \hat \alpha(h))=& \lio d(h_nN, g_nN)/d_n
=\lio d(g_n^{-1}h_nN, N)/d_n \\ & \\= & \lio |\sigma
(g_n^{-1}h_nN)|/d_n \le \lio (|g_n^{-1}h_n| + |u_n|)/d_n \\ & \\
= &
 \lio |g_n^{-1}h_n|/d_n=\dist(g, h).
\end{array}
$$
Thus the restriction of $\hat\alpha $ to any ball of radius $\e/4$
is an isometry.
\endproof

Recall that a subspace $Y$ of a metric space $X$  is said to be
{\it convex } if any geodesic path $p$ in $X$ such that $p_{\pm
}\in Y$ belongs to $Y$. In particular, if $X$ is geodesic, then
any convex subspace of $X$ is geodesic.

\begin{proposition}\label{section}
Let $d=(d_n)$ a scaling sequence, $\omega $ a non--principal
ultrafilter. Suppose that $T$ is a convex $\mathbb R$--tree in $\CG
$. Then for any point $y\in \CG $ such that $\hat\al (y)\in T$,
there is an isometric section $\sigma \colon T \to \CG $ of $\hat\al
$ such that $y\in \sigma (T)$.
\end{proposition}

To prove the proposition we need an auxiliary result.

\begin{lemma}\label{dist}
For any point $x\in \CG $ and any fiber $\mathcal F_b$, there is a
point $f\in \mathcal F_b$ such that $\dist (x, \mathcal F_b)=\dist
(x,f)=\dist (\hat \al (x), b)$.
\end{lemma}

\begin{proof}
By Theorem \ref{factors} (c), it suffices to prove the statement of
the lemma for $x=(1)^\omega $. Let $\mathcal F=\gamma \,\CN $, where
$\gamma = (g_n)\in \God $. We fix any section $s\colon G/N\to G$
that assigns to every element of $G/N$ a shortest preimage. Take
$f=(s (g_nN))^\omega $. By the choice of $s$, we have $|s
(g_nN)|=|\al (s (g_nN))|$. Hence
$$\dist (1,f)=\dist (1, (\al\circ s(g_nN))^\omega )= \dist (1,
\hat \al(f) )=\dist (1,b).$$ Note that $\gamma ^{-1}f= (g_n^{-1}s
(g_nN))^\omega \in \CN $ and thus $f\in \mathcal F$. Finally for any
point $f^\prime \in \mathcal F$ we have $f^\prime =(g_nu_n)^\omega$
for some $u_n\in N$. Hence
$$
\dist (1, f^\prime )=\lio\frac{|g_nu_n|}{d_n}\ge \lio\frac{|s
(g_nN)|}{d_n}=\dist (1, f).
$$
Thus $\dist (1, \mathcal F)=\dist (1,f)=\dist (1,b)$.
\end{proof}

\begin{proof}[Proof of Proposition \ref{section}]
We first consider a segment $I=[a,b]\subset T$ and any preimage
$x$ of $a$. Let $f$ be the point of $\mathcal F_b$ provided by
Lemma \ref{dist}. Recall that $\CG $ is a geodesic metric space.
Let $J_I=[x,f]$ be a geodesic segment in $\CG $. Note that
$\hat\al $ does not increase the distance. If $u,v\in J_I$ and
$\dist(\alpha (u),\alpha(v)) < \dist (u,v)$, then
$$
\begin{array}{rl}
\dist (a,b)\le & \dist (a, \alpha (u))+ \dist(\alpha
(u),\alpha(v)) +\dist(\alpha (v),b)\\ & \\
< &\dist (x,u)+ \dist (u,v) + \dist (v,f)=|J_I|
\end{array}
$$
that contradicts the choice of $f$. Thus $\hat \al$ isometrically
maps $J_I$ to its image in $\CGN $. Since $|J_I|=\dist (x,f)=\dist
(a,b)=|I|$ by the choice of $f$, $\hat\al (J_I)$ is a geodesic
segment in $\CGN $. As $T$ is convex, we have $\hat\al (J_I)\in T$.
Hence $\hat\al (J_I)=I$. Thus for any preimage $x$ of $a$, there is
an isometric section $\sigma _I\colon I\to \CG $ such that $\sigma
_I(a)=x$.

Let us fix a vertex $o$ of $T$. Suppose that we have already found
an isometric section $\sigma _{T_0}$ for a subtree $T_0\subset T$
containing $o$ such that $y\in T_0$. Let $b\in T\setminus T_0$. Then
there is a unique point $a\in T_0$ such that $[b,o]\cap T_0=[a,o]$.
Let $x=\sigma _{T_0} (a)$ and let $\sigma _I\colon I\to \CG $ be an
isometric section for $I=[a,b]$ such that $\sigma _I(a)=x$. We then
define an isometric section $\sigma_{T_1}\colon T_1\to \CG $, where
$T_1=T_0\cup [a,b]$, by the rule
$$
\sigma_{T_1} (t)= \left\{
\begin{array}{c}
\sigma _{T_0} (t), \; {\rm if\; } t\in T_0;
\\ \\
\sigma _I(t), \; {\rm if } \; t\in [a,b].
\end{array}
\right.
$$
Now we can complete the proof by transfinite induction.
\end{proof}

\begin{cor}\label{sectioncor}
Suppose that $\CGN$ is an $\mathbb R$--tree for some scaling
sequence $d=(d_n)$ and non--principal ultrafilter $\omega $. Then
there is an isometric section $\sigma \colon \CGN \to \CG $ of the
map $\hat\al$.
\end{cor}


\subsection{Central extensions of \ch groups}
\label{asoceochg}

We keep the notation from the previous section here.
Given a product $X\times Y$ of metric spaces $X$ and $Y$, by the
{\it product metric} we mean the metric on $X\times Y$ defined by
the rule
$$
\dist _{X\times Y} ((x_1,y_1), (x_2, y_2))=\dist _X(x_1,x_2)+\dist
_Y(y_1,y_2).
$$

Recall that every 2-dimensional cohomology class on a hyperbolic group can be represented by a bounded cocycle \cite{NR} (see also \cite{Min}, where it is proved for all dimensions $\ge 2$). This implies
that for any finitely generated group $G$ and any finitely generated
central subgroup $N\le G$ such that $G/N$ is hyperbolic, $G$ is
quasi--isometric to $N\times G/N$ \cite{NR}. (A particular case was also proved in \cite{Ger}.) Therefore, for any
$d$ and $\omega $, $\CG $ is bi--Lipschitz equivalent to $\CN\times
\CGN$ endowed with the product metric. In this section we generalize
this result to the class of \ch groups as follows.

\begin{theorem}\label{product}
Let $N$ be a central subgroup of a finitely generated group $G$
endowed with the induced metric. Suppose that for some
non--principal ultrafilter $\omega $ and some scaling sequence
$d=(d_n)$, $\CGN $ is an $\mathbb R$--tree. Then $\CG $ is
by--Lipschitz equivalent to $\CN\times \CGN $ endowed with the
product metric.
\end{theorem}

\begin{proof}
Let $\sigma \colon \CGN \to \CG $ be the isometric section provided
by Corollary \ref{sectioncor}. We define a map $\h \colon \CN \times
\CGN \to \CG $ as follows. Suppose that $x\in \CGN $ and
$g=(g_n)^\omega \in \CN $, where $(g_n)\in \Nod $. Then $\h(g,x)=
(g_n) \sigma (x)$.

First observe that $\h$ is well--defined. Indeed if $(g_n)^\omega
=(h_n)^\omega \in \CN $, then $g_n=h_nu_n$, where $u_n\in N$ and
$|u_n|=\oom (d_n)$. Since $N$ is central, for any $y=(y_n)^\omega
\in \CG$, we have
$$(g_n)(y_n)^\omega =(h_nu_ny_n)^\omega =(h_ny_nu_n)^\omega
=(h_ny_n)^\omega =(h_n)(y_n)^\omega .$$

Further observe that for any $y=(y_n)^\omega \in \CG $ and
$g_1=(g_{1n})^\omega, g_2=(g_{2n})^\omega \in \CN $ we have
$$
\dist (g_{1n},g_{2n})= \dist (y_ng_{1n}, y_ng_{2n})=\dist
(g_{1n}y_n, g_{2n}y_n).
$$
Therefore,
\begin{equation}\label{dg1g2}
\dist (g_1,g_2)=\dist ((g_{1n})y,(g_{2n})y).
\end{equation}

Suppose now that $(g_1,x_1)$, $(g_2, x_2)\in \CN \times \CGN $,
where $g_1=(g_{1n})^\omega, g_2=(g_{2n})^\omega \in \CN $, and
$k_1=\h(g_1,x_1)$, $k_2=\h(g_2,x_2)$. Applying (\ref{dg1g2}) we
obtain
\begin{equation}\label{dh1h2}
\begin{array}{rl}
\dist (k_1, k_2)= & \dist ((g_{1n})\sigma (x_1), (g_{2n})\sigma
(x_2))
\\&\\ \le & \dist ((g_{1n})\sigma (x_1), (g_{2n})\sigma (x_1))+\dist ((g_{2n})\sigma
(x_1), (g_{2n})\sigma (x_2))\\&\\= & \dist (g_1, g_2)+ \dist (\sigma
(x_1), \sigma (x_2)) \\&\\= & \dist (g_1, g_2)+\dist (x_1, x_2).
\end{array}
\end{equation}
Note that $ \hat\al(\h(g,x))=\hat\al(\sigma (x))=x$ since the action
of $\God $ preserves fibers. Hence
\begin{equation}\label{dx1x2}
\dist (x_1, x_2)=\dist (\hat\al(k_1), \hat\al (k_2))\le \dist
(k_1,k_2).
\end{equation}
Now reversing the inequality and replacing pluses with minuses in
(\ref{dh1h2}), we obtain
\begin{equation}\label{dg1g21}
\dist (g_1, g_2)\le \dist(k_1, k_2)+\dist (x_1,x_2)\le 2\dist(k_1,
k_2).
\end{equation}
Finally combining (\ref{dx1x2}) and (\ref{dg1g21}) we obtain
$$
\dist ((g_1,x_1), (g_2, x_2))=\dist (g_1,g_2)+\dist(x_1,x_2)\le
3\dist(k_1, k_2).
$$
This inequality together with (\ref{dh1h2}) shows that $\h$ is a
$3$--bi--Lipschitz map. To complete the proof it remains to note
that $\h$ is surjective by Theorem \ref{factors} (d).
\end{proof}

Our next goal is to prove Theorem \ref{conn}, which will be used in
the next section. We start with the following general lemma.

\begin{lemma}\label{path} Let $G$ be a non-virtually cyclic finitely generated
group, $\calc=\Con^\omega(G, d)$ its asymptotic cone. Let
$X\subseteq \calc$ be a finite subset. Then for every $x\in
\calc\setminus X$ there exists a path $p$ in $\calc\setminus X$ with
$p_-=x$, containing points arbitrary far away from $x$.
\end{lemma}

\proof Since $G$ is finitely generated, there exists an infinite
geodesic ray in the Cayley graph of $G$. Its ultralimit is an
infinite geodesic ray $r$ in $\calc$. Since $\calc$ is
homogeneous, we can assume that $r_-=x$. We can also assume that
$x=e=(1)^\omega$.

Let $X=\{x_1,...,x_m\}$, $x_i=(x_i(n))^\omega$. Let
$2l=\min\{\dist(e,x_i), i=1,...,m\}$, $l_i=\dist(e,x_i)$. For every
$g=(g_n)\in G^\omega_e(d)$ with $|g_n|\le ld_n$ consider the ray
$g\iv r$. The union of $g\iv r$ and any geodesic $[e, g\iv e]$ is a
path $p(g)$ starting at $e$ and containing points arbitrary far from
$e$.

If one of these paths does not contain any $x_i$, we are done.
Assume that every $p(g)$ contains $x_i$ from $X$. For every
$i=1,...,m$ let $M_i$ be the set of all $g=(g_n)^\omega$ with
$|g_n|\le ld_n$ such that $p(g)$ contains $x_i$. Then the union of
$M_i$ contains the limit $\lio\Ball_G(e,ld_n)$.

For every $g=(g_n)^\omega\in M_i$ $g\iv r$ must contain $x_i$
because $[e,g\iv e]$ is too short to contain $x_i$. Hence $gx_i\in
r$. Since $\dist(e,g e)\le l$, $\dist(gx_i, ge)=\dist(x_i, e)=l_i$,
we conclude that $gx_i$ belongs to the subgeodesic $r[l_i-l, l_i+l]$
of the ray $r$. Pick a number $N\ge 2$ and divide the interval
$r[l_i-l, l_i+l]$ into $N+1$ equal subintervals. Let $a_1,...,a_{N}$
be the division points, $a_j=(a_j(n))^\omega$. Then for every
$g=(g_n)^\omega\in M_i$ the point $g x_i$ is within distance
$\frac{2l}{N+1}$ from one of $a_j$. Hence $$\dist_G(g_nx_i(n),
a_j(n))\le \frac{2ld_n}{N+1}+o_\omega(d_n)\le \frac{2ld_n}{N}$$
$\omega$-almost surely. Therefore $g_nx_i(n)$ is in the ball $B_j$
of radius $\frac{2ld_n}{N}$ around $a_j(n)$ in $G$ $\omega$-almost
surely. Hence, for any fixed $i\in \{1,...,m\} $, the number of
elements $g_n$ such that $g_nx_i(n)\in B_j$ does not exceed $N$
times the number of elements of the ball $\Ball_G(\frac{2ld_n}{N})$.

Let $f(n)$ be the growth function of the group $G$. Since the union
of $M_i$ contains the limit $\lio\Ball_G(ld_n)$, we proved in the
previous paragraph that for every $N\ge 2$,

$$f(ld_n)\le mNf\left(\frac{ld_n}{N}\right)$$
$\omega$-almost surely.

By \cite[page 68]{Gro81} (see also \cite{DW}), then the asymptotic
cone $\Con(G,(ld_n))$ is locally compact. As shown in \cite{Point},
that asymptotic cone has Minkovski dimension 1 (for the definition,
see \cite{Point}). Then by \cite{Point} $G$ is virtually nilpotent
and $\Con^\omega(G,(ld_n))$ is homeomorphic to $\R$. But then
\cite[Proposition 6.1]{DS} implies that $G$ is virtually cyclic,
which contradicts the assumption of the lemma.
\endproof

For every metric space $X$, we define a {\it connectedness degree}
$c(X)$ as the minimal number of points of $X$ whose removal
disconnects $X$. If $X$ can not be disconnected by removing
finitely many points, we set $c(X)=\infty $. In particular, $X$
has cut points if and only if $c(X)=1$. By a {\it cut set} of $X$
we mean any subset of $X$ whose removal disconnects $X$.

\begin{theorem}\label{conn}
Let $N$ be a central subgroup of a finitely generated group $G$.
Suppose that $\CN $ consists of $m<\infty $ points for some
non--principal ultrafilter $\omega $ and some scaling sequence
$d=(d_n)$. Then $$c(\CG )=mc(\CGN).$$ Moreover, a finite subset
$C\subset \CG $ disconnects $\CG $ if and only if $C$ contains a
full preimage of a cut set of $\CGN$ under the map $\hat\al\colon
\CG \to \CGN $.
\end{theorem}

\begin{proof}
First note that $G$ is not virtually cyclic. Indeed otherwise either
$N$ is finite or $N$ contains an infinite cyclic group. In the first
case $\CG$ is $\R$, $m=1$, and the proposition is obvious. In the
second case $\CN$ contains infinitely many points that contradicts
$m<\infty$.

Take a finite set $C$ in $\CG$. Let $\hat C$ be the full preimage of
$\hat\al(C)$ under $\hat\al$. Then $\hat C$ is finite.

We shall need the following statement.

\begin{lemma}\label{sloi} Let $a, a'\in \CG\setminus \hat C$ and
$\hat\al(a)=\hat\al(a')$. Then $a$ and $a'$ are in the same
connected component of $\CG\setminus \hat C$.
\end{lemma}

\proof Let $\gamma\in\Nod$ be such that $\gamma a=a'$ (such $\gamma$
exists by Theorem \ref{factors}). By Lemma \ref{path}, there exists
a path $r$ in $\CG\setminus \hat C$ with $r_-=a$ containing points
arbitrary far away from $a$. Consider the path $r'=\gamma r$. Then
$r'_-=a'$. Note that $\hat C$ is closed under the action of $\Nod$.
Hence $r'$ does not contain points from $\hat C$. Consider a point
$z$ on $r$ such that $\dist(z,C)>\dist(a, a')$. Let $z'=\gamma z$
(using Theorem \ref{factors} again). Since $\CG$ is homogeneous,
there exists $\beta\in \God$ such that $z=\beta a$. Since $N$ is
central in $G$, we have
$$
\dist(z,z')=\dist (\beta a, \gamma \beta a)=\dist (\beta a, \beta
\gamma a)=\dist (a, a').
$$
Therefore any geodesic path from $z$ to $z'$ avoids $\hat C$.

Now consider the path $w$ that goes first from $a$ to $z$ along $r$,
then from $z$ to $z'$ along any geodesic $[z,z']$, then goes back to
$a'$ along $r'$. That path avoids points from $\hat C$. \endproof

Let us continue the {\em proof of Theorem \ref{conn}}.

Suppose that the finite set $C$ does not contain a full preimage
of a cut set of $\CGN$ under $\hat\al$. We need to show that $C$
is not a cut set of $\CG$ that is any two points $u, v\in
\CG\setminus C$ can be connected by a path avoiding $C$. Take any
geodesic $p$ connecting $u, v\in \CG\setminus C$. If this geodesic
does not contain points in $C$, we are done. Suppose that $p$
contains a point from $C$. Then it is enough to show how to
replace subpaths of $p$ connecting points close enough to points
in $C$. Thus without loss of generality we can assume that $u$ and
$v$ are from $\Ball(c,\zeta)\setminus \hat C$ for some $c\in C$
and some small enough $\zeta$.

Pick any $\zeta$ such that $\Ball(c,\zeta)\cap \hat C=\{c\}$. Pick
two points $u, v\in \Ball(c,\zeta)\setminus\{c\}$. Then $u,
v\not\in\hat C$.

Suppose that $\hat\al(C)$ is not a cut set of $\CGN$. Then
$\hat\al(u)$ can be connected with $\hat\al(v)$ by a path $q$ in
$\CGN\setminus\hat\al(C)$. Since $\hat\al$ is a covering map by
Theorem \ref{factors} (e), we can lift the path $p$ to a path
$\hat p$ in $\CG$ avoiding $\hat C$ and such that $\hat p_-=u$.
Note that $\hat\al(\hat p_+)=\hat\al(v)$. Hence by Lemma
\ref{sloi}, we can connect $\hat p_+$ with $v$ by a path $p_1$ in
$\CG\setminus\hat C$. The composition of $\hat p$ and $p_1$
connects $u$ and $v$ and avoids $C$ as required.

Now suppose that $\hat\al(C)$ is a cut set in $\CGN$. Let $C'$ be
the union of all fibers of $\hat\al$ contained in $C$. We can
assume that $\hat\al(C')$ is not a cut set. Therefore there exists
a path from $\hat\al(u)$ to $\hat\al(v)$ in
$\CGN\setminus\hat\al(C')$. That path must contain points
$\hat\al(c_1)$ for some $c_1\in C$ with $c_1\not\in C'$. Lifting
this path to $\CG$ and using Lemma \ref{sloi}, we can obtain a
path connecting $u$ and $v$ in $\CG\setminus C'$ but containing
points in $C\setminus C'$. If we could replace parts of this path
connecting points close the points in $C\setminus C'$ by paths
avoiding these points we would show that $C$ is not a cut set.
Thus without loss of generality we can assume that $u$, $v$ are
very close to a point in $C\setminus C'$, i.e. we can assume that
$c\in C\setminus C'$.

Then there exists a point $c'\not\in C$ with
$\hat\al(c')=\hat\al(c)$. Let $\gamma\in\Nod$ be such that $\gamma
c=c'$. We can assume that $\zeta$ is small enough so that
$\gamma\Ball(c,\zeta)$ does not intersect $C$. By Lemma
\ref{sloi}, there exist paths $r, r'$ in $\CG\setminus \hat C$
with $r_-=u, r_+=\gamma u$, $r'_-=v, r'_+=\gamma v$. Then one can
travel from $u$ to $v$ by first going to $\gamma u$ along $r$ then
by a path from $\gamma u$ to $\gamma v$ inside $\gamma\Ball(c,
\zeta)$, then back to $v$ along $r'$. That path avoids $C$. Hence
$C$ is not a cut set in $\CG$.

We have proved that every cut set of $\CGN$ contains a full
preimage of a cut set of $\CGN$. This implies $$c(\CG)\ge
mc(\CGN).$$

In order to prove the opposite inequality, note that if a finite set
$\{b_1,\dots,b_c\}$ disconnects $\CGN$ then the union of fibers
${\cal F}_{b_1},\dots,{\cal F}_{b_c}$ disconnect $\CN$ because
$\hat\al$ is continuous.
\end{proof}


\subsection{Applications}
\label{a}

Let us consider a family of central extensions of a \ch group
constructed as follows. By Corollary (\ref{subset}) and Remark
\ref{exist} there is a presentation $$H=\left\la a, b\;\left|
\;\bigcup_{i=1}^\infty \calr _i\right.\right\ra $$ that
simultaneously satisfies the \gsc condition and the classical small
cancellation condition $C'(1/24)$, where for every $i$, $\calr_i$
consists of cyclic shifts of a single word $R_i$ and its inverses.
Throughout this section we fix any such presentation and denote by
$r_n$ the length of the word $R_n$. Given a sequence of integers
$k=(k_n)$, where $k_n\ge 2$, we consider the central extension of
$H$ defined as follows.
\begin{equation}
G (k)=\left\langle a,\, b\; \left|\; [R_n, a]=1,\, [R_n, b]=1,\,
R_n^{k_n}=1,\, n=1,2,\dots\right.\right\rangle
\end{equation}

We begin with auxiliary results.

\begin{lemma}\label{lengthRn}
Let $U$ be a subword of a word $R_n^{k_n}$ of length at most
$k_nr_n/2$. Then the length of the element represented by the word
$U$ in $G$ is at least $|U|/8$.
\end{lemma}

\begin{proof}
Suppose that $U=V$ in $G$ where the word $V$ is geodesic in $G$. Let
$\Delta$ be a diagram over $G$ corresponding to this equality. Let
$\partial\Delta=pq\iv$ be the decomposition of the boundary of
$\Delta$, where $U\equiv\Lab(p)$, $V\equiv \Lab(q)$. We can turn
$\Delta$ into a diagram $\Delta^\prime$ over $H$ with the same
boundary label by the following procedure. Every cell corresponding
to the relation $R_i^{k_i}$ is replaced by a union of $k_i$ cells
each labeled by $R_i$ connected by a point. Every cell labeled by
$[R_i,a]$ or $[R_i,b]$ is replaced by the union of two cells labeled
by $R_i$ connected by an edge labeled by $a$ or $b$. In both cases
the union has the same boundary label as the original cell (see Fig.
\ref{fig5}).

\begin{figure}
\vspace{1mm}
  \hspace{25mm}
  \includegraphics{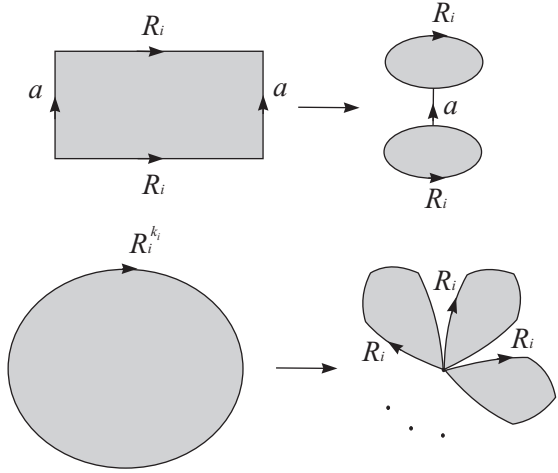}\\
  \vspace{-5mm}
  \caption{Transforming cells in $\Delta $}\label{fig5}
\end{figure}

To any cell $\Pi $ labeled by $R_n^{\pm 1}$ in $\Delta ^\prime $, we
assign the number $-1$ if $R_n$ reads along $\partial \Pi $ in the
clockwise direction and $1$ otherwise. By the {\it algebraic number}
of $R_n$--cells in $\Delta ^\prime $ we mean the sum of the assigned
numbers over all cells labeled by $R_n^{\pm 1}$. Notice that the
algebraic number of $R_n$--cells in $\Delta ^\prime $ is divisible
by $k_n$. If we reduce $\Delta^\prime $ by canceling pairs of cells
having a common edge and being mirror images of each other, we
obtain a diagram $\Delta^{\prime \prime}$ with the same property. We
keep the notation $pq^{-1}$ for the boundary of $\Delta ^{\prime
\prime}$. We have to show that
\begin{equation}\label{qp4}
|q|\ge |p|/8.
\end{equation}

Assume that the boundary of some cell $\Pi$ in $\Delta ^{\prime
\prime}$ has a common subpath of length at least $|\partial \Pi
|/24$ with $p$. Then $\partial \Pi $ is labeled by $R_n$ by the
$C^\prime (1/24)$--condition, and since $p$ is labeled by a power
of $R_n$ we may assume without loss of generality that $\partial
\Pi \subset p$. After cutting such a cell off we obtain a new
diagram with boundary $p^\prime q^{-1}$, where $\Lab (q^\prime )$
is again a power of $R_n$. Continuing this process we obtain a
subdiagram $\Sigma $ of $\Delta ^{\prime \prime}$ with boundary
$tq^{-1}$, where $t$ is labeled by a power of $R_n$, such that no
$R_n$--cell $\Pi $ of $\Sigma $ has a common subpath of lengths at
least $|\partial \Pi |/24$ with $t$. Now there are two cases to
consider.

{\it Case 1.} Suppose that $|t|\ge |p|/4$. Observe that any cell
$\Pi $ in $\Sigma $  satisfies the following condition:

\smallskip

\noindent{\it $(+)$ For any two vertices $x,y\in t\cap
\partial \Pi $, there is a common subpaths of $\partial \Pi $
and $t$ connecting $x$ and $y$.}

\smallskip

Indeed otherwise there is a subdiagram $\Xi $ of $\Sigma $ such that
$\Xi $ contains at least one cell and $\partial \Xi =s_1s_2$, where
$s_1$ is a subpath of $t$ and $s_2$ is a subpath of $\partial \Pi $.
Note that a common subpath of the boundary of any cell $\Omega $ in
$\Xi $ and $s_1$ (respectively $s_2$) has length less than
$|\partial \Omega |/24$ by the construction of $\Sigma $
(respectively since $\Delta ^{\prime\prime}$ is reduced). This
contradicts Lemma \ref{greendlin}.

In particular, $(+)$ and the $C^\prime (1/24)$--condition imply that
for any cell $\Pi $ of $\Sigma $, the intersection $\partial \Pi
\cap t$ is a path of length less than $|\partial \Pi |/24$. Note
that if the number of common edges of $t$ and $q$ is at least
$|p|/8$, the inequality (\ref{qp4}) is obvious. Hence we may assume
that more than $|t|-|p|/8\ge |p|/8$ edges of $t$ belong to cells of
$\Sigma $. Therefore the sum of perimeters of all cells in $\Sigma $
is greater than $24 |p|/8 =3|p|$. Applying Lemma \ref{sum-of-per} to
$\Delta ^{\prime\prime}$ we obtain $|p|+|q|\ge 9|p|/4$, which yields
(\ref{qp4}).

{\it Case 2.} Suppose that $|t|\le  |p|/4$. This means that we have
to cut at least $k=\frac{3|p|}{4r_n}$ $R_n$--cells to get $\Sigma $
from from $\Delta ^{\prime\prime}$. Note that $k\le |p|/r_n\le
k_n/2$ and all these cells have the same orientation by the
$C^\prime (1/24)$--condition. Since the algebraic number of
$R_n$--cells in $\Delta ^{\prime\prime}$ should be divisible by
$k_n$, the total number of $R_n$--cells in $\Delta ^{\prime\prime}$
is at least $2k$. Applying Lemma \ref{sum-of-per} again, we obtain
$|p|+|q|\ge 6kr_n/4\ge 9 |p|/8$. Hence $|q|\ge |p|/8$.
\end{proof}

\begin{theorem}\label{cext1} For every $m \ge 2$
there exists a finitely generated group $G$ such that for any
ultrafilter $\omega $ and any scaling sequence $d=(d_n)$, exactly
one of the following possibilities occurs and both of them can be
realized for suitable $\omega $ and $d$.
\begin{enumerate}
\item[(a)] $\CG $ is an $m$--fold cover  of a
circle--tree and $c(\CG )=m$. Moreover, a finite subset $C\subset \CG
$ disconnects $\CG $ if and only if $C$ contains a fiber of the map
$\hat\al\colon \CG \to \CGN $.

\item[(b)] $\CG $ is an $\mathbb R$--tree.
\end{enumerate}
In particular, in both cases $\CG $ is locally isometric to an
$\mathbb R$--tree.
\end{theorem}

\begin{proof}
Let $G=G(k)$ be the group corresponding to the sequence $k_n=m\ge 2$
for all $n$. The central subgroup $N=\langle R_1, R_2, \ldots
\rangle $ inherits a metric from $G$. Let $g=(g_n)^\omega\in \CN$.
If $g_n\ne 1$, then
$g_n=R_1^{\e_1}R_2^{\e_2}...R_{j(n)}^{\e_{j(n)}}$, where $0\le
\e_i\le m-1$ for all $i$ and $\e_{j(n)}\ne 0$. Note that
$$
\dist \left( g_n, R_{j(n)}^{\e_{j(n)}}\right)\le
\sum\limits_{i=1}^{j(n)-1} \left| R_i^{\e_i}\right| \le (m-1)
\sum\limits_{i=1}^{j(n)-1} \left| R_i\right| = \oom(R_{j(n)})
$$
by (${\bf Q}_3$), Lemma \ref{asc} (c), and Lemma \ref{lengthRn}. Hence
\begin{equation}\label{gnRin}
(g_n)^\omega =(R_{i_n}^s)^{\omega}
\end{equation}
for some $s\in \{ 0, 1, \ldots , m-1\} $. Now there are two cases to
consider.

{\it Case 1}. Suppose that $(r_n)$ is not $(\omega , d)$--visible.
Then (\ref{gnRin}) implies $(g_n)^\omega=(1)^\omega $, i.e., $\CN $
is a point. Applying the last assertion of Theorem \ref{factors} and
Theorem \ref{classG7}, we obtain (b).

{\it Case 2}. Assume now that $(r_n)$ is $(\omega , d)$--visible
(see Definition \ref{od-vis}). Let $(R_{i_n})$ be the sequence such
that $r_{i_n}=\To (d_n)$. Note any sequence $(i_n^\prime )$
satisfying $r_{i_n^\prime }=\To(d_n)$ is $\omega$-equal to $(i_n)$
(that is, $i_n=i_n^\prime $ $\omega$-almost surely) by $(Q_4)$ and
Lemma \ref{asc} (c). Therefore $\Con^\omega(N,d)$ contains at most
$m$ points.

On the other hand, the points $(R_{i_n}^s)^\omega $ are different
for different values of $s\in \{ 0, 1, \ldots , m-1\} $. Indeed if
$(R_{i_n}^s)^\omega =(R_{i_n}^t)^\omega $ for some  $s\ne t$,
$s,t\in \{ 0, 1, \ldots , m-1\} $, then $|R_{i_n}^{s-t}|=\oom
(d_n)$. Passing from $s-t$ to $m-(s-t)$ if necessary, we obtain
$|R_{i_n}^{l}|=\oom (d_n)=\oom (r_{i_n})$ for some $0<l\le m/2$.
However this contradicts Lemma \ref{lengthRn}. Hence $\CN $
consists of exactly $m$ points. Applying Theorem \ref{factors}
again, we obtain that the map $\CG \to \CGN $ induced by the
natural homomorphism is locally isometric and each fiber consists
of $m$ points, i.e., $\CG $ is an $m$--fold cover of $\CGN $. Note
that by Theorem \ref{classG7} $\CGN$ is a circle--tree. Hence it
contains cut points by Lemma \ref{cut}. Therefore a finite subset
$C\subset \CG $ disconnects $\CG $ if and only if $C$ contains a
fiber of the map $\hat\al\colon \CG \to \CGN $ by Theorem
\ref{conn}. In particular, $c(\CG )=m$.

Finally we note that ($Q_4$) and Lemma \ref{asc} (c) guarantees
that the second case occurs. The first case occurs when taking
$d_n=r_n$ and any $\omega $.
\end{proof}

Below we denote by $\mathbb S^1$ the unit circle with the lengths
metric.

\begin{theorem}\label{cext2}
There exists a finitely generated group $G$ and a scaling sequence
$d=(d_n)$ such that for any ultrafilter $\omega $, $\CG $ is
bi--Lipschitz equivalent to the product of an $\mathbb R$--tree and
$\mathbb S^1$. In particular, $\pi _1 (\CG )=\Z$.
\end{theorem}

\begin{proof}
Let $G(k)$ be the group corresponding to a sequence $k=(k_n)$ such
that
\begin{equation}\label{kn}
k_n\to \infty\;\;\; {\rm and}\;\;\; k_nr_n =o (r_{n+1}).
\end{equation}
The existence of such a sequence is guaranteed by the equality
$r_n=o(r_{n+1})$, which follows from ($Q_4$) and Lemma \ref{asc}
(c). Set $d_n=k_nr_n$, $d=(d_n)$. As in the proof of the previous
theorem let $g=(g_n)^\omega\in \CN$. If $g_n\ne 1$, then
$g_n=R_1^{\e_1}R_2^{\e_2}...R_{j(n)}^{\e_{j(n)}}$, where $0\le
\e_i\le k_i-1$ for all $i$ and $\e_{j(n)}\ne 0$. Note that
$$
\dist \left( g, R_{j(n)}^{\e_{j(n)}}\right)\le
\sum\limits_{i=1}^{j(n)-1} \left| R_i^{\e_i}\right| \le
\sum\limits_{i=1}^{j(n)-1} k_i \left| R_i\right| = \oom(R_{j(n)})
$$
by (\ref{kn}) and Lemma \ref{lengthRn}. Therefore $(g_n)^\omega
=\left(R_{j(n)}^{s(n)}\right)^\omega$ where $$s(n)\in \{ 0, 1,\ldots
, k_n-1\}.$$ Note that if $(g_n)^\omega \ne (1)^\omega $, then
$r_{j_n}=\To(|g_n|)=\To (d_n)$ and hence $j_n=n$ {\oas }. Thus
$(g_n)^\omega \in \lio p_n$, where $p_n$ is the cycle in $\Gamma (G,
\{ a, b\} )$ that begins and ends at $1$ and has label $R_n^{k_n}$.
Since $|R_n|=r_n =o(d_n)$, $\lio (p_n)$ coincides with the set of
points of type $\left(R_n^{s(n)}\right)^\omega\in \CN $ and we
obtain $\CN = \lio p_n$.

Observe that by Lemma \ref{lengthRn} $p_n$ is $(1/4,
1)$--quasi--isometric to a circle of lengths $|p_n|=k_nr_n=d_n$.
Hence $\lio p_n$ is bi-Lipschitz equivalent to the unit circle. Note
that by (\ref{kn}) the sequence $r_n$ is not $(\omega ,
d)$--visible. Hence $\CGN $ is an $\mathbb R$--tree by Theorem
\ref{classG7}. Applying now Theorem \ref{product}, we obtain that
$\CG $ is bi--Lipschitz equivalent to the product of an $\mathbb
R$--tree and $\mathbb S^1$.
\end{proof}

\begin{remark}
Arguing as in the proof of Theorem \ref{cext1} it is not hard to
classify all asymptotic cones of the group $G(k)$ from the proof of
Theorem \ref{cext2} as follows. For any scaling sequence $d_n$ and
any ultrafilter $\omega $, exactly one of the conditions (a)--(c)
below holds and all possibilities can be realized.
\begin{enumerate}
\item[(a)] $(r_n)$ is $(\omega , d)$--visible, $\CG $ is an infinite degree cover of a
circle--tree.

\item[(b)] $(k_nr_n)$ is $(\omega , d)$--visible and $\CG $ is
bi--Lipschitz equivalent to the product of an $\mathbb R$--tree and
a unit circle.

\item[(c)] There exists an $(\omega , d)$--visible sequence $(c_n)$ such that $c_n=o(k_nr_n)$ and $r_n=o(c_n)$. In this case $\CG $ is
bi--Lipschitz equivalent to the product of an $\mathbb R$--tree and
$\mathbb R$.

\item[(d)] There is no $(\omega , d)$--visible sequence $(c_n)$ such that $r_n\le c_n\le k_nr_n$. In this case $\CG $ is an $\mathbb R$--tree.
\end{enumerate}
In particular, even the finiteness of the connectedness number
$c(\CG )$ for a given group $G$ is not invariant under changing $d$
and $\omega $.
\end{remark}


\section{Lacunar hyperbolicity and divergence}
\label{chad}



\subsection{Divergence of non-constricted groups}


The following general statement allows one to estimate the
divergence function of a group with no cut points in some of its
asymptotic cones. Recall that given a path $p$ in a metric space, we
denote by $p_-$ and $p_+$ the beginning and the ending points of $p$
respectively. The length of $p$ is denoted by $|p|$.

\begin{theorem}\label{divergence} Let $G$ be a finitely generated group.
Suppose that for some sequence of scaling constants $d_n$ and
every ultrafilter $\omega$, the asymptotic cone
$\Con^\omega(G,(d_n))$ does not have cut points. Let $f(n)\ge n$
be a non-decreasing function such that $d_n\le f(d_{n-1})$ for
all sufficiently large $n$. Then the divergence
function $\dv(n)$ of $G$ does not exceed $Cf(n)$ for some
constant $C$ (and all $n$).
\end{theorem}

\proof Since the asymptotic cone does not change if we change a
finite subsequence of $(d_n)$, we can assume without loss of
generality that $d_1=\frac14$. Taking a constant multiple of $f$ if
necessary, we can assume that $d_n\le f(d_{n-1})$ for all $n\ge 2$.

Since $\calc$ does not have cut points for any choice of $\omega$,
by \cite[Theorem 2.1]{DMS}, we can conclude that there exists
a constant $C_1$ such that $\dv(n)$ is bounded
by $C_1n$ for every $n$ in any interval $[\frac{d_k}{18}, 18d_k]$
for every $k\ge 2$.

Let $\delta=\frac14, \lambda=2$. For every $n\ge 1$ choose elements $a_n, b_n, c_n$, in $G$ with
$\dist(a_n,b_n)\le n$ and such that  $\dv_\lambda(a_n,b_n,c_n;\delta)$ is
maximal possible, i.e. $$\dv(n)=\dv_\lambda(a_n,b_n,c_n;\delta).$$

Suppose, by contradiction, that $\dv(n)$ is not smaller than
$Cf(n)$ for some constant $C$ and all $n$. Then for every
$m\ge 1$ there exists $n=n(m)$ such that $\dv(n)>mf(n)$. Pick
$m> 12+18C_1$. Let $a=a_n, b=b_n, c=c_n$ where $n=n(m)$. Let
$r=\dist(c, \{a,b\})$. Let $B$ be the ball of radius $\delta r$
around $c$. Without loss of generality assume that $\dist(c,a)=r$.

Note that any geodesic $h$ connecting $a$ and $b$ passes through $B$
since $\dv_\gamma(a,b,c;\delta)>n$. Hence $r\le 2n$. (Indeed, if $r>2n$
every point in $h$ is at distance at least $r-|h|\ge
r-n>\frac{r}{2}>\delta r$ from $c$, and cannot belong to $B$.) Let
$c'$ be a point in $h\cap B$, so $c'\in h$ is at distance at most
$\delta r$ from $c$. Then $\dist(a,c')\le r(1+\delta)\le 2r$.

Let $b'$ be either $b$ or the point between $c'$ and $b$ at distance
$2r$ from $c'$. Let $h'$ be the part of $h$ between $c'$ and $b$.
Then $h'$ does not intersect $B$ (any point in $h'$ is at distance
$\ge (1-\delta)r > 2\delta r$ from $B$).

Since $\delta<\frac13$, $\dist(a,B), \dist(b,B)>2\delta r$ which
exceeds the diameter of $B$.

Since the Cayley graph $\Gamma$ of $G$ is infinite, homogeneous and
locally finite, for every vertex $x$ in $\Gamma$ there exists a
bi-infinite geodesic $q(x)$ passing through $x$.

Consider the geodesic $q(a)$. The point $a$ cuts $q(a)$ into two
geodesic rays $l(a)$ and $l(a)'$. Since $\dist(a,B)$ is greater than
the diameter of $B$, one of these rays does not pass through the
ball $B$. Let it be $l(a)$. Similarly, let $l(b')$ be a geodesic ray
starting at $b'$ and not passing through $B$.

Choose the smallest $k\ge 2$ such that $r/2\le d_k$. Note that then
$d_{k-1}\le r/2$ (even if $k=2$ since $d_1=\frac14$ and $r\ge 1$).
Let $x$ be the point in $l(a)$ at distance $5d_k$ from $a$, and let
$y$ be the point in $l(b')$ at distance $5d_k$ from $b'$.

{\bf Case 1.} Suppose that $\dist(x,y)<\frac{d_k}{2}$. Then consider
a geodesic $p$ connecting $x$ and $y$. Any point in $p$ is at
distance at least $5d_k-\frac{d_k}{2}>2r$ from $a$. Hence any point
in $p$ is at distance at least $2r-r>\delta r$ from $c$. Thus $p$
does not intersect $B$, so we found a path $[a,x]\cup p\cup
[y,b]\cup h'$ of length at most $11d_k+n$ connecting $a$ and $b$.
Then
$$11d_k+n\ge \dv_\lambda(a,b,c;\delta)=\dv(n)\ge mf(n).$$
But $11d_k+n\le 11f(d_{k-1})+n\le 11f(r/2)+n\le 12f(n)$ since $f$ is
a non-decreasing function and $f(n)\ge n$. Thus $12f(n)\ge mf(n)$, a
contradiction since $m> 12$.

\medskip

{\bf Case 2.} Suppose that $\dist(x,y)\ge \frac{d_k}{2}$. Since
$\dist(x,y)\le 10d_k+\dist(a,b')\le 10d_k+4r\le 18d_k$, the distance
$\dist(x,y)$ is in the interval $[\frac{d_k}{18}, 18d_k]$, and so
there exists a path $p$ of length at most $C_1\dist(x,y)$ avoiding
the ball of radius $\delta\dist(c, \{x,y\}) > \delta(5d_k-r)\ge
\delta(5d_k-2d_k)> \delta r$ around $c$. Then the path $[a,x]\cup
p\cup [y,b']\cup h'$ connects $a$ and $b$, avoids $B$, and has
length at most
$$10d_k+C_1(18 d_k)+n\le (10+18
C_1)f(d_{k-1})+f(n)\le (11+18C_1)f(n),$$ a contradiction since
$m>11+18C_1$.
\endproof


\subsection{Torsion groups with slow non-linear divergence}
\label{tgwsbnld}

Let $F_2=\la a,b\ra$ be the free group of rank $2$. We fix an
arbitrary odd prime $p$ and a large odd power $n_0$ of $p$, say,
$n_0>10^{80}$. Let $G(0)=F_2$, i.e. the set of relators ${\cal
R}_0$ of rank $0$ is empty. The set of {\em periods} of rank $0$
is empty by definition. Below we define the sets ${\cal R}_i$ of
defining relations of groups $G(i)$ and an increasing sequence
$d=(d_r)$ by induction. This sequence depends on a non-decreasing
function $\phi$ such that $\phi(0)=0$, $\phi(1)=1$, $\phi (r)\ge
2$ for every $r=2,3\dots$, and $\lim_{r\to\infty}\phi(r)=\infty$.

Two arbitrary segments from the set $\{ (d_r/\phi(r), \phi(r)
d_r]\mid r=1,2,\dots\}$ will have empty intersection, and $d_0=1$.
After $d_{r-1}$ is defined we introduce $d_r$ and then define all
groups $G(i)$ for $\phi(r-1)d_{r-1}<i\le \phi(r)d_r$.

Assume $r>0$ and $d_0,\dots d_{r-1}$ are already defined along
with hyperbolic groups $$G(0),G(1),\dots,G(i_{r-1}),$$ where
$i_{r-1}= [\phi(r-1)d_{r-1}]$. Let $\delta _{i_{r-1}}$ be the
hyperbolicity constant of $G_{i_{r-1}}$. We chose a minimal
integer $d_r$ such that
\begin{equation}\label{dr}
d_r\ge\max\{\phi(r)^2 d_{r-1}, \,  \phi(r)^2\delta _{i_{r-1}},\,
2\} .
\end{equation}
For example, $d_1=2$. We also define $i_r=[\phi(r)d_r]$.

Then we argue by induction on $i$ ($i_{r-1}<i\le i_r$). A word $A$
is called {\em simple in rank} $i-1$ if it is not conjugate in
rank $i-1$ (that is in $G(i-1)$) to $B^m$, where $|B|<|A|$ or $B$
is a period of rank $j\le i-1$.

Let ${\cal X}_i$ be a maximal set of simple in rank $i-1$ words of
length $i$ such that for two different $A,B\in {\cal X}_i$, we
have that $A$ is not conjugate of $B^{\pm 1}$ in rank $i-1$. All
the words from ${\cal X}_i$ are called {\em periods} of rank $i$.

For every period $A$ of rank $i$ ($i_{r-1}<i\le i_r$), we
introduce a large odd exponent $n_A$, where $n_A $ is a minimal
power of $p$ such that

$n_A\ge \max(n_0, d_r/i)$ for  $i_{r-1}<i< d_r/\phi(r)$, and

$n_A=n_0$ for $d_r/\phi(r)\le i\le i_r$.

The set ${\cal R}_i$  is, by definition, equal to ${\cal
R}_{i-1}\cup \{A^{n_A}| A\in {\cal X}_i\}$, and $G(i)=\la a,b |
{\cal R}_i\ra $. We will show in Lemma \ref{hyperb} that the group
$G(i)$ is hyperbolic. Finally, $G=G(n_0,\phi)=\la a,b
|\cup_{i=0}^{\infty}{\cal R}_i\ra$ .

Since $d_0=1$ and $d_1= 2$ we can chose ${\cal
R}_1=\{a^{n_0},b^{n_0}\}$, and $G(1)$ is the free product of two
cyclic groups of order $n_0$. Since $i_1\ge d_1=2$, one can set
${\cal R}_2=\{a,b,ab, ab^{-1}\}$. Thus $a$ and $b$ are periods of
rank $1$, and $ab$ and $ab^{-1}$ are periods of rank $2$. Hence,
for every word $w$ of length at most $2$, we have $w^{n_0}=1$ in
$G$.

The proof of the following Proposition is based on \cite{book} and
is contained in the next section.

\begin{proposition}\label{main1} The 2-generated group $G$ satisfies the
following properties
\begin{itemize}
\item[(a)]  The natural homomorphism $G(i_{r-1})\to G(i_{r})$ is
injective on the ball of radius $Kd_r/\phi(r)$, for a non-zero
constant $K$.

\item[(b)]{\bf [bounded torsion up to a small
deformation]}\label{wuton} There is a constant $c>1$ such that for
every large enough integer $r$ and every word $W$ with
$cd_r/\phi(r)<|W|_G<d_r \phi(r)/c$,  there exists a word $U$ of
length $\le 1$ such that  $(WU)^{n_0}=1$ in $G$. \item[(c)]

{\bf [relations are locally quasi-geodesic]} For every large
enough $r$ there exists a relation $u_r^{n_0}=1$ in
$\cup_{i=0}^{\infty}{\cal R}_i$ with $|u_r|=d_r$, such that no
non-empty subword $w$ of $u_r^{n_0}$ of length $<n_0|u_r|/2$ can
be equal to a word of length at most $|w|/2$ in $G$.

\item[(d)] {\bf [finite cyclic centralisers]} $G$ is an infinite
$p$-group in which the centralizer of every non-trivial element is
cyclic.

\item [(e)] The hyperbolic constant $\delta_{i_r}$ of the group
$G(i_r)$ is  $O(\phi(r)d_r)$.

\item [(f)] The order of arbitrary word $X$ in the generators of
$G$ is  $O(\phi(|X|)^3)$.
\end{itemize}
\end{proposition}

\begin{theorem}\label{th2} Let $G$ and $d=(d_n)$ be the group and the scaling sequence constructed above. Then $G$ is {\ch},  but for any
non--principal ultrafilter $\omega $, the asymptotic cone ${\frak
C}=\Con^\omega(G,d)$ does not have cut points.
\end{theorem}

First we note that by (\ref{dr}) we have $\delta _{r-1}=o
(d_r/\phi(r))$. Together with the first assertion of Proposition
\ref{main1} and Theorem \ref{dirlim} this implies the first
assertion of the Theorem.

Let us now prove that ${\frak C}$ does not have cut points. Indeed
suppose this is not so. Then ${\frak C}$ is a tree-graded space
with respect to the collection $\calc_1$ of maximal subsets
without cut points (see Lemma \ref{cutting}). Let $u_r$, $r=1,2,
\dots$ be the words given by Proposition \ref{main1} (c). Then the
limit of the loops $\pgot_r$ in the Cayley graph of $\cg$
corresponding to the relations $u_r^{n_0}=1$ is a non-trivial loop
in ${\frak C}$ (of length $n_0$). Indeed, if $\lio{\pgot_r}$ is
not a simple loop then $\omega$-a.s. there are two points $x_r$
and $y_r$ on $\pgot_r$ that are distance $\Theta(d_r)$ apart along
the loop $\pgot_r$ but $o(d_r)$-close in the Cayley graph of $G$.
But this would contradict Proposition \ref{main1} (c). Therefore
${\frak C}$ is not an $\R$-tree.

Hence some of the pieces in ${\frak C}$ contain infinitely many
points.

Since ${\frak C}$ is homogeneous, one of these pieces, $M$,
contains $O=(1)^\omega$. Let $A$ be another point in $M$.

As in the proof of \cite[Lemma 6.10]{DS} consider two cases.

{\bf Case 1.} Suppose that there are two pieces from $\cc_1$ that
intersect. Then every point is in two distinct pieces. Then we can
construct a geodesic $\g\colon [0,s]\to {\mathfrak C}_1 $ such
that $s=\Sigma_{i=1}^\infty s_i$ with $0<s_i<\frac{1}{i^2}$ and
$\g \left[ \Sigma_{i=0}^r s_i,\Sigma_{i=0}^{r+1}s_i\right]\subset
M_r$ for some pieces $M_r$, where $M_r\neq M_{r+1}$ for all $r\in
\N \cup \{ 0 \}$. Here $s_0=0$. Such a geodesic exists by Lemma
\ref{sir}. We call such a geodesic \textit{fractal at the arrival
point}. That geodesic with reverse orientation will be called
\textit{fractal at the departure point}. If $\g$ is fractal at the
departure point, $\g'$ is fractal at the arrival point,
$\g_+=\g'_-$, we can construct (using Lemma \ref{sir}) a geodesic
$\pgot$ which is a composition of an initial piece of $\g$ and the
terminal piece of $\g'$. The geodesic $\pgot$ is then fractal at
the departure and arrival points or \textit{bifractal}. By
homogeneity, every point in ${\mathfrak C}$ is the endpoint of a
bifractal geodesic.

Let $[A,B]$ be a bifractal geodesic.  Lemma \ref{strconv}, part
(2), implies that $[A,B]$ can intersect a piece $M$ containing $A$
in $A$ or in a non-trivial sub-geodesic $[A,B']$. Since $[A,B]$ is
fractal at the departure point the latter case cannot occur. It
follows that the intersection of $[A,B]$ and $M$ is $\{A\}$. There
exists an isometry $\gamma=(x_r)\in G^\omega_e(d)$ such that
$\gamma(O)=B$. Since $[A,B]$ is fractal at the arrival point also,
it follows that $[A,B]\cap \gamma M=\{B\}$. For every $Z\in \gamma
M$ we have that $[A,B]\cup [B,Z]$ is a geodesic, by Lemma
\ref{sir}. In particular $A$ is the projection of $\gamma M$ onto
$M$. A symmetric argument gives that $B=\gamma(O)$ is the
projection of $M$ onto $\gamma M$.

Note that in the argument of the previous paragraph we only used
the fact that $\gamma O=B$. Let us change $\gamma$ a little bit
preserving the property. Then the conclusions of the previous
paragraph will still be true.

Let $w_r$ be the shortest word  representing $x_r\iv$ in $G$. Note
that $\frac1C d_r\le |w_r|\le Cd_r$ for $\omega$-almost every $r$
where $C$ is a constant. By Proposition \ref{main1} (b),
$\omega$-a.s. there exist words $u_r$ with $|u_r|\le 1$ such that
$(w_ru_r)^{n_0}=1$ is in ${\cal R}$. Let $\beta=(y_rx_r)\in
G^\omega_e(d)$ where $y_r\in G$ is represented by $u_r\iv$. Then
again $\beta(O)=B$. Notice that $\beta^{n_0}=1$.

For every $k\ge 2$ consider the following piecewise geodesic path
from $\beta^kO\in \beta^k M$ to $A$:
\begin{equation}\label{eq7}
\g=[\beta^k O,\beta^{k-1}A]\cup [\beta^{k-1}A,\beta^{k-1}O]\cup
[\beta^{k-1}O, \beta^{k-2}A]\cup...\cup [\beta O,A]\end{equation}
where every odd numbered segment is a bifractal geodesic and every
even segment is a non-trivial geodesic inside a piece. By Lemma
\ref{sir}, $\g$ is a geodesic which is not inside $M$. By the
strong convexity of pieces in a tree-graded space (Lemma
\ref{strconv}), we conclude that $\beta^k O\not\in M$, so $\beta^k
M\ne M$. This contradicts the equality $\beta^{n_0}=1$.

{\bf Case 2.} Now suppose that all pieces in $\cc_1$ are disjoint.
Note that we could repeat the argument from Step 1 if we found an
isometry $\delta$ from $G^\omega_e(d)$ such that
$B=\delta(O)\not\in M$, the projection of $M$ onto $\delta(M)$ is
$B$ and the projection of $\delta(M)$ onto $M$ is $A$. Indeed, by
slightly changing $\delta$, we can find an isometry $\beta$ with
the same property and, in addition, $\beta^{n_0}=1$. On the other
hand, for every $k$, consider the piecewise geodesic curve
(\ref{eq7}). Every even numbered geodesic segment $\g_i$ in it is
non-trivial and inside a piece $M_i$, and every odd numbered
geodesic segment connects a point $U_i$ in a piece $M_i$ and a
point $U_{i+1}$ in a piece $M_{i+1}$ such that $U_i$ is the
projection of $M_{i+1}$ onto $M_i$ and $U_{i+1}$ is a projection
of $M_i$ onto $M_{i+1}$. Hence by Lemma \ref{projA}, the odd
numbered geodesic segments $\g_i$ intersect $M_i$ (resp.
$M_{i+1}$) in exactly one point. By Lemma \ref{sir} the curve
(\ref{eq7}) is a geodesic, and so $\beta^k(O)\ne O$ for any $k$, a
contradiction.

Thus our goal is to find such $\delta$.

Lemma \ref{cutting}, part (b), implies that $A$ is the projection
of a point $B\in {\mathfrak C}\setminus M$. Let $\gamma=(x_r)$ be
an isometry from $G^{\omega}_e$ such that $\gamma(O)=B=[x_r\iv]$.
If $[A,B]$ intersects $\gamma(M)$ in $B$ then we have found the
desired $\delta=\gamma$ by Lemma \ref{projA} (since the isometries
of $\frak C$ permute the pieces of $\calc_1$, $\gamma(M)\in
{\calc_1}$).

Assume $[A,B]\cap \gamma(M)=[B',B]$, $B'\ne B$ Since all the
pieces are disjoint, $B'\neq A$. We have $B'=\gamma(A')$ for some
$A'\in M$. Since the space ${\frak C}$ is homogeneous, and all
pieces of ${\cal C}_1$ are disjoint, the stabilizer of $\gamma(M)$
in $G^\omega_e(d)$ acts transitively on $\gamma(M)$. Then there
exists $\gamma'$ in it such that $\gamma'(B)=B'$. We have that
$\gamma'\gamma(M)=\gamma(M)$ projects onto $M$ in $A$ and $M$
projects onto $\gamma'\gamma(M)$ in $B'=\gamma'\gamma(O)$, so we
can take $\delta=\gamma'\gamma$.\endproof

Proposition \ref{main1} and Theorem \ref{th2} imply that the
divergence and the orders of elements of $G(n_0, \phi)$ can grow
arbitrarily slow.

\begin{cor}\label{cormain}
For any positive function $f$ with $f(r)/r$ non-decreasing and
$\lim_{r\to\infty}f(r)/r=\infty$, and for any non-decreasing
function $g(r)$ with $\lim_{r\to\infty}g(r)=\infty$, there is a
function $\phi$ such that

(a) for some $\lambda$, the divergence function $\dv(r,\lambda)$
of the group $G=G(n_0,\phi)$ is $O(f(r))$, is not linear, but does
not exceed a linear function on an infinite subset of $\N$;

(b) The order  of any element $x\in G(n_0,\phi)$ is $O(g(|x|))$.
\end{cor}

\proof Consider the functions $f$ and $g$ from the formulation of
the corollary. Then we can choose a non-decreasing function $\phi$
such that $\lim_{r\to\infty}\phi(r)=\infty$, $\phi(r)\ge 2$ for
any $r\ge 2$, and $\phi(r)^3<  \min (cf(r-1)/(r-1), g(r))$ for
every $r\ge 2$ and a constant $c$. Then, by (\ref{dr}) and
condition (e) of Proposition \ref{main1}, we have $d_r \le
c'\phi(r)^3d_{r-1}$ for some constant $c'$ and every $r>0$. The
right-hand side is less than $c'cf(r-1)d_{r-1}/(r-1)\le C
f(d_{r-1})$ for $C\ge c'c$ since the function $f(r)/r$ is
non-decreasing and $d_{r-1}\ge r-1 $ by (\ref{dr}), if $r\ge 2$.
Besides, $d_1<Cf(d_0)$ if $C$ is large enough. Now the  statement
on the divergency follows from Theorem \ref{divergence}. Condition
(f) obviously implies the second statement of the corollary.
\endproof

\subsection{The proof of Proposition \ref{main1}}
\label{tpop}

  Here we present the proof of Proposition \ref{main1}. We can apply
lemmas from \cite{book} to the construction of group $G$ from that
proposition because it obviously satisfies the $R$-condition from
$\S 25$ \cite{book} since we do not use relations of the 'second
type' here. In particular, {\em every reduced diagram of rank $i$
arising below is a $B$-map} by Lemma 26.5 \cite{book}. (See the
definition in subsection 20.4 \cite{book}.) The contiguity
  diagrams we use now are more particular than those in the previous
sections.
  (See their definition in subsection 20.1 of\cite{book}).

\begin{notation}
As in \cite[Chapter 7] {book} we fix certain positive numbers
$\eta <<\zeta<<\epsilon<<\gamma << \beta <<\alpha$
 between 0 and 1 where "$<<$" means ``much smaller". Here ``much"
means enough to satisfy all the inequalities in Chapters 7 and 8
of \cite{book}. We also have $n_0^{-1}<<\eta.$
\end{notation}

Denote by $P(i)$ the maximum of $n_A|A|$ for the periods $A$ of
rank at most $i$, that is the maximum length of relations of rank
at most $i$.

\begin{lemma}\label{hyperb}
The group $G(i)$ is $\delta_i$-hyperbolic for arbitrary $i\le
i_r=[\phi(r)d_r]$, where $\delta_i=n_0P(i)$.
\end{lemma}

\proof

{\bf Step 1.} First we want to prove that a geodesic subpath  $p$
of a boundary $\partial\Delta$, where $\Delta$ is a reduced
diagram of rank $i$, is a smooth section of rank $k$ (see the
formulations of the smoothness conditions $S1-S5$ in subsection
20.4 of \cite{book}) if we define $\rank(p)=
k=[(1-2\beta)^{-1}P(i)+1]$.

Condition $S1$ holds since $p$ is geodesic.

Assume that $\Gamma$ is a contiguity diagram of a cell $\Pi$ to a
a geodesic subpath $p$ in a reduced diagram $\Delta$ over $G(i)$
with $(\Pi,\Gamma,p)\ge\varepsilon$. Let
$\partial(\Pi,\Gamma,p)=p_1q_1p_2q_2$ . Then by Lemma  21.2
\cite{book}, $|q_1|>(1-2\beta)|q_2|$. Indeed, the proof of Lemma
21.2 \cite{book} does not change if one replaces 'smooth $p$' by
'geodesic $p$'.

Let $A^{n_A}$ be the label of $\partial\Pi$. Since $|q_1|\le
P(i)$, we have $|q_2|<(1-2\beta)^{-1}P(i)$. This implies condition
$S2$ from $\S 20.4$ \cite{book} for $p$ since then
$|q_2|<(1+\gamma)k$. Since our relations are of the first type
(i.e., of the form $A^n_A$, and so the boundary of any cell is
just one 'long section'; see section 25 \cite{book}), $S2$ implies
$S3$, and $S4$ is obvious because $1<\alpha^{-1}$. Condition $S5$
automatically holds since $k>i$, and therefore there are no cells
of rank $k$ in $\Delta$. We conclude that  $p$ is a {\em smooth
section of rank $k$} in $\partial\Delta$.

{\bf Step 2.} Now let  $xyz$ be a triangle in the Cayley graph of
$G(i)$ with geodesic paths $x$, $y$, and $z$. To prove that $x$
belongs to $n_0P(i)$-neighborhood of $y\cup z$, we introduce an
inscribed geodesic hexagon $\Psi=t_1x't_2y't_3z'$, where (1) $x',
y'$, and $z'$ are subpaths of $x,y$, and $z$, respectively, (2)
$\max (|t_1|, |t_2|, |t_3|)\le \eta^{-1}k$, and the sum
$|x'|+|y'|+|z'|$ is minimal for hexagons satisfying (1) and (2).
(It follows that if $|x'|>0$, then the distance between any point
of $x'$ and $y'\cup z'$ is at least $\eta^{-1}k$. Similar
properties hold for $y'$ and $z'$.)

Then we have decompositions $x=x_1x'x_2$, $y=y_1y'y_2$, and
$z_1z'z_2$. If $|x_1|< 2\gamma^{-1}\eta^{-1}k$, then every point
$o$ of $x_1$ is at the distance at most $2\gamma^{-1}\eta^{-1}k$
from $z_2$. If $|x_1|\ge 2\gamma^{-1}\eta^{-1}k\ge
2\gamma^{-1}|t_1|$, then we can apply Lemma 22.4 \cite{book} to a
reduced diagram with boundary $x_1t_1^{-1}z_2$ where, according to
Step 1, sections $x_1$ and $z_2$ are smooth of rank $k$. By Lemma
22.4 \cite{book}, again, every point $o$ of $x_1$ can be connected
to a point of $z_2$ by a path of length
$<2\gamma^{-1}\eta^{-1}k<3\gamma^{-1}\eta^{-1}P(i)<\frac12
n_0P(i)$.

Similarly, the distance from every point of  $x_2$ to $y_1$ is
less than $\frac12 n_0P(i)$. Thus to complete the proof, it
suffices to show that $|x'|<n_0P(i)$. Proving by contradiction, we
suppose $|x'|\ge n_0P(i)>\frac12 n_0 k$ and consider a reduced
diagram $\Delta$ of rank $i$ with boundary $z't_1x't_2y't_3$.
First let us check that $\Delta$ is a $C$-map in the meaning of
section 23.1 \cite{book} with 3 long sections of the first type
having rank $k$, namely, $s_0=z'$, $s_1=x', s_2=y'$, with short
sections $t_1,t_2,t_3$, and with paths $p_1, p_2, q$ of zero
length (in the notation of \cite{book}).

Since $|s_1|=|x'|\ge \frac12 n_0 k$, the diagram $\Delta$
satisfies condition $C1$ and $C2$ with $j=k$ because $l=3-1=2$ for
3 sections $s_1,s_2,s_3$. Conditions $C3$ and $C4$ hold since
$|q|=|p_1|=|p_2|=0$ and $s_1, s_2, s_3$ are geodesic in $\Delta$.
Condition $C5$ holds since $\max\{|t_1|, |t_2|,
|t_3|\}<\eta^{-1}k$ (and $\eta^{-1}=d$ in the book \cite{book}).
It follows from the choice of $\Psi$ that there are no contiguity
subdiagrams between either $x'$ and $y'$ or $x'$ and $z'$, or $y'$
and $z'$, since, by Lemma 21.1(1) \cite{book}, the side arcs of
such a subdiagram would be of  length $<\eta^{-1}k$. Therefore
$\Delta$ satisfy condition $C6$. The condition $C7$ holds for the
same reason as $S4$ at Step 1. Thus $\Delta$ is a $C$-map.

Since $\Delta$ is a $C$-map, there must be a contiguity subdiagram
between a pair of sides from $\{x', y', z'\}$ by Lemma 23.15
\cite{book}. But this is impossible as was shown in the previous
paragraph; a contradiction.
\endproof

\begin{lemma}\label{d1'}
The natural homomorphism $G(i_{r-1})\to G(i_{r})$ is injective on
the ball of radius $0.4 d_r/\phi(r)$.
\end{lemma}

\proof Let $w$ be a word equal to $1$ in $G(i_r)$ but not in
$G(i_{r-1})$. Then there is a reduced diagram $\Delta$ such that
its boundary label is $w$ and it contains a cell $\Pi$ of a rank
$j>i_{r-1}$. It follows from the construction of defining words of
rank $j>i_{r-1}$ that the perimeter $|\partial\Pi|$ of $\Pi$ is at
least $d_r/\phi(r)$. Therefore, by Lemma 23.16 of \cite{book},
$$|w|=|\partial\Delta|>(1-\alpha)|\partial\Pi|\ge (1-\alpha)d_r/\phi(r)>0.8
d_r/\phi(r).$$ The inequality $|w|>0.8 d_r/\phi(r)$ gives the
injectivity radius at least $0.4 d_r/\phi(r)$.
\endproof

The following lemma seems to be known. The short proof of it has
been communicated to the authors by Ian Agol.

\begin{lemma}\label{puti}
Let $o_1$ and $o_2$ be two distinct points on the boundary of a
double punctured disk $D$, and $x$ the boundary cycle starting
(and ending) at $o_1$. Let $y$ be a simple path connecting $o_1$
and $o_2$ in $D$ and separating the punctures, and $z$ a simple
loop starting at $o_1$ and going around exactly one of the
punctures. Then every simple path connecting $o_1$ and $o_2$ in
$D$ is homotopic to $x^s p_1 x^{-s}p_2$ where each of $p_1$, $p_2$
is a product of at most 2 factors from $\{x, y,z^{\pm 1}\}$.
\end{lemma}

\proof Let $\cal D$ be the diffeomorphism group of $D$. (The
elements of $\cal D$ leave the boundary $\partial D$ invariant and
may permute the punctures.) Denote by ${\cal D}_0$ the subgroup of
$\cal D$ that fixes every point on the boundary of the disc
$\partial D$. It is well known that, modulo diffeomorphisms
isotopic to the identity element, ${\cal D}_0$ is the cyclic braid
group $B_2=\langle\sigma\rangle$ .

\begin{figure}[!ht]
\centering
\unitlength .6mm 
\linethickness{0.4pt}
\ifx\plotpoint\undefined\newsavebox{\plotpoint}\fi 
\begin{picture}(0.88,100)(60,-10)
\put(99.63,44.75){\line(0,1){1.365}}
\put(99.6,46.11){\line(0,1){1.363}}
\multiput(99.52,47.48)(-.042,.453){3}{\line(0,1){.453}}
\multiput(99.4,48.84)(-.0441,.3383){4}{\line(0,1){.3383}}
\multiput(99.22,50.19)(-.04526,.26919){5}{\line(0,1){.26919}}
\multiput(99,51.54)(-.04599,.22277){6}{\line(0,1){.22277}}
\multiput(98.72,52.87)(-.04646,.18936){7}{\line(0,1){.18936}}
\multiput(98.39,54.2)(-.04675,.16407){8}{\line(0,1){.16407}}
\multiput(98.02,55.51)(-.04693,.1442){9}{\line(0,1){.1442}}
\multiput(97.6,56.81)(-.04701,.12813){10}{\line(0,1){.12813}}
\multiput(97.13,58.09)(-.04702,.11482){11}{\line(0,1){.11482}}
\multiput(96.61,59.35)(-.04696,.10359){12}{\line(0,1){.10359}}
\multiput(96.05,60.6)(-.04686,.093948){13}{\line(0,1){.093948}}
\multiput(95.44,61.82)(-.046712,.085568){14}{\line(0,1){.085568}}
\multiput(94.78,63.01)(-.046523,.078195){15}{\line(0,1){.078195}}
\multiput(94.09,64.19)(-.046298,.071644){16}{\line(0,1){.071644}}
\multiput(93.35,65.33)(-.04604,.065771){17}{\line(0,1){.065771}}
\multiput(92.56,66.45)(-.045752,.060465){18}{\line(0,1){.060465}}
\multiput(91.74,67.54)(-.047958,.058731){18}{\line(0,1){.058731}}
\multiput(90.88,68.6)(-.047462,.05392){19}{\line(0,1){.05392}}
\multiput(89.97,69.62)(-.046953,.04952){20}{\line(0,1){.04952}}
\multiput(89.03,70.61)(-.048754,.047749){20}{\line(-1,0){.048754}}
\multiput(88.06,71.57)(-.050487,.045912){20}{\line(-1,0){.050487}}
\multiput(87.05,72.49)(-.054896,.046329){19}{\line(-1,0){.054896}}
\multiput(86.01,73.37)(-.059716,.046725){18}{\line(-1,0){.059716}}
\multiput(84.93,74.21)(-.065016,.047099){17}{\line(-1,0){.065016}}
\multiput(83.83,75.01)(-.070885,.047453){16}{\line(-1,0){.070885}}
\multiput(82.69,75.77)(-.077432,.047783){15}{\line(-1,0){.077432}}
\multiput(81.53,76.48)(-.0848,.048091){14}{\line(-1,0){.0848}}
\multiput(80.34,77.16)(-.086522,.04492){14}{\line(-1,0){.086522}}
\multiput(79.13,77.79)(-.094904,.044894){13}{\line(-1,0){.094904}}
\multiput(77.9,78.37)(-.10454,.0448){12}{\line(-1,0){.10454}}
\multiput(76.64,78.91)(-.11578,.04462){11}{\line(-1,0){.11578}}
\multiput(75.37,79.4)(-.12908,.04433){10}{\line(-1,0){.12908}}
\multiput(74.08,79.84)(-.14515,.04392){9}{\line(-1,0){.14515}}
\multiput(72.77,80.24)(-.16501,.04333){8}{\line(-1,0){.16501}}
\multiput(71.45,80.58)(-.19029,.04251){7}{\line(-1,0){.19029}}
\multiput(70.12,80.88)(-.22368,.04134){6}{\line(-1,0){.22368}}
\multiput(68.78,81.13)(-.27007,.03964){5}{\line(-1,0){.27007}}
\multiput(67.43,81.33)(-.3392,.037){4}{\line(-1,0){.3392}}
\put(66.07,81.48){\line(-1,0){1.361}}
\put(64.71,81.57){\line(-1,0){1.364}}
\put(63.35,81.62){\line(-1,0){1.365}}
\put(61.98,81.62){\line(-1,0){1.364}}
\put(60.62,81.56){\line(-1,0){1.361}}
\multiput(59.26,81.46)(-.339,-.0386){4}{\line(-1,0){.339}}
\multiput(57.9,81.3)(-.26988,-.04089){5}{\line(-1,0){.26988}}
\multiput(56.55,81.1)(-.22349,-.04238){6}{\line(-1,0){.22349}}
\multiput(55.21,80.85)(-.19009,-.04339){7}{\line(-1,0){.19009}}
\multiput(53.88,80.54)(-.16481,-.04409){8}{\line(-1,0){.16481}}
\multiput(52.56,80.19)(-.14494,-.04459){9}{\line(-1,0){.14494}}
\multiput(51.26,79.79)(-.12888,-.04493){10}{\line(-1,0){.12888}}
\multiput(49.97,79.34)(-.11557,-.04515){11}{\line(-1,0){.11557}}
\multiput(48.7,78.84)(-.10433,-.04528){12}{\line(-1,0){.10433}}
\multiput(47.45,78.3)(-.094695,-.045333){13}{\line(-1,0){.094695}}
\multiput(46.21,77.71)(-.086313,-.04532){14}{\line(-1,0){.086313}}
\multiput(45.01,77.08)(-.078938,-.045251){15}{\line(-1,0){.078938}}
\multiput(43.82,76.4)(-.07721,-.048141){15}{\line(-1,0){.07721}}
\multiput(42.66,75.67)(-.070664,-.04778){16}{\line(-1,0){.070664}}
\multiput(41.53,74.91)(-.064798,-.0474){17}{\line(-1,0){.064798}}
\multiput(40.43,74.1)(-.059499,-.047){18}{\line(-1,0){.059499}}
\multiput(39.36,73.26)(-.054682,-.046582){19}{\line(-1,0){.054682}}
\multiput(38.32,72.37)(-.050274,-.046145){20}{\line(-1,0){.050274}}
\multiput(37.32,71.45)(-.048532,-.047974){20}{\line(-1,0){.048532}}
\multiput(36.35,70.49)(-.046724,-.049737){20}{\line(0,-1){.049737}}
\multiput(35.41,69.5)(-.047212,-.054139){19}{\line(0,-1){.054139}}
\multiput(34.51,68.47)(-.047686,-.058952){18}{\line(0,-1){.058952}}
\multiput(33.66,67.41)(-.048146,-.064245){17}{\line(0,-1){.064245}}
\multiput(32.84,66.31)(-.045736,-.065983){17}{\line(0,-1){.065983}}
\multiput(32.06,65.19)(-.045967,-.071857){16}{\line(0,-1){.071857}}
\multiput(31.32,64.04)(-.046161,-.07841){15}{\line(0,-1){.07841}}
\multiput(30.63,62.87)(-.046315,-.085783){14}{\line(0,-1){.085783}}
\multiput(29.98,61.67)(-.046425,-.094164){13}{\line(0,-1){.094164}}
\multiput(29.38,60.44)(-.04648,-.1038){12}{\line(0,-1){.1038}}
\multiput(28.82,59.2)(-.04649,-.11504){11}{\line(0,-1){.11504}}
\multiput(28.31,57.93)(-.04642,-.12835){10}{\line(0,-1){.12835}}
\multiput(27.85,56.65)(-.04626,-.14442){9}{\line(0,-1){.14442}}
\multiput(27.43,55.35)(-.046,-.16429){8}{\line(0,-1){.16429}}
\multiput(27.06,54.03)(-.04558,-.18957){7}{\line(0,-1){.18957}}
\multiput(26.74,52.71)(-.04496,-.22298){6}{\line(0,-1){.22298}}
\multiput(26.47,51.37)(-.04401,-.26939){5}{\line(0,-1){.26939}}
\multiput(26.25,50.02)(-.0425,-.3385){4}{\line(0,-1){.3385}}
\put(26.08,48.67){\line(0,-1){1.36}}
\put(25.96,47.31){\line(0,-1){1.363}}
\put(25.89,45.94){\line(0,-1){1.365}}
\put(25.88,44.58){\line(0,-1){1.364}}
\put(25.91,43.21){\line(0,-1){1.362}}
\multiput(25.99,41.85)(.0441,-.4528){3}{\line(0,-1){.4528}}
\multiput(26.12,40.49)(.0456,-.3381){4}{\line(0,-1){.3381}}
\multiput(26.3,39.14)(.0465,-.26898){5}{\line(0,-1){.26898}}
\multiput(26.54,37.8)(.04702,-.22256){6}{\line(0,-1){.22256}}
\multiput(26.82,36.46)(.04733,-.18914){7}{\line(0,-1){.18914}}
\multiput(27.15,35.14)(.04751,-.16385){8}{\line(0,-1){.16385}}
\multiput(27.53,33.83)(.0476,-.14398){9}{\line(0,-1){.14398}}
\multiput(27.96,32.53)(.0476,-.12791){10}{\line(0,-1){.12791}}
\multiput(28.43,31.25)(.04755,-.1146){11}{\line(0,-1){.1146}}
\multiput(28.96,29.99)(.04744,-.10337){12}{\line(0,-1){.10337}}
\multiput(29.53,28.75)(.047295,-.09373){13}{\line(0,-1){.09373}}
\multiput(30.14,27.53)(.047107,-.085351){14}{\line(0,-1){.085351}}
\multiput(30.8,26.34)(.046885,-.077979){15}{\line(0,-1){.077979}}
\multiput(31.5,25.17)(.046629,-.071429){16}{\line(0,-1){.071429}}
\multiput(32.25,24.02)(.046344,-.065557){17}{\line(0,-1){.065557}}
\multiput(33.04,22.91)(.046031,-.060253){18}{\line(0,-1){.060253}}
\multiput(33.87,21.83)(.045691,-.055429){19}{\line(0,-1){.055429}}
\multiput(34.73,20.77)(.047711,-.0537){19}{\line(0,-1){.0537}}
\multiput(35.64,19.75)(.047182,-.049303){20}{\line(0,-1){.049303}}
\multiput(36.59,18.77)(.048974,-.047523){20}{\line(1,0){.048974}}
\multiput(37.56,17.82)(.053367,-.048082){19}{\line(1,0){.053367}}
\multiput(38.58,16.9)(.05511,-.046074){19}{\line(1,0){.05511}}
\multiput(39.63,16.03)(.059932,-.046448){18}{\line(1,0){.059932}}
\multiput(40.7,15.19)(.065234,-.046798){17}{\line(1,0){.065234}}
\multiput(41.81,14.4)(.071103,-.047124){16}{\line(1,0){.071103}}
\multiput(42.95,13.64)(.077652,-.047425){15}{\line(1,0){.077652}}
\multiput(44.12,12.93)(.085022,-.047698){14}{\line(1,0){.085022}}
\multiput(45.31,12.26)(.0934,-.047944){13}{\line(1,0){.0934}}
\multiput(46.52,11.64)(.10304,-.04816){12}{\line(1,0){.10304}}
\multiput(47.76,11.06)(.10475,-.04431){12}{\line(1,0){.10475}}
\multiput(49.01,10.53)(.11598,-.04408){11}{\line(1,0){.11598}}
\multiput(50.29,10.04)(.12929,-.04373){10}{\line(1,0){.12929}}
\multiput(51.58,9.61)(.14535,-.04324){9}{\line(1,0){.14535}}
\multiput(52.89,9.22)(.16521,-.04256){8}{\line(1,0){.16521}}
\multiput(54.21,8.88)(.19048,-.04163){7}{\line(1,0){.19048}}
\multiput(55.55,8.59)(.22387,-.04031){6}{\line(1,0){.22387}}
\multiput(56.89,8.34)(.3378,-.048){4}{\line(1,0){.3378}}
\multiput(58.24,8.15)(.4525,-.0473){3}{\line(1,0){.4525}}
\put(59.6,8.01){\line(1,0){1.362}}
\put(60.96,7.92){\line(1,0){1.364}}
\put(62.32,7.88){\line(1,0){1.365}}
\put(63.69,7.89){\line(1,0){1.364}}
\put(65.05,7.95){\line(1,0){1.36}}
\multiput(66.41,8.06)(.3388,.0402){4}{\line(1,0){.3388}}
\multiput(67.77,8.22)(.26969,.04214){5}{\line(1,0){.26969}}
\multiput(69.12,8.43)(.22329,.04341){6}{\line(1,0){.22329}}
\multiput(70.46,8.69)(.18988,.04427){7}{\line(1,0){.18988}}
\multiput(71.78,9)(.1646,.04485){8}{\line(1,0){.1646}}
\multiput(73.1,9.36)(.14474,.04526){9}{\line(1,0){.14474}}
\multiput(74.4,9.77)(.12867,.04552){10}{\line(1,0){.12867}}
\multiput(75.69,10.22)(.11536,.04569){11}{\line(1,0){.11536}}
\multiput(76.96,10.72)(.10412,.04576){12}{\line(1,0){.10412}}
\multiput(78.21,11.27)(.094484,.045771){13}{\line(1,0){.094484}}
\multiput(79.44,11.87)(.086102,.045719){14}{\line(1,0){.086102}}
\multiput(80.64,12.51)(.078728,.045616){15}{\line(1,0){.078728}}
\multiput(81.82,13.19)(.072174,.045467){16}{\line(1,0){.072174}}
\multiput(82.98,13.92)(.070442,.048106){16}{\line(1,0){.070442}}
\multiput(84.11,14.69)(.064578,.047699){17}{\line(1,0){.064578}}
\multiput(85.2,15.5)(.059281,.047275){18}{\line(1,0){.059281}}
\multiput(86.27,16.35)(.054465,.046835){19}{\line(1,0){.054465}}
\multiput(87.31,17.24)(.05006,.046377){20}{\line(1,0){.05006}}
\multiput(88.31,18.17)(.04831,.048198){20}{\line(1,0){.04831}}
\multiput(89.27,19.13)(.046493,.049953){20}{\line(0,1){.049953}}
\multiput(90.2,20.13)(.046961,.054357){19}{\line(0,1){.054357}}
\multiput(91.1,21.16)(.047412,.059172){18}{\line(0,1){.059172}}
\multiput(91.95,22.23)(.047848,.064467){17}{\line(0,1){.064467}}
\multiput(92.76,23.32)(.04543,.066194){17}{\line(0,1){.066194}}
\multiput(93.53,24.45)(.045634,.072069){16}{\line(0,1){.072069}}
\multiput(94.26,25.6)(.045798,.078622){15}{\line(0,1){.078622}}
\multiput(94.95,26.78)(.045918,.085996){14}{\line(0,1){.085996}}
\multiput(95.59,27.99)(.045989,.094378){13}{\line(0,1){.094378}}
\multiput(96.19,29.21)(.046,.10402){12}{\line(0,1){.10402}}
\multiput(96.74,30.46)(.04595,.11525){11}{\line(0,1){.11525}}
\multiput(97.25,31.73)(.04582,.12856){10}{\line(0,1){.12856}}
\multiput(97.71,33.01)(.04559,.14463){9}{\line(0,1){.14463}}
\multiput(98.12,34.32)(.04523,.1645){8}{\line(0,1){.1645}}
\multiput(98.48,35.63)(.0447,.18978){7}{\line(0,1){.18978}}
\multiput(98.79,36.96)(.04393,.22319){6}{\line(0,1){.22319}}
\multiput(99.06,38.3)(.04276,.26959){5}{\line(0,1){.26959}}
\multiput(99.27,39.65)(.0409,.3387){4}{\line(0,1){.3387}}
\put(99.43,41){\line(0,1){1.36}}
\put(99.55,42.36){\line(0,1){2.387}}
\put(62,81.75){\circle*{1.5}} \put(64,8){\circle*{1.5}}
\put(42.75,44){\circle{4.5}} \put(82.25,44){\circle{5.15}}
\multiput(61.75,82)(.0478723,-1.5585106){47}{\line(0,-1){1.5585106}}
\put(65,61.5){\makebox(0,0)[cc]{$y$}}
\put(62.25,86){\makebox(0,0)[cc]{$o_1$}}
\put(64,3.25){\makebox(0,0)[cc]{$o_2$}}
\qbezier(61.5,82)(14.75,35.88)(38,32.25)
\qbezier(38,32.25)(56.13,29.38)(61.75,82)
\put(42.5,61.38){\vector(-1,-1){.1}}\multiput(45.25,64.5)(-.04782609,-.05434783){115}{\line(0,-1){.05434783}}
\put(63.13,37.13){\vector(0,-1){.1}}\multiput(63,43)(.04167,-1.95833){6}{\line(0,-1){1.95833}}
\put(48.75,65.5){\makebox(0,0)[cc]{$z$}}







\put(92.5,66.9){\vector(-1,2){.1}}


\put(95.25,70){\makebox(0,0)[cc]{$x$}}
\put(76.25,23.5){\vector(1,1){.1}}\qbezier(50,24.25)(60.88,9.13)(76.25,23.5)
\put(55.25,22.5){\makebox(0,0)[cc]{$\sigma$}}
\end{picture}
\caption{The paths $x,y,z$ and the diffeomorphism $\sigma$.}
\label{fig2.1}
\end{figure}
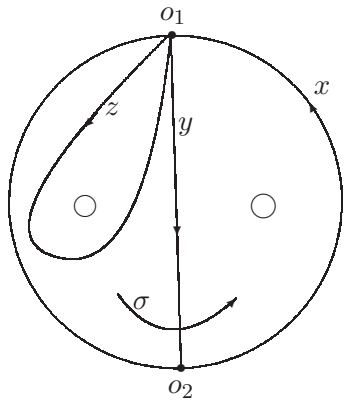

  Up to the action of $\cal D$, there are only two arcs with
boundary on the disk: one inessential isotopic into the boundary,
and one essential separating the two punctures. It is not hard to
see that any essential arc connecting $o_1$ and $o_2$ is the image
of $y$ under some diffeomorphism from ${\cal D}_0$.

One computes that: $\sigma^{2k}(z)= x^{-k} z x^k,
\sigma^{2k}(x)=x$ (which is the well-known action of $B_2$ on the
fundamental group of the punctured disk). There are two boundary
parallel arcs, $zy$ and $xzy$, both of which are fixed by
$\sigma$. We may also compute that $\sigma(y)=zxzy$.

  Now, we compute that $$\sigma^{2k}(y)=\sigma^{2k}(z^{-1}zy)
=\sigma^{2k}(z^{-1}) zy = x^{-k} z^{-1} x^k zy$$ We may also
compute
$$\sigma^{2k+1}(y)=\sigma^{2k}(zxzy) =\sigma^{2k}(z) x zy= x^{-k} z
x^k x z y = x^{-k-1} x z x^{k+1} z y$$ In each case, we see that
the arcs are in the normal form $x^k p_1 x^{-k} p_2$, where each
of $p_1$, $p_2$ is a product of at most 2 factors from $\{x,
y,z^{\pm 1}\}$.

\endproof

\begin{lemma}\label{wu} Let $W$ be a word with $|W|_G>12n_0$. Then
for at most one word $U$ of length $\le 1$, the word $WU$ is
conjugate in $G$ to a word of length $\le |W|_G/15$.
\end{lemma}

\proof By contradiction, assume that there are two distinct words
$U_1$ and $U_2$ of length at most $1$ such that $WU_1$ and $WU_2$
are conjugate in $G$ to some $V_1$ and $V_2$, respectively, with
$|V_1|, |V_2|\le |W|_G/15$. Then we have two annular diagrams
$\Delta_k$ ($k=1,2$) over $G$ with contours $w_ku_k$ and $v_k$
such that  $\Lab (w_k)\equiv W$, $\Lab (u_k)\equiv U_k$, and $\Lab
(v_k)\equiv V_k$.

Since $\Lab (w_1)\equiv \Lab (w_2)$ we can identify the paths
$w_1$ and $w_2$ and obtain a diagram $\Gamma$ on a disk with two
holes. The holes are bounded by paths $v_1$ and $v_2^{-1}$, and
the third boundary component of $\Gamma$ is $x=u_1u_2^{-1}$. There
are two vertices $o_1$ and $o_2$ on the boundary of $\Gamma$
connected by three paths, namely, by $u_1^{-1}$, by $u_2$, and by
$w=w_1=w_2$. Let $\Delta$ be a reduced diagram with the same
boundaries obtained from $\Gamma$.

Note that the words $WU_k$ are not non-trivial in $G$ since
$|W|_G>|U_k|_G$. So are the words $V_1$ and $V_2$. The word
$X\equiv U_1U_2^{-1}$ is also non-trivial since every non-trivial
in $F_2$ word of length at most $2$ is also non-trivial in $G$ by
Lemma 23.16 \cite{book}. Therefore Lemma 22.2 \cite{book} is
applicable to $\Delta$, and there is a simple path $t$ in $\Delta$
connecting some vertices $O_1$ on $x$ and $O_2$ on $v_1$ such that
$|t|<(1/2+4\gamma)(|v_1|+|v_2|+|x|)$. Hence $o_1$ is connected
with $O_2$ by a simple path $s$ of length
$|s|<(1/2+4\gamma)(2|W|_G/15 +2)+1$ because $|V_k|\le |W|_G/15$
and $|U_k|\le 1$. Therefore there is a simple loop $z$ starting at
$o_1$ and surrounding the hole bounded by $v_1$ such that $$|z|\le
2|s|+|v_1|< (1+8\gamma)(2|W|_G/15+2)+2+|W|_G/15 < 0.205|W|_G$$
since $|W|_G>12n_0$ and $\gamma$ is very small.

On the one hand, vertices $o_1$ and $o_2$ are connected by the
simple path $w$ labeled by $W$ in $\Gamma$. The reduction process
(the cancellations of cells) preserves these properties of $w$ in
the following sense (see $\S 13$ \cite{book}): There is a simple
path $w'$ in $\Delta$ connecting $o_1$ and $o_2$ such that its
label $W'$ is equal to $W$ in $G$.

On the other hand, $y=u_1^{-1}$ also connects $o_1$ and $o_2$ in
$\Delta$, and $|y|\le 1$. By Lemma \ref{puti}, the path $w'$ is
homotopic for some $e$, to $p_1x^{e}p_2x^{-e}$ where $|p|_1,
|p_2|<2\times 0.205|W|_G=0.41 |W|_G$.

Hence $|W|_G= |W'|_G \le |p_1|+|p_2|+2|x^e|_G<0.82
|W|_G+2|X^e|_G$. But $|X^e|_G\le n_0$ since we have $X^{n_0}=1$ in
$G$ for every word $X$ of length at most $2$. Therefore we obtain
$0.18|W|_G\le 2n_0$ against the assumption of the lemma. The lemma
is proved.
\endproof

The following lemma gives condition (b) of Proposition \ref{main1}

\begin{lemma}\label{wuton1}
There is a constant $c>1$ such that for every integer $r\ge 1$ and
every word $W$ with $\max (12n_0, cd_r/\phi(r))<|W|_G<d_r
\phi(r)/c$, there exists a word $U$ of length $\le 1$ such that
$(WU)^{n_0}=1$ in $G$.
\end{lemma}

\proof We set $c=150\zeta^{-1}$ where $\zeta$ is the small
positive number from Chapter 7 of the book \cite{book}. Then we
may assume that $W$ is a geodesic word in $G$, that is
$|W|=|W|_G$.

By Lemma \ref{wu}, there exist four different words $U_k$ ($1\le
k\le 4$) of length at most $1$ such that $WU_k$ is not conjugate
in $G$ to a word of length $\le |W|/15$. On the other hand, since
$|WU|< 2|W|$, it follows from the definition of $G(i)$ that $WU_k$
must be conjugate in rank $i=2|W|$ to a power $A_k^{m_k}$ where
$A_k$ is either period of some rank $j_k\le i$ or a simple in rank
$i$ word. Respectively, we have four reduced annular diagrams
$\Delta_k$ of rank $i$ with boundary paths $w_ku_k$ and $p_k$
where $\Lab(w_k)\equiv W$, $\Lab(u_k)\equiv U_k$, and $\Lab
(p_k)\equiv A_k^{m_k}$.

Further one may assume that if $A_k$ is a period of some rank,
then $\Delta_k$ has no cells compatible with $p_k$, since, by
Lemma 13.3 \cite{book}, one may delete such a cell and replace
$A_k^{m_k}$ by $A_k^{m_k \pm n_{A_k}}$. Then we may assume that
$A_k\equiv A_l$ if $A_k$ is conjugate to $A_l^{\pm 1}$ in $G$.
Since $|W|>|U_k|$, we also have $m_k\ne 0$.

{\bf Case 1.} Assume that $m=|m_k| \le 10\zeta^{-1}$ for some $k$.
Then $|A_k|=|A_k^m|/m \ge 0.1\zeta |W|/15=|W|/c>d_r/\phi(r)$. On
the other hand, by Theorem 22.4 and Lemma 26.5 \cite{book},
$|WU_k|\ge \bar\beta m|A_k|$, where $\bar \beta$ is the constant
form Chapter 7 of \cite{book} which is close to $1$. Hence
$|A_k|\le \bar\beta^{-1}(|W|+1)<2|W|=2|W|_G< \phi(r)d_r$.

Thus $|A_k|\in (d_r/\phi(r), d_r\phi(r))$ and $|A_k|<2|W|=i$. It
follows that $A_k$ cannot be simple in rank $i$. (Indeed otherwise
it must be simple in all smaller ranks, and so it is conjugate in
rank $j-1=|A_k|-1$ to a period (or to its inverse) of rank $j$ by
the definition of the set ${\cal X}_j$.) Hence $A_k$ is a period
of rank $j_k=|A_k|$, and $n_A=n_0$ because $|A_k|\in (d_r/\phi(r),
d_r\phi(r))$. By definition of ${\cal R}_{j_k}$, we have
$A_k^{n_0}=1$ in $G$. Since $WU_k$ is conjugate to a power of
$A_k$ in $G$, we also have $(WU_k)^{n_0}=1$, as desired.

{\bf Case 2.} We may now assume that $|m_k|>10\zeta^{-1}$ for
$k=1,\dots,4$. Then for every pair $(k,l)$ ($1\le k< l\le 4$), we
identify diagram $\Delta_k$ with the mirror copy of $\Delta_l$
along the subpaths labeled by $W$ as we did this in the proof of
Lemma \ref{wu}. We obtain diagrams $\Gamma_{kl}$ with 3 boundary
components $x_{kl}$, $v^1_{kl}$, $v^2_{kl}$ labeled by
$X_{kl}\equiv U_kU_l^{-1}$, $A_k^{m_k}$, and $A_l^{-m_l}$,
respectively. Denote by $\Delta_{kl}$ the reduced forms of
diagrams $\Gamma_{kl}$.

Diagrams $\Delta_{kl}$ satisfy all conditions of $E$-maps defined
in subsection 24.2 \cite{book} since $|v^1_{kl}|>
10\zeta^{-1}|A_k|$, $|v^2_{kl}|>10\zeta^{-1}|A_l|$, and
$|x_{kl}|\le 2 < \zeta\min(|v^1_{kl}|,|v^2_{kl}|)$. Lemma 24.6
\cite{book} says that for some $s, s'\in \{1,2\}$, the $E$-map
$\Delta_{kl}$ contains a contiguity submap $\Gamma^{kl}$ of
$v^s_{kl}$ to $v^{s'}_{kl}$ with contiguity degree
$(v^s_{kl},\Gamma^{kl}, v^{s'}_{kl})>0.1$, and also
$(v^{s'}_{kl},\Gamma^{kl}, v^{s}_{kl})>0.1$ if $s\ne s'$.

If $s=s'$ then using lemmas 21.1(1) and 25.8 from \cite{book} for
$\Gamma^{kl}$ we obtain a a contradiction since
$0.1|m_{t}|>\zeta^{-1}$ for $t=k,l$. Then we may assume that $s=1$
and $s'=2$. Lemma 25.10 \cite{book} implies in turn that
$A_k\equiv A_l$ and the cycles $v^1_{kl}$ and $v^2_{kl}$ are
$A_k$-compatible, i.e., the word $X_{kl}=\Lab (x_{kl})$ is
conjugate in $G$ to a power of $A_k$. These non-trivial in $F_2$
words $X_{kl}$ are non-trivial in $G$ by Lemma 23.16 \cite{book}.

Now we have $A_1\equiv\dots\equiv A_4\equiv A$, and all $X_{kl}$
are conjugate to some powers $A^{m_{kl}}$ in $G$. By Lemma 26.5
and Theorem 22.4 \cite{book},
$|m_{kl}||A|<\bar\beta^{-1}|X_{kl}|\le 2\bar\beta^{-1}<3$. It
follows from this estimate and Lemmas 22.1, 23.16 \cite{book} that
all $X_{kl}$ are conjugate to the powers of the same word $A$ in
the free group $F_2$. But this is impossible because it is easy to
see, that for four different words $U_k$ of length at most $1$,
there are two words in the set $\{X_{kl}=U_kU_l^{-1}\}$ which are
conjugate to $(ab)^{\pm 1}$ and $(ab^{-1})^{\pm 1}$, respectively.
Hence Case 2 is also impossible, and the lemma is proved.\endproof

Since every word in ${\cal X}_i$ has length $i$ and the set of
defining relations of $G$ contains $\{A^{n_0} | A\in {\cal
X}_{d_r}\}$, the following lemma gives Condition (c) of
Proposition \ref{main1}.

\begin{lemma}\label{c}
The set ${\cal X}_{d_r}$ is non-empty for every large enough $r$.
If $A\in {\cal X}_{d_r}$, then there exists no non-empty subword
$w$ of the word $A^{n_0}$ such that $|w|\le n_0d_r/2$ and $w$ is
equal in $G$ to a word of length $\le |w|/2$.
\end{lemma}

\proof The number of positive words of length $i$ in $\{ a,b\}$
containing no non-empty subwords of the form $B^6$ is at least
$(3/2)^i$ for $i\ge 1$ (\cite{book}, Theorem 4.6) Since the number
of all words of length $i$ in $\{a,b\}$ is $2^i$ there is such
6-aperiodic word $A$ of arbitrary length $i>>1$ which is not a
proper power. (Indeed, $(3/2)^i>
2^{[i/2]}+2^{[i/3]}+2^{[i/4]}+2^{[i/5]}$ for all large enough
$i$.) Let $i=d_r$. It suffices to prove that $A$ is simple in rank
$d_r-1$.

   Arguing by contradiction, we have from the inductive definition
that $A$ is conjugate in rank $i-1$ to a power of a period $B$ or
some rank $j<i$ or to a power of a simple in rank $i-1$ word $B$
with $|B|<|A|$. In both cases we have a reduced diagram $\Delta$
of rank $i-1$ whose contours $p$ and $q$ are labeled by words $A$
and $B^s$, respectively. In the second case $\Delta$ has no cell
compatible with $q$ (\cite{book}, $\S 13.3$). Hence $q$ is a
smooth section of the B-map $\Delta$ by Lemma 26.5 \cite{book}.

  If $\Delta$ has at least one cell, then there is a cell $\Pi$ in $\Delta$
  with contiguity degree to $p$ greater than $1/2-\alpha-\gamma>\varepsilon$
(Lemma 21.7 and
  Corollary 22.2 of \cite{book}). According to Theorem 22.2 this implies
that
  there is a contiguity subdiagram $\Gamma$ of rank 0 with
$(\Pi,\Gamma,p)\ge\varepsilon$.
  This means than $A$ has a subword $B^t$, where $t=[\varepsilon n_B/2]$.
This contradict
  the 6-aperiodicity of $A$ because $\varepsilon n_B\ge \varepsilon n_0
>12$.

  Thus $\Delta$ has no cells, that is the positive word $A$ is a power of
some word of
  length $|B|<|A|$ in the
  free group $F_2$. This contradicts the choice of $A$, and
  the first statement of the lemma is proved.

  Then assume that the word $w$ is equal to $v$ in $G$. If a reduced diagram
$\Delta$
  for this equality has rank $\ge d_r$, then its perimeter $|w|+|v|$ is at
least
  $(1-\alpha)n_0 d_r>\frac34 n_0d_r$ by Lemma 23.16 \cite{book}. Since
$|w|\le n_0d_r/2$
  we obtain $|v|\ge |w|/2$ as desired.

  If $\rank(\Delta)<d_r$, then the section of the boundary $\partial\Delta$
labeled
  by $w$ is smooth by Lemma 26.5 \cite{book}. It follows from Theorem 22.4
\cite{book}
  that $|v|>(1-\beta)|w|>|w|/2$, and the lemma is proved.

  \endproof

  The group $G$ is infinite by Theorem 26.1 \cite{book}. The order
  of arbitrary element of $G$ divides some $n_A$ by Theorem 26.2
  \cite{book}, and so $G$ is a $p$-group according to our choice of the
  exponents $n_A$. The centralizers of non-trivial elements of $G$
  are cyclic by Theorem 26.5 \cite{book}. Thus we obtain
  condition (d) of Proposition \ref{main1}.

  By Lemma \ref{hyperb}, $\delta_{i_r}\le n_0P(i_r)$. If
$i_{r-1}<i<d_r/\phi(r)$,
  then, for a period $A$ of rank (and length) $i$ , we have by definition of
$n_A$,
   that $n_A|A|\le p\max \{n_0, d_r\}$. If $d_r/\phi(r)\le i\le
i_r=[d_r\phi(r)]$, then
   $n_A|A|=n_0|A|\le n_0d_r\phi(r)$. Hence the obvious induction on $r$
shows that
   $P(i_r)=O(d_r\phi(r))$. Thus $n_0 P(i_r)
  =O(\phi(r)d_r)$, and the Property (e) of Proposition \ref{main1}
  is obtained too.

  It is shown in Theorem 26.2 \cite{book}, that the order of every
  word $X$ does not exceed the  order $n_A$ of a period $A$ of rank (and
  length) $i\le |X|$. As above, we have that if $i_{r-1}<i<d_r/\phi(r)$,
then
   $n_A=O(d_r/i_{r-1})=O(d_r/d_{r-1})= O(\phi(r)^3)$ by the
   definition of $d_r$ and Lemma \ref{hyperb}. If $d_r/\phi(r)\le i\le
i_r=[d_r\phi(r)]$,
   then $n_A=n_0=O(1)$. Since $i\ge i_{r-1}\ge r$, we have by
   induction on $r$ that $n(X) = n_A=o(\phi(|X|)^3)$, and
   property (f) is obtained.
  Since the conditions (a), (b), and (c) are provided
  by lemmas \ref{d1'}, \ref{wuton1}, and \ref{c}, respectively, the proof of
Proposition
  \ref{main1} is complete.


\section{Open problems}
\label{op}

\subsection{Algebraic properties of \ch groups}

Since the class of \ch groups is very large, it would be interesting
to establish more common properties of the groups in this class
except those established in Sections \ref{acochg} and \ref{twor}.

Here is a concrete problem.

\begin{prob} Is it true that the growth of every non-elementary
\ch group is (a) exponential? (b) uniformly exponential?
\end{prob}

Inspired by Theorem \ref{sub} and Corollary \ref{corsub}, it is
natural to ask what kind of subgroups can \ch groups have. In
particular, we formulate the following

\begin{prob}\label{expgr} Can a finitely generated non virtually
cyclic
subgroup of exponential growth of a
\ch group satisfy a non-trivial law?
\end{prob}

\begin{remark} \label{grcond}
The answer to Problem \ref{expgr} is ``no" for \ch groups for which,
using the notation of Remark \ref{dlrem}, the injectivity radii
$r_i$ are ``much larger" than the hyperbolicity constants
$\delta_i$. More precisely, let $G$ be a direct limit of groups
$G_i$ and homomorphisms $\alpha_i\colon G_i\to G_{i+1}$ such that
$G_i=\la S_i\ra$, $\alpha_i(S_i)=S_{i+1}$, $\Gamma(G_i,S_i)$ is
$\delta_i$-hyperbolic, and the induced homomorphism $G_i\to G$ is
injective on a ball of radius $r_i=\exp\exp(C\delta_i)$ for a large
enough constant $C$. We claim that then a subgroup of exponential
growth in $G$ cannot satisfy a non-trivial law.
\end{remark}
\proof Indeed, let $H=\la x_1,...,x_n\ra $ be a subgroup of $G$
having exponential growth. Let $b$ be the maximal length of an
element $x_i$ in generators $S$ of $G$. Let $H_i=\la
x_1(i),...,x_n(i)\ra $ be a pre-image of $H$ in $G_i$. We can assume
that $x_j(i)$ have length at most $b$ in $G_i$.

By \cite[Propositions 3.2 and 5.5]{Kou}, for every $i\ge 1$,  there
exists a pair of elements $u_i, v_i$
of length at most
$\exp\exp(C_0\delta_i)$ (for some uniform constant $C_0$)
in the non virtually cyclic subgroup $H_i$,  generating
a free subgroup of $H_i$. Let $a_i, b_i$ be the images of $u_i, v_i$
in $G$. Let $l_i$ be the length of the shortest word in $\{a_i,
b_i\}^{\pm 1}$ that is equal to 1 in $G$. Then

$$
l_i > \frac{r_i}{\exp\exp(C_0\delta_i)}.
$$

Hence if we assume that $C>C_0$ we deduce that $\lim l_i= \infty$.
Hence $H$ cannot satisfy any non-trivial law.\endproof

\begin{remark} It is easy to see that the \ch groups from
examples in Sections 3-5 can be chosen to satisfy the growth
condition of Remark \ref{grcond}. Thus there are elementary amenable
\ch groups as well as groups with proper subgroups cyclic, torsion
groups, groups with non-trivial centers, etc. satisfying this
condition.
\end{remark}

\begin{prob} It is easy to construct a \ch group with undecidable
word problem (one can use a small cancellation non-recursive
presentation as in Proposition \ref{sparse}). But suppose that the
word problem in a \ch group $G$ is decidable. Does it imply that the
conjugacy problem is decidable as well?
\end{prob}

Using the known facts about solvability of the conjugacy problem in
hyperbolic groups \cite{Gr,Short} it is easy to deduce that the
answer is ``yes" if the growth condition of Remark \ref{grcond}
holds.

It is also interesting to study linearity of \ch groups. We do not
know the answer to the following basic question.

\begin{prob} Is every linear \ch group hyperbolic?
\end{prob}

\subsection{Asymptotic cones and finitely presented groups}

Theorems \ref{classG7} and \ref{cext1} proved in this paper provide
us with a reach source of finitely generated groups all of whose
asymptotic cones are locally isometric, but not all of them are
isometric. Similar methods can be used to show that the groups from
\cite{TV} and from \cite[Section 7]{DS} also satisfy this property.
However all these groups are infinitely presented. Moreover, in all
our examples asymptotic cones are locally isometric to an $\mathbb
R$--tree, which implies hyperbolicity for finitely presented groups
by Proposition \ref{fp}. However the following problem is still
open.

\begin{prob}
Does there exist a finitely presented group all of whose asymptotic
cones are locally isometric, but not all of them are isometric?
\end{prob}

Note that finitely presented groups with different asymptotic cones
were constructed in \cite{OS2005} (earlier, in \cite{KSTT}, such
groups were found under the assumption that the Continuum Hypothesis
does not hold).

\subsection{Asymptotic cones and amenability}

Another interesting problem is to find a characterization of
groups all of whose asymptotic cones are locally isometric to an
$\mathbb R$--tree in the spirit of Theorem \ref{dirlim}. In
particular, do such groups satisfy a suitable small (graded)
cancellation condition? The affirmative answer to this question
and the Kesten-Grigorchuk criterion for amenability would give an
approach to the following.

\begin{prob}\label{prob3}
Suppose that all asymptotic cones of a non-virtually cyclic group
$G$ are locally isometric to an $\mathbb R$--tree. Does it follow
that $G$ is non--amenable?
\end{prob}

Below is another problem about asymptotic cones of amenable groups,
which is still open.

\begin{prob}
Is there a finitely generated (resp. finitely presented) amenable
non--virtually cyclic group all (resp. some) of whose asymptotic
cones have cut--points?
\end{prob}

In particular, we do not know whether our groups from Section
\ref{chag} have cut points in all asymptotic cones (for some
choice of parameters).

\subsection{Divergence and Floyd boundary}

\begin{prob} Is there a finitely presented group with divergence
function $\dv(n,\delta)$ strictly between linear and quadratic for
some $\delta$?
\end{prob}

Recall that if the Floyd boundary $\partial G$ of a finitely
generated group $G$ is nontrivial, $G$ acts on $\partial G$ as a
convergence group \cite{Kar}. On the other hand, geometrically
finite convergence groups acting on non--empty perfect compact
metric spaces are hyperbolic relative to the set of the maximal
parabolic subgroups \cite{Y}.

\begin{prob}\label{prob2}
Suppose that a finitely generated group $G$ has a non--trivial
Floyd boundary. Is $G$ hyperbolic relative to a collection of
proper subgroups?
\end{prob}

Note that if $G$ is hyperbolic relative to a collection of proper
subgroups, then all asymptotic cones of $G$ are tree--graded with
respect to some proper subsets. In particular, all asymptotic cones
of $G$ have cut points. Thus Proposition \ref{Floyd} may be
considered as an evidence towards the positive solution of Problem
\ref{prob2}.

\subsection{Fundamental groups of asymptotic cones}

The example of a group $G$ such that $\pi _1(\CG )=\Z$ for some $d$
and $\omega $ allows us to realize any finitely generated free
Abelian group as the fundamental group of $\CG $ for a suitable $G$
by taking direct products of groups. On the other hand if $1\to N\to
G\to H\to 1$ is a finitely generated central extension and $N$ is
endowed with the metric induced from $G$, then $\CN $ has the
structure of an Abelian topological group. Hence $\pi _1(\CN )$ is
Abelian. Thus there is no hope to construct asymptotic cones with
countable non--Abelian groups by generalizing our methods. This
leads to the following.

\begin{prob}
Does there exist a finitely generated group $G$ such that $\pi_1(\CG
)$ is countable (or, better, finitely generated) and non--Abelian
for some (any) $d$ and $\omega$? Can $\pi_1(\CG)$ be finite and
non-trivial?
\end{prob}

Note that for every countable group $C$ there exists a finitely
generated group $G$ and an asymptotic cone $\CG$ such that
$\pi_1(\CG)$ is isomorphic to the uncountable free power of $C$
\cite[Theorem 7.33]{DS}.


\section*{Appendix: Finitely presented groups whose asymptotic cones are
$\mathbb R$-trees. By M. Kapovich and B. Kleiner.}


\addcontentsline{toc}{section}{Appendix: Finitely presented groups
whose asymptotic cones are $\mathbb R$-trees\\ {\rm By M.~Kapovich
and B.~Kleiner}} \setcounter{section}{8} \setcounter{theorem}{0}

The main result of this appendix is the following

\begin{theorem}\label{main}
Suppose that $G$ is a finitely-presented group such that {\em
some} asymptotic cone of $G$ is an $\R$-tree.  Then $G$ is
Gromov-hyperbolic.
\end{theorem}

This theorem will be an easy application of (a slightly modified
version of) Gromov's local-to-global characterization of
hyperbolic spaces.

Before proving Theorem \ref{main}, we will need several
definitions and auxiliary results.

\medskip
\noindent {\bf 1. Metric notions.} Given a metric space $Z$, let
$B_R(z)$ denote the closed $R$-ball centered at $z$ in $Z$. A
geodesic triangle $\Delta\subset Z$ is called {\em $R$-thin} if
every side of $\Delta$ is contained in the $R$-neighborhood of the
union of two other sides. A geodesic metric space $Z$ is called
{\em $\delta$-Rips-hyperbolic}  if each geodesic triangle in $Z$
is $\delta$-thin. (Rips was the first to introduce this
definition.)

Let $Z$ be a metric space (not necessarily geodesic).
 For a basepoint $p\in Z$ define a number
$\delta_{p}\in [0,\infty]$ as follows. For each $x\in Z$ set
$|x|_p:= d(x,p)$ and
$$
(x,y)_p:= \frac{1}{2}( |x|_p +  |y|_p - d(x,y)).
$$
Then
$$
\delta_p:= \inf_{\delta\in [0,\infty]} \{ \delta | \forall x,y,
z\in Z, (x,y)_p\ge \min( (x,z)_p, (y,z)_p) -\delta\}.
$$
We say that $Z$ is {\em $\delta$-Gromov-hyperbolic}, if
$\infty>\delta\ge \delta_p$ for some $p\in X$. We note that if $Z$
a geodesic metric space which is $\delta$-Gromov-hyperbolic  then
$Z$ is $4\delta$-Rips-hyperbolic and vice-versa (see
\cite[6.3C]{Gr}).

\medskip
A metric space $Z$ is said to have {\em bounded geometry} if there
exists a function $\phi(r)$ such that every $r$-ball in $Z$
contains at most $\phi(r)$ points. For instance, every
finitely-generated group $G$ with a word-metric has bounded
geometry.

\medskip
\noindent {\bf 2. Rips complexes.} Given a metric space $Z$, let
$P_d(Z)$ denote the $d$-Rips complex, i.e., the complex whose
$k$-simplices are $k+1$-tuples of points in $Z$, which are within
distance $\le d$ from each other.
We equip the Rips complex $P_d(Z)$ with a
path metric for which each simplex is path-isometric to a regular
Euclidean simplex of side length $d$.

Given a cell complex $X$, we let $X^i$ denote the $i$-skeleton of
$X$.

\begin{lemma}
\label{sc} Let $G=\la A|\mathcal R\ra $ be a finitely presented
group, $D$ the length of the longest relation in $\mathcal R$.
Then $P_d(G)$ is simply connected for all $d\ge D$.
\end{lemma}
\proof Let $Y$ be the Cayley complex of this presentation, i.e.
the universal cover of the presentation complex of $\la A|R\ra $.
Then $Y^0=G$ and $Y^1$ is the Cayley graph of $G$ (with respect to
the generating set $A$).

First of all, $P_d(G)$ is connected for each $d\ge 1$. We note
that $Y^1=P_1(G)$. Since $\pi_1(Y^1)$ is generated by the
boundaries of the 2-cells in $Y$, it is clear that the map
$$
\pi_1(Y^1)\to \pi_1(P_d(G))
$$
is trivial for $d\ge D$. Vanishing of $\pi_1(P_d(G))$ however is
slightly  less obvious.

Let $d\ge 1$. Consider a loop $\gamma :S^1\to  P_d^{1}(G)$. After
homotoping $\gamma $ if necessary, we may assume that it is a
simplicial map with respect to some triangulation $\calt$ of
$S^1$. Define a map $\gamma _1:S^1\to  Y^{1}$ as follows. For each
vertex $v$ of $\calt$, let $\gamma _1(v)\in G= P_d^0(G)$ be equal
to $\gamma (v)$. For each edge $e=[v_1v_2]$ of $\calt$, let
$\gamma _1\mbox{\Large \(|\)\normalsize}_e$ be a geodesic in
$Y^{1}$ between $\gamma _1(v_1)$ and $\gamma _1(v_2)$. There is a
natural map
$$
Y^{1}\stackrel{i_1}{\to } P_d^1(G)$$
 which takes each $v\in G=Y^0$ to the corresponding vertex of
$P_d^{0}(G)$ and maps each edge of $Y^{1}$ at constant speed to
the corresponding edge of $P_d^{1}(G)$. Let $\gamma _2{:=} i_1\circ
\gamma _1$.

If  $d\geq D$  then $i_1$ can be extended to a map
$$
Y\stackrel{i_2}{\to } P_d^{2}(G).
$$
Since $Y$ is simply-connected, this implies that $\gamma _2$ is
null-homotopic in $P_d^2(G)$.

On the other hand, we claim that $\gamma _2$ is homotopic to $\gamma $ in
$P_d^{2}(G)$. To see this, for each edge $e=[v_1v_2]$ of $\calt$,
let $y_0=\gamma (v),y_1,...,y_m=\gamma (w)$ be the vertices of $Y^{1}$ on
$\gamma _1(e)$  so that $\gamma _2(e)$ is the concatenation of the edges
$$
[y_0 y_1],...,[y_{m-1} y_m]\subset P_d^{2}(G).
$$
Since $\gamma _1(e)$ is a geodesic between $y_0, y_m$ and
$d_{Y^{1}}(y_0, y_m)\le d$, we get:
$$
d_{Y^{1}}(y_0, y_i)\le d, i=1,...,m-1.
$$
Hence each triple of vertices $y_0, y_i, y_m$ spans a 2-simplex
$\Delta_i$ in $P_d^{2}(G)$. Together these simplices define a
homotopy between $\gamma (e)$ and $\gamma _2(e)$ (rel. the end-points).
Thus the loops $\gamma $ and $\gamma _2$ are homotopic. \qed

\medskip
{\bf 4. Coarse Cartan--Hada\-mard theorem.} Our main technical result is the
following coarse Cartan--Hada\-mard theorem for Gromov-hyperbolic spaces:

\begin{theorem}
\label{localglobal} (Cf. \cite{Gr}, \cite[Theorem
8.1.2]{Bowditch-notes}) There are  constants $C_1$, $C_2$, and
$C_3$ with the following property. Let $Z$ be a metric space of
bounded geometry. Assume that for some $\delta$, and $d\geq C_1\delta$,
every ball of radius $C_2d$ in $Z$ is $\delta$-Gromov-hyperbolic, and
$P_d(Z)$ is $1$-connected.  Then $Z$ is $C_3d$-Gromov-hyperbolic.
\end{theorem}

\noindent One can give a direct proof of this theorem modeled on
the proof of the Cartan-Hadamard theorem.  Instead of doing this,
we will use 6.8.M and 6.8.N from \cite{Gr}.
In brief, the idea of the proof is to translate
Gromov's local-to-global result in \cite{Gr}, which is expressed
using isoperimetric
information, into one using $\delta$-hyperbolicity.

Consider the $d$-Rips complex $P_d(Z)$ of $Z$. Given a polygonal loop $c$
$$
c: S^1\to P_d^{1}(Z),
$$
let $L(c)$ denote the length of $c$ and let $A(c)$ be the
least area of a simplicial disk
$$
f: D^2\to P^{2}_d(Z)
$$
so that $f|\partial D^2=c$. If such disk does not exist, we set $A(c)=\infty$.
Note that in order to retain the proper scaling
behavior,  the length and area are computed
here using the metric on $P^1_d(Z)$ and $P^2_d(Z)$ rather than the
combinatorial length and area.

\medskip
Taking $A_0'=500d^2$ in \cite[6.8.M]{Gr} we get:

\begin{proposition}[6.8.M, adapted version]
\label{6.8M} Suppose that $Z$ is a metric space of bounded
geometry, such that for some $d\geq 0$ every simplicial circle
$S'$ in $P^1_d(Z)$ with
$$
500d^2\leq A(S')\leq 64(500d^2)$$
 satisfies
\begin{equation}
\label{locallinear} L(S')\geq d\sqrt{(4000)(64)(500)}
\end{equation}
and $P_d(Z)$ is $1$-connected.  Then $P^1_d(Z)$ is
$(400)\sqrt{500}d$-Rips-hyperbolic (see \cite[6.8.J]{Gr}) and
$Z$ is $(400)\sqrt{500}d$-Gromov-hyperbolic.
\end{proposition}

\noindent Theorem 6.8.N from \cite{Gr} states

\begin{proposition}[6.8.N]
\label{6.8N} If $Z$ is $\delta$-Gromov-hyperbolic and $d\geq 8\delta$,
then every simplicial circle $S'\subset P^1_d(Z)$ satisfies
$L(S')\geq \frac{d}{4\sqrt{3}}A(S')$.
\end{proposition}

\bigskip
\noindent {\em Proof of Theorem \ref{localglobal}.}
Since the statement
of the theorem is scale invariant, after rescaling the metric we may
assume  that

\begin{equation}
\label{d0def} \frac{500d^2}{4\sqrt{3}}\geq
\sqrt{(4000)(64)(500)}.
\end{equation}
Let $C_1{:=} 32$ and $C_2{:=} 64\cdot 500$. Let $S'\subset
P^1_d(Z)$ be a simplicial circle with

\begin{equation}
\label{areapinch} 500d^2\leq A(S') \leq 64(500d^2)
\end{equation}
and let $f:D\to  P^2_d(Z)$ be a least area simplicial $2$-disk
filling $S'$.  By (\ref{areapinch}), there are at most $(64)(500)$
triangles in the triangulated $2$-disk $D$ which are mapped
isomorphically by $f$. Therefore, if we look at
$\operatorname{Im}(f)\subset P^2_d(Z)$, and let $W\subset
\operatorname{Im}(f)$ be the closure of the union of $2$-simplices
contained in $\operatorname{Im}(f)$, then connected components
$W_i$ of $W$ have diameter $\leq (64)(500)d$.  This means that we
can decompose $D$ along disjoint arcs as the union of disks $D_i$,
$1\leq i\leq k+1$ and regions $E_j$, so that each $f(E_j)$ is at
most 1-dimensional and the diameter of each ``minimal 2-disk''
$f(D_i)$ is at most $(64)(500)d$.  This decomposition corresponds
to the Van-Kampen diagram associated with $f$.

By assumption, every ball of radius $C_2d=64\cdot 500d$ is
$\delta$-hyperbolic and
$$
d\geq C_1\delta=32\delta,
$$
so by applying Proposition  \ref{6.8N} to $f(\partial D_1),\ldots f(\partial
D_{k+1})$ and adding up the results, we obtain
$$
L(S')\geq\frac{d}{4\sqrt{3}}A(S')\geq \frac{d}{4\sqrt{3}}500d^2
\geq d\sqrt{(4000)(64)(500)}
$$
where the last inequality comes from (\ref{d0def}). By Proposition
\ref{6.8M} we conclude that $Z$ is $C_3d$-Gromov-hyperbolic where
$C_3{:=} (400)(\sqrt{500})$. \qed

\begin{cor}
\label{localcor} There exist a constant $0<c<\infty$ such that for
each finitely-presented group $G$ there exists a constant $\rho$
(depending on the presentation) with the property:

Suppose that for some $R\ge \rho$, each ball $B_R(y)\subset Y^1$ is
$cR$-hyperbolic, where $Y^1$ is the Cayley graph of
$G$. Then $G$ is Gromov-hyperbolic.
\end{cor}
\proof The complex $P_d^{2}(G)$ is simply-connected for each $d\ge
D$, see Lemma \ref{sc}, where $D$ is the length of the
longest relator in the presentation of $G$.
Let $C_1, C_2$   be the constants from Theorem
\ref{localglobal}, where $Z=G$ with the word metric. Choose $\rho$
so that $\rho/C_2 = D$. Let $c:= \frac{1}{4C_1 C_2}$.

For $R\ge \rho$ set $d:= R/C_2$ and $\delta:= cR$.  Since
$$
d= \frac{R}{C_2}\ge \frac{\rho}{C_2}=D,
$$
the complex $P_d^{2}(G)$ is
1-connected, see Lemma \ref{sc}. We now verify the assumptions of Theorem \ref{localglobal}.

First, by our choice of the constant $c$,
$$
d= \frac{R}{C_2} \ge 4C_1 cR= C_1(4\delta),
$$
in fact, the equality holds.

Next, by the assumption of Corollary \ref{localcor}, each ball
$B_R(y)\subset Y^1$ is $cR=\delta$-Rips-hyperbolic. Therefore for
each $x\in G$, the ball $B_R(x)\subset G$ is
$4\delta$-Gromov-hyperbolic. Since  $C_2d=R$,  every $C_2d$-ball
in $G$ is $4\delta$-Gromov-hyperbolic.

Theorem \ref{localglobal} now implies that $Y^1$ is Gromov-hyperbolic.
\qed

\medskip
{\bf 5. Proof of Theorem  \ref{main}.} Let $\omega$ be a nonprincipal
ultrafilter on $\N$, let $R_j$ be a sequence of positive real
numbers such that $\lio R_j =\infty$. Let $Y$ be the Cayley complex
of a finite presentation of $G$.

Let $y_j\in G$ be a sequence.  By our assumption, the asymptotic
cone $\lio \frac{1}{R_j}(Y^1, y_j)$ is a tree for some choice of
$\omega$ and $(R_j)$. Thus each geodesic triangle in an $R_j$-ball
$B(y_j, R_j)\subset Y^1$ is $\delta_j$-thin, where
$\lio\frac{\delta_j}{R_j}=0.$
 Hence the same is true for each ball
$B(y, R_j)\subset Y^1$, $y\in G$. For sufficiently large $j$,
$R_j\ge \rho=\rho(Y)$ and $\frac{\delta_j}{R_j}<c$, where $\rho, c$
are the constants from the previous corollary. Hence, by Corollary
\ref{localcor}, the graph $Y^1$ is Gromov-hyperbolic and therefore
$G$ is too. \qed

\vspace{5mm}

\noindent Michael Kapovich:\\
{\small\sc Department of Mathematics, University of California, Davis, CA 95616.\\}
{\it E-mail:} {\tt kapovich@math.ucdavis.edu}

\vspace{3mm}
\noindent Bruce Kleiner:\\
{\small\sc Department of Mathematics, Yale University, New Haven, CT 06520-8283.\\}
{\it E-mail:} {\tt bruce.kleiner@yale.edu}

 \vspace{3mm}

\noindent Alexander Yu. Olshanskii:\\
{\small \sc Department of Mathematics, Vanderbilt University , Nashville, TN 37240.\\
Department of Mathematics, Moscow State University, Moscow, 119899, Russia.\\}
{\it E-mail:} {\tt alexander.olshanskiy@vanderbilt.edu}

\vspace{3mm}

\noindent Denis V. Osin: {\small\sc \\ Department of Mathematics, The City College of New York,
New York, NY 10031.\\}
{\it E-mail:} {\tt denis.osin@gmail.com}

\vspace{3mm}

\noindent Mark V. Sapir:
{\small\sc \\ Department of Mathematics, Vanderbilt University, Nashville, TN 37240.\\}
{\it E-mail: } {\tt m.sapir@vanderbilt.edu}

\end{document}